\documentclass[a4paper,12pt]{amsart}
\usepackage[english]{babel}
\usepackage{amsmath, amsfonts, amssymb, amsthm,mathptmx}
\usepackage{color}
\usepackage{graphicx}
\usepackage{multirow}
\allowdisplaybreaks

\usepackage{rotating,pdflscape,color,amssymb,hyperref,subfigure,psfrag,amsmath,eufrak,bbm,epsfig}

\newtheorem{lemme}{Lemma}[section]

\newtheorem{theo}[lemme]{Theorem}

\newtheorem{cor}[lemme]{Corollary}

\newtheorem{propo}[lemme]{Proposition}

\newenvironment{dem*}{\underline{\textbf{Proof }}}{\hfill$\square$ $$$$}

\theoremstyle{definition}

\newtheorem{rem}[lemme]{Remark}
\newtheorem{ex}[lemme]{Example}
\DeclareMathOperator{\im}{\mathfrak{I}}

\DeclareMathOperator{\eloi}{\overset{(d)}{=}}

\DeclareMathOperator{\defeloi}{\overset{(d)}{:=}}
\DeclareMathOperator{\cloi}{\overset{(d)}{\longrightarrow}}
\DeclareMathOperator{\cvp}{\overset{(\mathbb{P})}{\longrightarrow}}
\DeclareMathOperator{\E}{\mathbb{E}}
\DeclareMathOperator{\EE}{\mathcal{E}}
\DeclareMathOperator{\pro}{\mathbb{P}}

\DeclareMathOperator{\op}{op}

\DeclareMathOperator{\diag}{\text{diag}}

\DeclareMathOperator{\M}{\mathcal{M}}
\DeclareMathOperator{\C}{\mathbb{C}}
\DeclareMathOperator{\R}{\mathbb{R}}

\DeclareMathOperator{\NN}{\mathcal{N}_{\mathbb{C}}}

\newcommand{\ind}{\textbf{       }\textbf{     }\textbf{     }}
\newcommand{\ope}{\operatorname}
\newcommand{\tto}{\longrightarrow}
\newcommand{\ol}{\overline}
\newcommand{\wt}{\widetilde}

\newcommand{\saut}{\text{ } \\}
\newcommand{\ep}{\varepsilon}

\newcommand{\wg}{\operatorname{Wg}}
\newcommand{\Moeb}{\operatorname{Moeb}}
\newcommand{\card}{\operatorname{Card}}

\newcommand{\tr}{\operatorname{Tr}}
\newcommand{\bA}{\mathbf{A}}
\newcommand{\bb}{\mathbf{b}}
\newcommand{\bc}{\mathbf{c}}
\newcommand{\bV}{\mathbf{V}}
\newcommand{\bU}{\mathbf{U}}
\newcommand{\bT}{\mathbf{T}}
\newcommand{\bt}{\mathbf{t}}
\newcommand{\bi}{\mathbf{i}}

\newcommand{\bP}{\mathbf{P}}
\newcommand{\bB}{\mathbf{B}}
\newcommand{\bC}{\mathbf{C}}
\newcommand{\bD}{\mathbf{D}}
\newcommand{\bI}{\mathbf{I}}
\newcommand{\bG}{\mathbf{G}}
\newcommand{\bQ}{\mathbf{Q}}
\newcommand{\bR}{\mathbf{R}}

\newcommand{\be}{\mathbf{e}}

\newcommand{\bM}{\mathbf{M}}
\newcommand{\bN}{\mathbf{N}}
\newcommand{\bX}{\mathbf{X}}

\newcommand{\bW}{\mathbf{W}}

\newcommand{\bH}{\mathbf{H}}

\newcommand{\bJ}{\mathbf{J}}

\newcommand{\dI}{\mathrm{I}}
\newcommand{\dV}{\mathrm{V}}

\newcommand{\bpr}{\begin{pr}}
\newcommand{\epr}{\end{pr}}
\newcommand{\trm}{\textrm}
\newcommand{\f}{\frac}
\newcommand{\ff}{\frac{1}}
\newcommand{\one}{\mathbbm{1}}
\newcommand{\cd}{\cdots}
\newcommand{\ds}{\displaystyle}
\newcommand{\lf}{\left}
\newcommand{\ri}{\right}
\newcommand{\st}{such that }

\newcommand{\bPo}{\mathbf{Po}}

\newcommand{\Tr}{\operatorname{Tr}}

\newcommand{\bes}{\begin{equation*}}
\newcommand{\ees}{\end{equation*}}
\newcommand{\beqy}{\begin{eqnarray}}
\newcommand{\eeqy}{\end{eqnarray}}
\newcommand{\beq}{\begin{eqnarray*}}
\newcommand{\eeq}{\end{eqnarray*}}
\newcommand{\bbe}{\begin{equation}}
\newcommand{\ee}{\end{equation}}
\newcommand{\bbm}{\begin{bmatrix}}
\newcommand{\ebm}{\end{bmatrix}}
\newcommand{\bpm}{\begin{pmatrix}}
\newcommand{\epm}{\end{pmatrix}}
\newcommand{\bdet}{\begin{vmatrix}}
\newcommand{\edet}{\end{vmatrix}}
\newcommand{\la}{\label}
\newcommand{\eqre}{\eqref}
\newcommand{\ti}{\times}
\newcommand{\eg}{\ = \ }
\newcommand{\egd}{\ := \ }

\newcommand{\OO}[1]{O \!\left( #1 \right)}
\newcommand{\oo}[1]{o \!\left( #1 \right)}

 \newcommand{\lto}{\longrightarrow}
 \newcommand{\ovl}{\overline}
\newcommand{\alp}{\alpha}
\newcommand{\al}{\alp}
\newcommand{\ka}{\kappa}
\newcommand{\bet}{\beta}
\newcommand{\ga}{\gamma}
\newcommand{\tta}{\theta}
\newcommand{\lam}{\lambda}
\newcommand{\si}{\sigma}
\newcommand{\ov}{\ovl{v}}
\newcommand{\ld}{\ldots}

\newcommand{\lexp}[2]{{\vphantom{#2}}^{#1}#2}
\newcommand{\lbinom}[3]{{\vphantom{#3}}^{\;\;#1}_{#2}#3}

\newcommand{\Ec}[1]{\E \left[ #1 \right]}
\newcommand{\inve}[1]{ \left( #1 \right)^{-1}}

\newcommand{\bgt}{\begin{itemize}}
\newcommand{\ent}{\end{itemize}}
\newcommand{\ste}{\; ;\;}
\newcommand{\ite}{\item}
\newcommand{\ccdots}{\cdots\cdot}
\newcommand{\ninf}{\underset{n\to\infty}{\longrightarrow}}

\newenvironment{pr}{\noindent {\bf Proof. }}{\hfill $\square$\\}

\begin{document}

\title{Outliers in the Single Ring Theorem}
\author{Florent Benaych-Georges and Jean Rochet}



\thanks{FBG and JR: MAP5,
Universit\'e Paris Descartes,
45, rue des Saints-P\`eres
75270 Paris Cedex 06, France. florent.benaych-georges@parisdescartes.fr, jean.rochet@parisdescartes.fr.}

\maketitle

\begin{abstract}This text is about spiked models of non-Hermitian random matrices.   More specifically, we consider matrices of the type $\bA+\bP$, where the rank of $\bP$ stays bounded as the dimension goes to infinity and where the matrix $\bA$ is a non-Hermitian random matrix, satisfying an isotropy hypothesis: its distribution is invariant under the left and right actions of the unitary group. The macroscopic eigenvalue distribution of such matrices is governed by the so called  \emph{Single Ring Theorem}, due to Guionnet, Krishnapur and Zeitouni. We first prove that if $\bP$ has some eigenvalues out of the maximal circle of the single ring, then $\bA+\bP$ has some eigenvalues (called  \emph{outliers}) in the neighborhood of those of $\bP$, which is not the case for the eigenvalues of $\bP$ in the inner cycle of the single ring. Then, we study the fluctuations of the outliers of $\bA$ around the eigenvalues of $\bP$ and prove that they are distributed as the eigenvalues of some finite dimensional random matrices. Such kind of fluctuations had already been shown for Hermitian models. More surprising facts are that outliers can here have very various rates of convergence to their limits (depending on the Jordan Canonical Form of $\bP$) and that some correlations can appear between outliers at a macroscopic distance from each other (a fact already noticed by Knowles and Yin in \cite{KYIOLW} in the Hermitian case, but only for   non Gaussian models, whereas spiked Gaussian  matrices belong to our model and can have such correlated outliers). Our first result generalizes a   result by Tao proved specifically for matrices with i.i.d. entries, whereas the second one (about the fluctuations) is new. 
\end{abstract}

\section{Introduction}
\begin{figure}[ht]
\centering
\includegraphics[scale=.3]{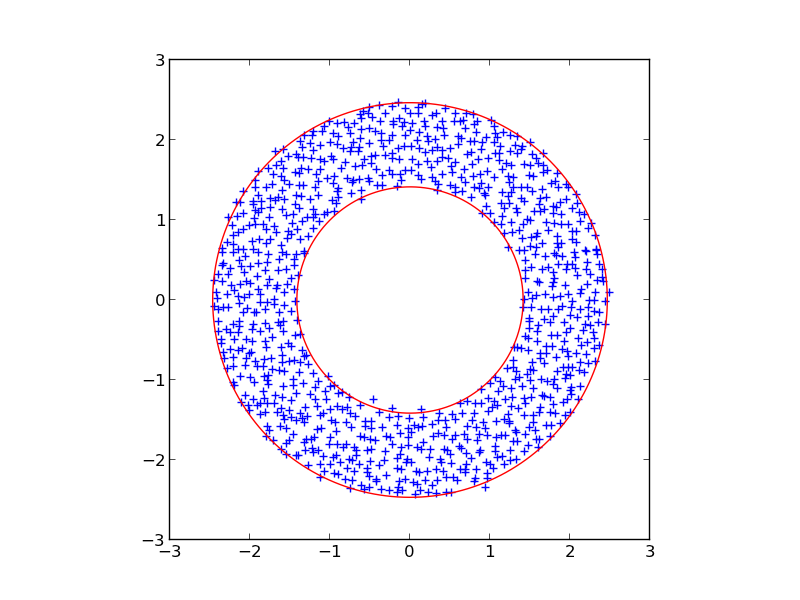}\includegraphics[scale=.3]{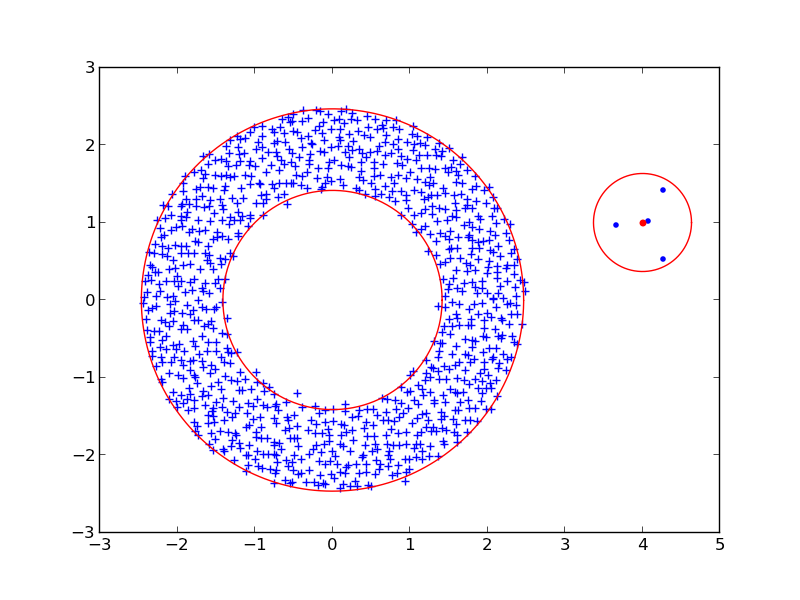}\includegraphics[scale=.3]{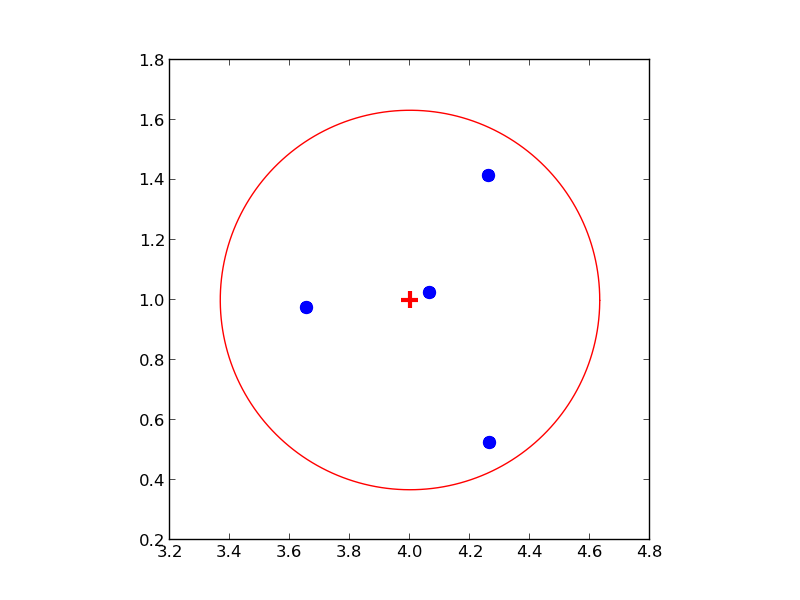}
\caption{Spectrums of $\bA$ (left), of $\bA+\bP$ (center) and zoom on a part of the spectrum of $\bA+\bP$ (right), for the same matrix $\bA$ (chosen as in \eqre{2721410h21} for $s_i$'s uniformly distributed on $[0.5,4]$ with $n=10^3$) and  $\bP$ with rank 4 having one block $\bR_3(\tta)$ and one block $\bR_1(\tta)$ in its Jordan Canonical Form ($\tta=4+i$). We see, on the right, four outliers around $\tta$ ($\tta$ is  the red cross): three of them are at distance $\approx n^{-1/6}$ and one of them, much closer, is at distance $\approx n^{-1/2}$. One can notice that the three ones draw an approximately equilateral triangle. This phenomenon will be explained by   Theorem \ref{t3}.}\label{figure_intro}
\end{figure}

\ind We know that, most times, if one adds to a large random matrix, a finite rank perturbation, it barely modifies its spectrum. However, we observe that the extreme eigenvalues may be altered and deviated away from the bulk. This phenomenon has already been well understood in the Hermitian case. It was shown under several hypotheses in \cite{SP06,FP07, CDF09, CDF09b, BEN1, BEN4, BEN2, BEN3, CDFF, KYIsoSCL, KYIOLW} that for a large random Hermitian matrix, if the strength of the added perturbation is above a certain threshold, then the  extreme eigenvalues of the perturbed matrix deviate at a macroscopic distance from the bulk (such eigenvalues are usually called \emph{outliers}) and have well understood fluctuations, otherwise they stick to the bulk and fluctuate as those of the non-perturbated matrix (this phenomenon is called the \emph{BBP phase transition}, named after the authors of \cite{BAI}, who first brought it to light for empirical covariance matrices). Also, Tao, O'Rourke, Renfrew, Bordenave and Capitaine  studied a non-Hermitian case:  in \cite{TAO1,ROR12,BORCAP1} they considered spiked i.i.d. or elliptic random matrices  and proved that for large enough spikes, some outliers also appear at precise positions. 
In this paper, we study finite rank perturbations for another natural model of non-Hermitian random matrices, namely the \emph{isotropic} random matrices, i.e. the random matrices invariant, in law, under the left and right actions of the unitary group.  Such matrices can be written 
\bbe\la{2721410h21}
\bA = \bU \left(\begin{array}{lll} s_1 & & \\ & \ddots & \\ & & s_n \\
\end{array}\right) \bV, 
\ee 
with $\bU$ and $\bV$ independent Haar-distributed random matrices and  the $s_i$'s some positive numbers which are   independent from $\bU$ and $\bV$. We suppose  that 
the  empirical distribution of   the $s_i$'s tends to a probability measure $\nu$ which is compactly supported on $\R^+$. We know that the singular values of a random matrix with i.i.d. entries satisfy this last condition (where $\nu$ is the Mar$\check{\operatorname{c}}$enko-Pastur quarter circular law with density $\pi^{-1}\sqrt{4-x^2}\mathbf{1}_{[0,2]}(x)dx$, see for example \cite{agz,bai-silver-book,TAO2,BOR1}), so one can see this model as a generalization of the Ginibre matrices (i.e. matrices with i.i.d.  standard complex Gaussian entries). In \cite{GUI}, Guionnet, Krishnapur and Zeitouni showed that the eigenvalues of $\bA$ tend to spread over a single annulus centered in the origin as the dimension tends to infinity. Furthermore in \cite{GUI2}, Guionnet and Zeitouni proved the convergence in probability of the support of its ESD (Empirical Spectral Distribution) which shows the lack of natural outliers for this kind of matrices (see Figure \ref{figure_intro}). This result has been recently improved in \cite{BEN} with exponential bounds for the rate of  convergence.\\ \indent 

In this paper, we prove that, for a finite rank perturbation $\bP$ with bounded operator norm, outliers of $\bA+\bP$ show up close to the   eigenvalues of $\bP$   which are outside the annulus   whereas no outlier appears inside the inner circle of the ring. Then we show (and this is the main difficulty of the paper) that the outliers have fluctuations which are not necessarily Gaussian and whose convergence rates depend on the shape of the perturbation, more precisely on its Jordan Canonical Form\footnote{\label{10101415h44}Recall that any matrix $\bM$ in the set   $ \M_n(\C)$ of $n\ti n$   complex matrices is similar to a   square block diagonal matrix
$$
\left( 
\begin{array}{lllll}
\bR_{p_1}(\tta_1) &  &  & (0)  \\
 & \bR_{p_2}(\tta_2) &  &  \\
 &  & \ddots &  \\
 (0)& &  & \bR_{p_r}(\tta_r) \\
\end{array}
 \right)
\quad \text{ where } \quad
\bR_{p}(\tta) \ = \  \left( 
\begin{array}{lllll}
\tta & 1 &  & (0)  \\
 & \ddots & \ddots &  \\
 &  & \ddots & 1 \\
 (0)& &  & \tta  \\
\end{array}
 \right) \in \M_p(\C), 
$$
which is called the \emph{Jordan Canonical Form} of $\bM$, unique up to the order of the diagonal blocks  \cite[Chapter 3]{HORNJOHNSON}.}.  Let us denote
 by $a<b$ the radiuses of the circles bounding the support of the limit spectral law of $\bA$.
We prove that for any eigenvalue $\tta$ of $\bP$ \st $|\tta|>b$, if one denotes by $$\underbrace{p_{1}, \ld, p_{1}}_{\bet_1\trm{ times}}>\underbrace{p_{2}, \ld, p_{2}}_{\bet_2\trm{ times}}>\cd >\underbrace{p_{\al}, \ld, p_{\al}}_{\bet_\al\trm{ times}}$$ the sizes of the blocks of type $\bR_p(\tta)$ (notation introduced in Footnote \ref{10101415h44}) in the Jordan Canonical Form of $\bP$, then there are exactly $\bet_1p_1+\cd+\bet_{\al}p_\al$ outliers of $\bA+\bP$ tending to $\tta$ and among them, $\bet_1p_1$ go to $\tta$ at rate $n^{-1/(2p_1)}$, $\bet_2p_2$ go to $\tta$ at rate $n^{-1/(2p_2)}$, etc... (see Figure \ref{figure_intro}). Moreover, we give the precise limit distribution of the fluctuations of these outliers around their limits. This limit distribution is not always Gaussian but corresponds to the law of the eigenvalues of some Gaussian matrices (possibly with correlated entries, depending on the eigenvectors of $\bP$ and $\bP^*$). A surprising fact is that some correlations can appear between the fluctuations of outliers with different  limits. In \cite{KYIOLW}, for spiked Wigner matrices, Knowles and Yin had already brought to light some  correlations between outliers at a macroscopic distance from each other   but it was for  non Gaussian models, whereas spiked Ginibre matrices belong to our model and can have such correlated outliers.\\ 

The motivations behind the study of outliers in non-Hermitian models comes mostly from the general effort toward the understanding of the effect of a perturbation with small rank on the spectrum of a large-dimensional operator. The Hermitian case is now quite well understood, and this text provides a review of the question as far as outliers of isotropic non-Hermitian models are concerned. Besides, isotropic non-Hermitian matrix models also appear in wireless networks (see e.g. the recent preprint \cite{zhangqiu}).


\section{Results}
\subsection{Setup and assumptions}
\ind Let, for each $n\geq 1$, $\bA_n$ be a random matrix which admits the   decomposition $\bA_n = \bU_n \bT_n \bV_n$ with $\bT_n = \diag\left(s_1,\ldots,s_n\right)$ where the $s_i$'s are non negative  numbers  (implicitly depending on $n$)  and where $\bU_n$ and $\bV_n$ are two  independent random unitary matrices which are Haar-distributed and independent from the matrix $\bT_n$. We make (part of) the assumptions of the \emph{Single Ring Theorem} \cite{GUI} : \\
~--~ \textbf{Hypothesis 1}: There is a deterministic number $b\ge 0$ \st as $n\to\infty$, we have the convergence in probability $$\ds \ff{n}\Tr (\bT_n^2) \lto  b^2,$$
   
\noindent --~ \textbf{Hypothesis 2}: There exists $M>0$, such that $\pro(\|\bT_n\|_{\op} > M) \tto 0$, \\
~--~ \textbf{Hypothesis 3}: There exist a constant $\kappa>0$ such that
\begin{eqnarray*} 
\im(z) \ > \ n^{-\kappa} \implies \left|\im\left(G_{\mu_{\bT_n}}(z)\right)\right| \ \leq \ \ff{\kappa},
\end{eqnarray*}
where for $\bM$ a matrix, $\mu_\bM$ denotes the empirical spectral distribution (ESD) of $\bM$  and for $\mu$ a probability measure, $G_{\mu}$ denotes the \emph{Stieltjes transform} of $\mu$, that is $\ds G_{\mu}(z) = \int \frac{\mu(dx)}{z-x}$.\\


\begin{ex}
\indent Thanks to \cite{GUI}, we know that our hypotheses are satisfied for example in the model of random complex matrices $\bA_n$ distributed according to the law
$$
\frac{1}{Z_n} \exp \left( -n \tr  V(\bX \bX^*)\right) d \bX,
$$
where $d\bX$ is the Lebesgue measure of the $n \times n$ complex matrices set, $V$ is a polynomial with positive leading coefficient and $Z_n$ is a normalization constant. It is quite a natural unitarily invariant model. One can notice that  $V(x) = \frac{x}{2  }$ gives  the renormalized  \emph{Ginibre matrices}. 
\end{ex}

\begin{rem} \label{1410141} If one strengthens Hypothesis {\bf 1} into the convergence in probability of the ESD $\mu_{\bT_n}$   of $\bT_n$ to a limit probability measure $\nu$, then by the Single Ring Theorem \cite{GUI,RUD},  we know that the ESD $\mu_{\bA_n}$ of $\bA_n$ converges, in probability, weakly to a deterministic probability measure whose support is $\left\{z \in \C, \ a \leq |z| \leq b \right\}$ where
$$
\begin{array}{lll}
a &=& \big(\int x^{-2} \nu (dx)\big)^{-1/2},\\
b &=& \big(\int x^2 \nu (dx)\big)^{1/2}.\\
\end{array}  
$$
\end{rem}

\begin{rem} \label{apos} According to \cite{GUI2}, with a bit more work (this works consists in extracting subsequences within  which the ESD of $\bT_n$ converges, so that we are in the conditions of the previous remark), we know that there is no natural outlier outside the   circle centered at zero with radius $b$ as long as $\| \bT_n \|_{\op}$ is bounded, even if $\bT_n$ has his own outliers. In Theorem \ref{t2}, to make also sure there is no natural outlier inside the inner circle (when $a>0$), we may suppose in addition that $\sup_{n \geq 1} \|\bT_n^{-1}\|_{\op} < \infty$. \end{rem}

\begin{rem}In the case where the matrix $\bA$ is a \emph{real} isotropic matrix (i.e. where $\bU$ and $\bV$ are Haar-distributed on the \emph{orthogonal}  group), despite the facts that the Single Ring Theorem still holds, as proved in \cite{GUI}, and that the Weingarten calculus works quite similarly, our proof does not work anymore: the reason is that we use in a crucial way the bound of Lemma \ref{lemmeFlorent}, proved in \cite{BEN} thanks to an explicit formula for the Weingarten function of the unitary group, which has no analogue for the orthogonal group. However,   numerical simulations  tend to show that similar behaviors occur, with the difference that the radial invariance of certain limit distributions is replaced by the invariance under the action of some discrete groups, reflecting the transition from the unitary group to the orthogonal one. 
\end{rem}

\subsection{Main results}

\begin{figure}[h!]
\centering
\includegraphics[scale=.5]{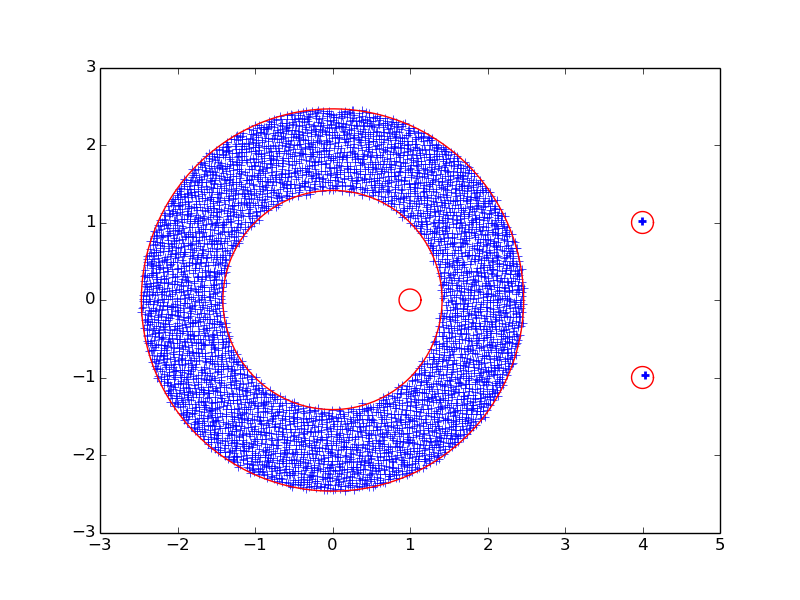}
\caption{Eigenvalues of $\bA_{n}+\bP_n$ for $n=5.10^3$,  $\nu$ the uniform law on $[0.5,4]$ and $\bP_n = \diag(1,4+i,4-i,0,\ldots,0)$. The small circles are centered at $1$,$4+i$ and $4-i$, respectively, and each have a radius $\frac{10}{\sqrt{n}}$ (we will see later that in this particular  case, the rate of convergence of $\lambda_i(\bA_n+\bP_n)$ to $  \lambda_i(\bP_n)$ is $\ff{\sqrt{n}}$).}\label{figure1}
\end{figure}

\ind Let us now consider a sequence of  matrices $\bP_n$ (possibly random, but independent of $\bU_n,\bT_n$ and $\bV_n$) with rank lower than a fixed integer $r$ such that $\|\bP_n\|_{\op}$ is also bounded.  
Then, we have the following theorem (note that in its statement, $r_b$, as the $\lam_i(\bP_n)$'s, can  possibly depend on $n$ and be random):   

\begin{theo}[{Outliers for finite rank perturbation}] \label{t1} Suppose Hypothesis {\bf 1} to hold. Let $\ep >0$ be fixed and suppose that $\bP_n$ hasn't any eigenvalues in the band $\left\{ z \in \C, \ b + \ep < |z| < b + 3 \ep \right\}$ for all sufficiently large $n$, and has $r_b$  eigenvalues counted with multiplicity\footnote{To sort out misunderstandings: we call the \emph{multiplicity} of an eigenvalue its order as a root of the characteristic polynomial, which is greater than or equal to the dimension of the associated eigenspace.}  $\lambda_1(\bP_n), \ldots, \lambda_{r_b}(\bP_n)$ with modulus higher than $b + 3 \ep$. \\
\indent Then, with a probability tending to one, $\bA_n+\bP_n$ has exactly $r_b$ eigenvalues with modulus higher than $b + 2 \ep$. Furthermore, after labeling properly, 
$$
\forall i \in \{1,\ldots,r_b \}, \ \ \lambda_i(\bA_n+\bP_n) - \lambda_i(\bP_n) \ \cvp \ 0 .
$$
\end{theo}
\indent This first   result is a generalization  of  Theorem 1.4 of  Tao's paper \cite{TAO1}, and so is its proof. However, things are different inside the annulus. Indeed, the following result establishes the lack of small outliers: 

\begin{theo}[No outlier inside the bulk] \label{t2} 
Suppose that there exists $M'>0$ \st   $$\pro(\|\bT_n^{-1}\|_{\op} > M') \tto 0$$ and that there is $a>0$ deterministic \st we have the convergence in probability $$\ff{n}\sum_{i=1}^n s_i^{-2}\lto \ff{a^2}.$$ Then for all $\delta \in ]0,a[$, with a probability tending to one, 
$$
\mu_{\bA_n+\bP_n} \left( \left\{ z \in \C, |z| < a - \delta \right\} \right) \ = \ 0,
$$
where $\mu_{\bA_n+\bP_n}$ is the Empirical Spectral Distribution of $\bA_n+\bP_n$.
\end{theo}
\indent Theorems \ref{t1} and \ref{t2} are illustrated in Figure \ref{figure1} (see also Figure \ref{figure_intro}). We drew circles around each eigenvalues of $\bP_n$ and we do observe the lack of outliers inside the annulus.

\indent Let us now consider the fluctuations of the outliers. 
We need to be more precise about the perturbation matrix $\bP_n$. Unlike Hermitian matrices, non-Hermitian matrices are not determined, up to a conjugation by a unitary matrix, only by their spectrums. A key parameter here will be the  Jordan Canonical Form (JCF) of $\bP_n$. 
From  now on, we consider a deterministic perturbation $\bP_n$ 
of rank $\le r/2$ with $r$ an integer independent of $n$ (denoting the upper bound on the rank of $\bP_n $ by $r/2$ instead of $r$ will lighten the notations in the sequel).

As $\ope{dim}(\ope{Im}\bP_n+(\ker\bP_n)^\perp)\le r$,   one can find a unitary matrix $\bW_n$ and an $r\ti r$ matrix $\bPo$ \st \bbe\la{192140}\bP_n\ =\ \bW_n \bpm\bPo  &0\\0&0\epm \bW_n^*. \ee To simplify the problem, we shall suppose that $\bPo$ does not depend on $n$ (even though most of what follows can be extended to the case where $\bPo$ depends on $n$ but converges to a fixed $r\ti r$ matrix as $n\to\infty$).

Let us now introduce the \emph{Jordan Canonical Form} (JCF)
of $\bPo $ : we know that up to a basis change, one can write $\bPo$ as a direct sum of   \emph{Jordan blocks}, i.e. blocks of the type \bbe\la{2811410h21}
\bR_{p}(\tta) \ = \  \left( 
\begin{array}{lllll}
\tta & 1 &  & (0)  \\
 & \ddots & \ddots &  \\
 &  & \ddots & 1 \\
 (0)& &  & \tta  \\
\end{array}
 \right) \in \M_{p}(\C)\qquad \trm{($\tta\in \C$, $p\ge 1$)} .
\ee 
Let us denote by $\tta_1, \ld, \tta_q$ the distinct eigenvalues of $\bPo $ which are in $\{|z|>b+ 3\ep\}$ (for $b$ as in Hypothesis {\bf 1} and  $\ep$ as in the hypothesis of Theorem \ref{t1})  and for each $i=1, \ld, q$, introduce a positive integer $\al_i$, some positive integers      $p_{i,1}> \cdots> p_{i,\al_i}$ corresponding to the distinct sizes of the blocks relative to the eigenvalue $\tta_i$  and $\bet_{i,1}, \ld, \bet_{i, \al_i}$ \st  for all $j$, $\bR_{p_{i,j}}(\tta_i)$ appears $\bet_{i,j}$ times, so that, for a certain $\bQ\in \ope{GL}_r(\C)$, we have:\\ 
\begin{equation}\label{EqDecP27114} \bJ=\bQ^{-1}\bPo\bQ= \lf(\trm{Matrix with spec.  $\subset\{|z|\le b+3\ep\}$}\ri) \ \bigoplus \ \bigoplus_{i=1}^q \ \bigoplus_{j=1}^{\alp_i} \! \underbrace{\begin{pmatrix}\bR_{p_{i,j}}(\tta_i)&&\\ &\ddots&\\ &&\bR_{p_{i,j}}(\tta_i)\end{pmatrix}}_{
 \bet_{{i,j}} \trm{ blocks}}\end{equation}
where $\oplus$ is defined, for square block matrices, by $\mathbf{M}\oplus \mathbf{N}:=\bpm \mathbf{M}& 0\\ 0&\mathbf{N}\epm$.

\indent The asymptotic orders of the fluctuations of the eigenvalues of $\wt{\bA}_n : = \bA_n  + \bP_n $ depend  on the sizes $p_{i,j}$ of the blocks. Actually, for each $\tta_i$,
 we know, by Theorem \ref{t1},  there are $\sum_{j=1}^{\alp_i}p_{ij}\ti\bet_{i,j}$ eigenvalues of $\wt{\bA}_n$  which tend  to $\tta_i$ : we shall write them with a tilda and a $\tta_i$ 
 on the top left corner: $\lexp{\tta_i}{\wt{\lambda}}$. Theorem \ref{t3} below will state that  for each block with size 
 $p_{i,j}$  corresponding to $\tta_i$  of the JCF of $\bPo$, there are $p_{i,j}$ eigenvalues (we shall write them with $p_{i,j}$ on the bottom left corner : $\lbinom{\tta_i}
 {p_{i,j}}{\wt{\lambda}}$) whose  convergence rate will be $n^{-1/(2p_{i,j})}$. As there are $\bet_{{i,j}}$ blocks of size 
 $p_{i,j}$, there are actually $p_{i,j}\times \bet_{{i,j}}$ eigenvalues tending to $\tta_i$ with  convergence rate    
 $n^{-1/(2p_{i,j})}$ (we shall write them $\lbinom{\tta_i}{p_{i,j}}{\wt{\lambda}_{s,t}}$ with $s \in \{1,\ldots,p_{i,j}\}$ and 
 $t \in \{1,\ldots,\bet_{{i,j}}\}$). It would be convenient to denote by $\Lambda_{i,j}$ the vector with size $p_{i,j}\times \bet_{{i,j}}$ defined by 
 \beqy \label{defLambda}
 \Lambda_{i,j} \egd  \displaystyle \left(n^{1/(2p_{i,j})} \cdot \Big(\lbinom{\tta_i}{p_{i,j}}{\wt{\lambda}_{s,t}} - \tta_i \Big) \right)_{\substack{1\le s\le p_{i,j}\\ 1\le t\le \bet_{i,j}}}.
\eeqy 

\indent Let us now define the family of random matrices that we shall use to characterize the limit distribution of the $\Lambda_{i,j}$'s. 
For each $i=1, \ldots, q$, let $I(\tta_i)$ (resp. $J(\tta_i)$) denote the set, with cardinality $\sum_{j=1}^{\al_i}\bet_{i,j}$, of indices in $\{1, \ld, r\}$ corresponding to the first (resp. last) columns of the blocks $\bR_{p_{i,j}}(\tta_i)$ ($1\le j\le \al_i$) in \eqre{EqDecP27114}.
\begin{rem} \la{rem163129092014}
 Note that the columns of $\bQ$ (resp. of $(\bQ^{-1})^*$) whose index belongs to $I(\tta_i)$ (resp. $J(\tta_i)$) are   eigenvectors of $\bPo$  (resp. of $\bPo^*$) associated to $\tta_i$ (resp. $\ol{\tta_i}$). Indeed, if $k\in I(\tta_i)$ and 
 $\be_k$ denotes the $k$-th vector of the canonical basis, then  
 $\bJ \be_{k} = \tta_i \be_{k}$, so that $\bPo (\bQ \be_{k}) = \tta_i \bQ \be_{k}$. 
 \end{rem}

Now, let  \bbe\la{2121414h2}\lf({m}^{\tta_i}_{k,\ell}\ri)_{\displaystyle^{i=1, \ldots, q,}_{(k,\ell)\in J(\tta_i)\ti I(\tta_i)}} \ee be the random centered complex Gaussian vector with covariance   \bbe\la{2121414h3}\E\lf( {m}^{\tta_i}_{k,\ell} \; {m}^{\tta_{i'}}_{k',\ell'}\ri)\ =\ 0, \qquad \E\lf( {m}^{\tta_i}_{k,\ell} \; \ovl{{m}^{\tta_{i'}}_{k',\ell'}}\ri)\eg  \f{b^2}{\tta_i\ovl{\tta_{i'}}-b^2} \;\be_{ k}^*\bQ^{-1}(\bQ^{-1})^*\be_{ {k'}}\;\be_{ \ell'}^*\bQ^*\bQ\, \be_{ {\ell}},\ee where    $\be_1, \ld, \be_r$ are  the column vectors of the canonical basis of $\C^r$.
Note that each entry of this vector  has  a rotationally invariant Gaussian distribution on the complex plane.

 For each $i,j$, let $K(i,j)$ (resp. $K(i,j)^-$) be the set, with cardinality $\bet_{i,j}$ (resp. $\sum_{j'=1}^{j-1}\bet_{i,j'}$), of  indices in $J(\tta_i)$    corresponding to a block of the type $\bR_{p_{i,j}}(\tta_i)$ (resp. to a block of the type  $\bR_{p_{i,j'}}(\tta_i)$ for  $j'<j$). In the same way, let $L(i,j)$ (resp. $L(i,j)^-$) be the set, with the same cardinality as $K(i,j)$ (resp. as $K(i,j)^-$), of indices in $I(\tta_i)$ corresponding to a block of the type $\bR_{p_{i,j}}(\tta_i)$ (resp. to a block of the type  $\bR_{p_{i,j'}}(\tta_i)$ for  $j'<j$).  Note that $K(i,j)^-$ and $L(i,j)^-$ are empty if $j=1$.  Let us define the random matrices 
 \beqy \label{section2.3}
 \ope{M}^{\tta_i,\mathrm{I}}_{j}\egd[m^{\tta_i}_{k,\ell}]_{\ds^{k\in K(i,j)^-}_{\ell \in L(i,j)^-}} &\qquad \qquad \qquad& \ope{M}^{\tta_i,\dI\dI}_{j}\egd[m^{\tta_i}_{k,\ell}]_{\ds^{k\in K(i,j)^-}_{\ell\in L(i,j)}} \nonumber\\
 &&\\
 \ope{M}^{\tta_i,\dI\dI\dI}_{j}\egd[m^{\tta_i}_{k,\ell}]_{\ds^{k\in K(i,j)}_{\ell\in L(i,j)^-}} &\qquad \qquad \qquad& \ope{M}^{\tta_i,\dI\dV}_{j}\egd[m^{\tta_i}_{k,\ell}]_{\ds^{k\in K(i,j)}_{\ell \in L(i,j)}} \nonumber
 \eeqy
 and then let us define the  matrix 
${\bM}^{\tta_i}_{j}$   as 
\bbe\la{21021415h}\bM^{\tta_i}_{j}\egd  \tta_i \left(\ope{M}^{\tta_i,\dI\dV}_{j}-\ope{M}^{\tta_i,\dI\dI\dI}_{j}\inve{\ope{M}^{\tta_i,\dI}_{j}}\ope{M}^{\tta_i,\dI\dI}_{j} \right)
\ee
\begin{rem} \label{rem230214}
It follows from the fact that the matrix $\bQ$ is invertible,  that $\ope{M}^{\tta_i,\dI}_{j}$ is a.s.  invertible and so is $\bM^{\tta_i}_j$. 
\end{rem}

\begin{rem}
\indent From the Remark \ref{rem163129092014} and \eqref{2121414h3}, we see that each matrix   $\bM^{\tta_i}_{j}$  essentially depends  on the eigenvectors of $\bP_n$ and of $\bP_n^*$ associated to blocks $\bR_{p_{i,j}}(\tta_i)$ in \eqre{EqDecP27114} and the correlations between several $\bM^{\tta_i}_{j}$'s depend essentially on the scalar products of such vectors. 
 \end{rem}
 \indent Now, we can formulate our main result.

\begin{theo} \label{t3}
\begin{enumerate}
\item As $n$ goes to infinity, the random vector  $$\displaystyle \left(\Lambda_{i,j} \right)_{\displaystyle^{1 \leq i \leq q}_{1 \leq j\leq \alp_i}} $$ defined at \eqre{defLambda} converges jointly   to  the distribution of  a  random vector $$\displaystyle \left(\Lambda^\infty_{i,j} \right)_{\displaystyle^{1 \leq i \leq q}_{1 \leq j\leq \alp_i}} $$ with joint distribution defined by the fact that 
for each $1 \leq i \leq q$ and $1 \leq j \leq \al_i$,  $\Lambda_{i,j}^\infty$ is the collection of the  ${p_{i,j}}^{\trm{th}}$ roots of the  eigenvalues of  $\bM^{\tta_i}_{j}$ defined at \eqre{21021415h}. 
\item The distributions of the  random matrices $\bM^{\tta_i}_{j}$ are absolutely continuous with respect to the Lebesgue measure  and  none of the coordinates of the  random vector $\displaystyle \left(\Lambda^\infty_{i,j} \right)_{\displaystyle^{1 \leq i \leq q}_{1 \leq j\leq \alp_i}} $   has   distribution  supported by a single point. 
\end{enumerate}
\end{theo} 
 
\begin{rem} Each non zero complex number has exactly $p_{i,j}$ ${p_{i,j}}^{\trm{th}}$ roots, drawing a regular $p_{i,j}$-sided polygon. Moreover, by the second part of the theorem, the spectrums of the $\bM^{\tta_i}_{j}$'s almost surely do not contain $0$,  so  each $\Lambda_{i,j}^\infty$ is actually a complex random vector with $p_{i,j}\ti \bet_{i,j}$ coordinates, which draw $\bet_{i,j}$ regular $p_{i,j}$-sided polygons.
%
%
\end{rem}

\begin{ex}\la{252141}
\indent  For example, suppose that   $\bP_n$ has only one   eigenvalue $\tta$ with modulus $>b+2\ep$ (i.e. $q=1$), with multiplicity $4$ (i.e.  $r_b=4$). Then five cases can occur (illustrated by   simulations in  Figure \ref{Fig1321413h58}, see also Figure \ref{figure_intro}, corresponding to the case (b)):\\ 
\begin{enumerate} \ite[(a)] The JCF of $\bP_n$ for $\tta $ has  one block with size $4$  (so that $\al_1=1$, $(p_{1,1},\bet_{1,1})=(4,1)$) : then the $4$ outliers of $\wt{\bA}_n$ are the vertices of a  square   with center $\approx \tta$ and size $\approx n^{-1/8}$ (their limit distribution is the one of the four fourth roots of the complex      Gaussian variable $\tta m^\tta_{1,1}$      with covariance  given by \eqre{2121414h3}).\\
\ite[(b)] The JCF of $\bP_n$ for $\tta$ has  one block with size $3$ and  one block with size $1$ (so that $\al_1=2$, $(p_{1,1},\bet_{1,1})=(3,1)$, $(p_{1,2},\bet_{1,2})=(1,1)$) : then the $4$ outliers of $\wt{\bA}_n$ are the vertices of an equilateral triangle with center $\approx \tta$ and size $\approx n^{-1/6}$ plus a point at distance $\approx n^{-1/2}$ from $\tta$ (the three first ones  behave like the three third roots of the variable $\tta m_{1,1}^\tta$ and the last one behaves like $\tta(m_{4,4}^{\tta} - m_{1,4}^{\tta}m_{4,1}^{\tta}/m_{1,1}^{\tta})$ where $m_{1,1}^{\tta},m_{1,4}^{\tta},m_{4,1}^\tta,m_{4,4}^\tta$ are   Gaussian variables with correlations given by \eqre{2121414h3}). \\
\ite[(c)] The JCF of $\bP_n$ for $\tta$ has two   blocks with size $2$  (so that $\al_1=1$, $(p_{1,1},\bet_{1,1})=(2,2)$) : then the $4$ outliers of $\wt{\bA}_n$ are the extremities of two crossing segments with centers  $\approx \tta$ and size $\approx n^{-1/4}$ (their limit distribution is the one of the   square roots of the eigenvalues of the matrix 
$$
M^{\tta}_1 \ = \ \tta\bpm
m_{1,1}^{\tta } & m_{1,3}^{\tta }\\
m_{3,1}^{\tta } & m_{3,3}^{\tta } \\
\epm
$$
where $m_{1,1}^{\tta},m_{1,3}^{\tta},m_{3,1}^\tta,m_{3,3}^\tta$ are   Gaussian variables with correlations given by \eqre{2121414h3}).\\
\ite[(d)] The JCF of $\bP_n$ for $\tta$ has one   block  with size $2$ and two blocks with size $1$  (so that $\al_1=2$, $(p_{1,1},\bet_{1,1})=(2,1)$, $(p_{1,2},\bet_{1,2})=(1,2)$ ) : then the $4$ outliers of $\wt{\bA}_n$ are  the extremities of a segment  with center    $\approx \tta$ and size $\approx n^{-1/4}$ plus two points at distance    $\approx n^{-1/2}$ from $\tta$ (the two first ones  behave like the square roots of   $\tta m_{1,1}^\tta$ and the two last ones  behave like the eigenvalues of the matrix
$$
M_2^\tta \ = \ \tta \bpm m_{3,3}^\tta & m_{3,4}^\tta \\ m_{4,3}^\tta & m_{4,4}^\tta\epm -  \f{\tta}{m_{1,1}^\tta} \bpm m_{3,1}^\tta \\ m_{4,1}^\tta\epm \bpm m_{1,3}^\tta & m_{1,4}^\tta\epm
$$ 
where the $m_{i,j}^\tta$'s are   Gaussian variables with correlations given by \eqre{2121414h3}).\\
\ite[(e)] The JCF of $\bP_n$ for $\tta$ has four   blocks with size $1$  (so that $\al_1=1$, $(p_{1,1},\bet_{1,1})=(1,4)$) : then the $4$ outliers of $\wt{\bA}_n$ are four points at distance    $\approx n^{-1/2}$ from $\tta$  (their limit distribution is the one of the eigenvalues of the matrix
$$
M_1^\tta \ = \ \tta \bpm m_{1,1}^\tta & m_{1,2}^\tta & m_{1,3}^\tta & m_{1,4}^\tta \\
m_{2,1}^\tta & m_{2,2}^\tta & m_{2,3}^\tta & m_{2,4}^\tta \\
m_{3,1}^\tta & m_{3,2}^\tta & m_{3,3}^\tta & m_{3,4}^\tta \\
m_{4,1}^\tta & m_{4,2}^\tta & m_{4,3}^\tta & m_{4,4}^\tta \\
\epm
$$ 
where the $m_{i,j}^\tta$'s are    Gaussian variables with correlations given by \eqre{2121414h3}).\\
\end{enumerate}
\begin{figure}[ht]
\centering
\subfigure[The blue dots draw a \textbf{square} with center $\approx \tta$ at distance $\approx n^{-1/8}$ from $\tta$]{
\includegraphics[scale=.3]{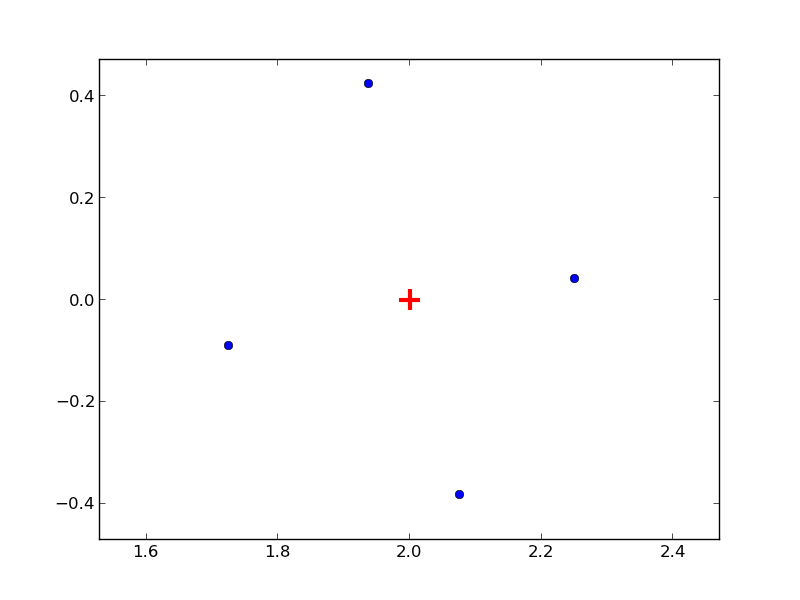} 
\label{Fig2Cas1}}\qquad
\subfigure[The blue dots draw an \textbf{equilateral triangle} with center $\approx \tta$ at distance $\approx n^{-1/6}$  from $\tta$  plus a point at distance $\approx n^{-1/2}$ from $\tta$]{
\includegraphics[scale=.3]{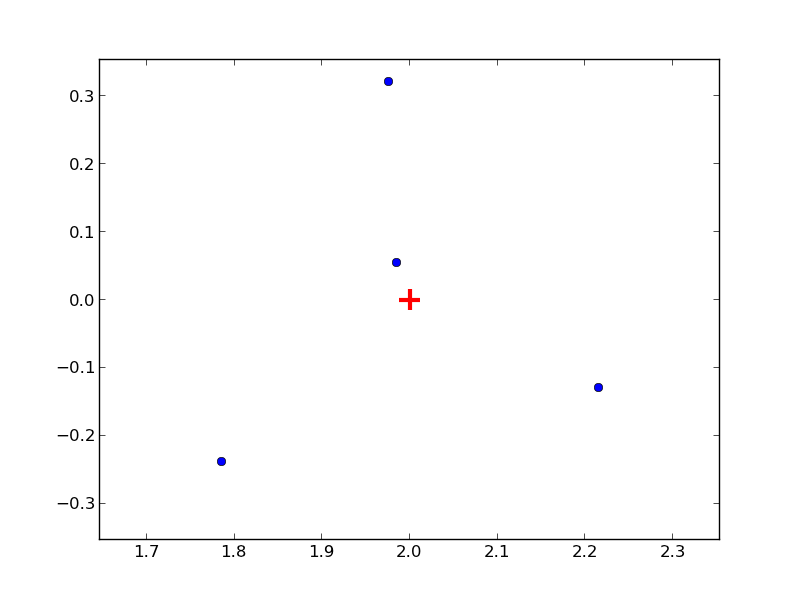} 
\label{Fig2Cas2}}
\subfigure[The blue dots draw two crossing  \textbf{segments} with centers  $\approx \tta$ and lengths $\approx n^{-1/4}$]{
\includegraphics[scale=.3]{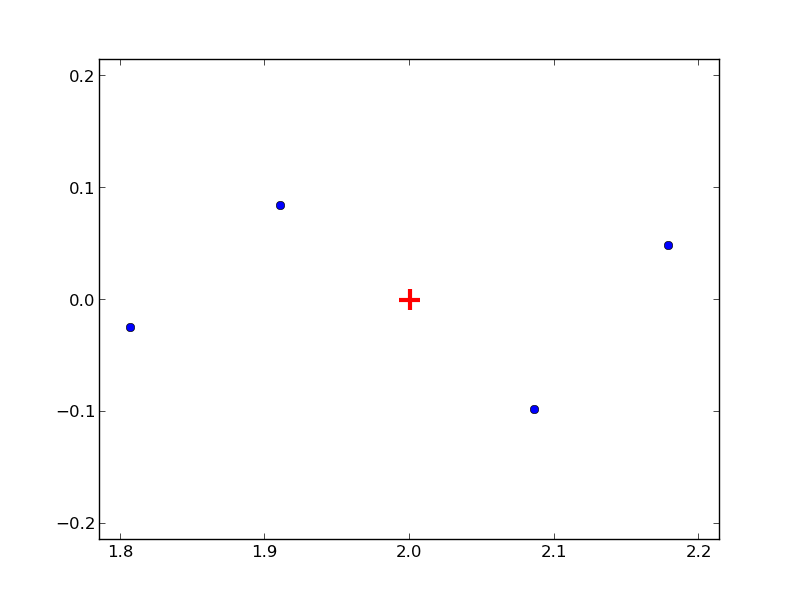} 
\label{Fig2Cas3}}\qquad
\subfigure[The blue dots draw a      \textbf{segment} with center   $\approx \tta$ and length  $\approx n^{-1/4}$  plus two points   at distance $\approx n^{-1/2}$ from $\tta$]{
\includegraphics[scale=.3]{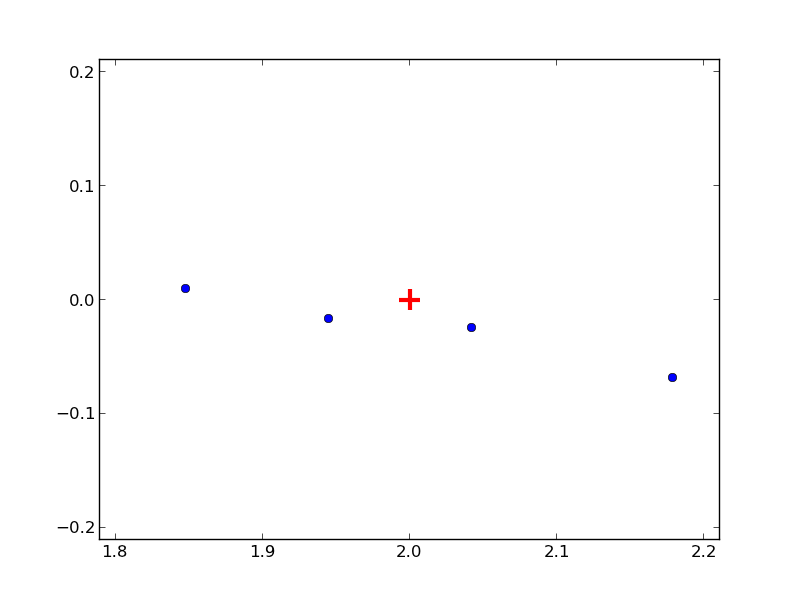} 
\label{Fig2Cas4}}
\caption{The four first cases of Example \ref{252141} (the fifth one, less visual, does not appear here): the red cross is $\tta$ and the blue circular dots are the outliers of $\bA_n+\bP_n$
  tending to $\tta_1$. Each figure is made with the   simulation   of   $\bA_n$  a renormalized Ginibre matrix with size $2.10^3$ plus $\bP_n$    (whose  choice   depends of course of the case) with $\tta=2$. 
}\label{Fig1321413h58} 
\end{figure}
\end{ex}

\subsection{Examples}
\subsubsection{Uncorrelated case}
Let us suppose that\\ \begin{eqnarray}\la{2221420h}&&\forall i,i'=1, \ld,q,\; \forall (k,\ell, k', \ell')\in J(\tta_i)\ti I(\tta_i)\ti J(\tta_{i'})\ti I(\tta_{i'}),\\  \nonumber && \\  \nonumber && \qquad \qquad \qquad \qquad \qquad \qquad \qquad \qquad \be_{ k}^*\bQ^{-1}(\bQ^{-1})^*\be_{ {k'}}\;\cdot\;\be_{ \ell'}^*\bQ^*\bQ\, \be_{ {\ell}}\eg\one_{k=k',\, \ell=\ell'}\\  \nonumber && \end{eqnarray}
Note that it is the case when in  \eqre{2121414h3}, $\bQ$ is unitary,  i.e. when $\bP$ is unitarily conjugated to   $\bpm \bJ&0\\ 0&0\epm$, with $\bJ$ as in \eqre{EqDecP27114}.

By \eqre{2121414h3}, Hypothesis \eqre{2221420h} implies that the entries $m_{k,\ell}^{\tta_i}$ of the random vector of \eqre{2121414h2} are independent and that each $m_{k,\ell}^{\tta_i}$ has a distribution which depends only on $\tta_i$. 
Let us introduce some notation. For $\bet$ a positive integer, we define\footnote{For any $\si>0$,   $\NN(0, \si^2)$ denotes the centered Gaussian law on $\C$ with covariance $\ff{2}\bpm\si^2&0\\0&\si^2\epm$.} \beqy\la{222148h44Flo}&&\ope{Ginibre}(\bet)\egd\trm{$\bet\ti \bet$ random matrix  with i.i.d. $\NN(0,1)$ entries},\\ \nonumber \\ &&
\ope{Ginibre}(\bet, \bet')\egd\trm{$\bet\ti \bet'$ random matrix with i.i.d. $\NN(0,1)$ entries},
\eeqy
 and we get the following corollary:
 \begin{cor}\la{2221421h30}If Hypothesis \eqre{2221420h} holds, then :\begin{enumerate}\ite the collection of random vectors $\displaystyle \left(\Lambda_{i,1},\Lambda_{i,2},\ldots,\Lambda_{i,\al_i} \right)$, indexed by $i=1, \ld, q$, i.e. by the distinct  limit outliers $\tta_i$, is asymptotically independent,\\
 \ite for each $i=1, \ld, q$ and each $j=1, \ld, \al_i$, the matrix $\bM_j^{\tta_i}$ is distributed as:\bgt\ite if $j=1$, then 
$$\bM_j^{\tta_i}\sim \f{\tta_i\, b}{ \sqrt{|\tta_i|^2-b^2}}\ope{Ginibre}(\bet_{i,j}),$$\ite if $j>1$, then  $$\bM_j^{\tta_i}\sim \f{\tta_i\, b}{ \sqrt{|\tta_i|^2-b^2}}\lf(\ope{Ginibre}(\bet_{i,j})-\ope{Ginibre}(\bet_{i,j},\rho_{i,j})\ti\ope{Ginibre}(\rho_{i,j})^{-1}\ti\ope{Ginibre}(\rho_{i,j},\bet_{i,j})\ri),$$\ent
where the four  Ginibre matrices involved if $j>1$ are independent and where $\rho_{i,j}=\sum_{j'=1}^{j-1}\bet_{i,j'}$. 
 \end{enumerate}
 \end{cor}
 
\begin{rem}\la{2221421h31}\text{ }
\bgt\ite The first part of this corollary means that under Hypothesis \eqre{2221420h}, the fluctuations of outliers of $\wt{\bA}_n $  with different limits  are independent. We will see below that it is not always true anymore if Hypothesis \eqre{2221420h} does not hold.
\ite In the second part of this corollary, $j=1$ means that $p_{i,j}=\max_{j'} p_{i,j'}$, i.e. that we consider the outliers of $\wt{\bA}_n $ at the largest possible distance ($\approx n^{-1/(2p_{i,1})}$) from $\tta_i$.
\ite In the second part of the corollary, for $j>1$, the four matrices involved   are independent, but the $\bM_j^{\tta_i}$'s are not independent as $j$ varies (the reason is that the matrix $M_j^{\tta_i,\dI}$ of \eqre{21021415h} contains   $M_{j'}^{\tta_i,\dI\dV}$ as a submatrix as soon as $j'<j$).
\ite If one weakens Hypothesis \eqre{2221420h} by supposing it to hold only for $i=i'$ (resp. $i\ne i'$), then only the second (resp. first) part of the corollary stays true. 
\ent
\end{rem}
 
 The $i=i'$ case of the last point of the previous remark implies the following corollary.
 \begin{cor}\la{2221421h32}
If, for a certain $i$,  $\al_i=\bet_{i,1}=1$  (i.e. if $\tta_i$ is an eigenvalue of $\bP$ with multiplicity\footnote{Let us recall that what is here called the \emph{multiplicity} of an eigenvalue its order as a root of the characteristic polynomial, which is not smaller  than the dimension of the associated eigenspace.} $p_{i,1}$ but with   associated eigenspace having dimension one), then the random vector $$\displaystyle \left(n^{1/(2p_{i,1})} \cdot \left(\lbinom{\tta_i}{p_{i,1}}{\wt{\lambda}_{s,1}} - \tta_i \right) \right)_{1\le s\le p_{i,1}} $$  converges in distribution to the vector of the ${p_{i,1}}^{\trm{th}}$ roots of a $\NN(0, \f{b^2}{|\tta_i|^2(|\tta_i|^2-b^2)})$ random variable. 
 \end{cor}
 
 \subsubsection{Correlated case}
 If Hypothesis \eqre{2221420h} does not hold anymore, then the individual and joint distributions of the random matrices $\bM^{\tta_i}_j$ are not anymore  related to Ginibre matrices  as in Corollary \ref{2221421h30}: the entries of the matrices $M_j^{\tta_i,\dI,\dI\dI,\dI\dI\dI,\dI\dV}$ of \eqre{21021415h} can have non uniform variances, even be correlated, and one can also have  correlations between the entries of two matrices $\bM^{\tta_i}_j$, $\bM^{\tta_{i'}}_{j'}$ for $\tta_i\ne \tta_{i'}$. 
This last case has the surprising consequence that  outliers  of $\wt{\bA}_n $ with different limits can be asymptotically correlated. Such a situation had so far only been brought to light,  by Knowles and Yin in \cite{KYIOLW}, for deformation of non Gaussian Wigner matrices. Note that in our model no restriction on the distributions of the deformed matrix $\bA_n$ is made ($\bA_n$ can for example be a renormalized Ginibre matrix). The following corollary gives an example of a simple situation where such correlations occur. This simple situation corresponds to the following case : we suppose that  for some $i\ne i'$ in $\{1, \ld, q\}$, we have $\bet_{i,1}=\bet_{i',1}=1$. We let   $\ell$ and $\ell'$ (resp. $k$ and $k'$) denote the indices in $\{1, \ld, r\}$ corresponding to the last (resp. first) columns of the block $\bR_{p_{i,1}}(\tta_i)$ and of the block  $\bR_{p_{i',1}}(\tta_{i'})$ and set  \bbe\la{2321413h01}K\egd\be_{ k}^*\bQ^{-1}(\bQ^{-1})^*\be_{ {k'}}\;\cdot\;\be_{ \ell'}^*\bQ^*\bQ\, \be_{ {\ell}} .\ee We will see in the next corollary that as soon as $K\ne 0$,  the fluctuations of outliers at macroscopic distance from each other (i.e. with distinct limits) are not independent.  Set  \bbe\la{2321413h01bis}\si^2\egd \f{|\tta_i|^2b^2}{|\tta_i|^2-b^2}\be_{ k}^*\bQ^{-1}(\bQ^{-1})^*\be_{ {k}}\;\cdot\;\be_{ \ell}^*\bQ^*\bQ\, \be_{ {\ell}}\; ,\qquad {\si'}^2\egd \f{|\tta_{i'}|^2b^2}{|\tta_{i'}|^2-b^2}\be_{ k'}^*\bQ^{-1}(\bQ^{-1})^*\be_{ {k'}}\;\cdot\;\be_{ \ell'}^*\bQ^*\bQ\, \be_{ {\ell'}}.\ee

 \begin{cor}\la{2221421h33}  Under this hypothesis,  for any $1\le s\le p_{i,1}$ and any $1\le s'\le p_{i',1}$, as $n\to\infty$, the random vector \bbe\la{262141}(Z_n, Z_n')\egd\lf(\sqrt{n}\left(\lbinom{\tta_i}{p_{i,1}}{\wt{\lambda}_{s,1}} - \tta_i \right)^{p_{i,1}} , \;\sqrt{n} \left( {\lbinom{\tta_{i'}}{p_{i',1}}{\wt{\lambda}_{s',1}} - \tta_{i'} }\right)^{p_{i',1}}\ri)\ee  converges in distribution to a complex centered Gaussian vector $(Z,Z')$ defined by 
   \bbe\la{2321413h01ter}Z\sim\NN(0, \si^2)\,,\quad  Z'\sim\NN(0, {\si'}^2) \,,\quad\E[ZZ']\eg 0\,,\quad\E[Z\ovl{Z'}]\eg  \f{\tta_i\ovl{\tta_{i'}}\, b^2 \,K}{\tta_i\ovl{\tta_{i'}}-b^2} .\ee
   \end{cor}
   
  \begin{ex}\la{252141PFCor} Let us illustrate this corollary (which is already an example) by a still more particular example. Suppose   that $\bA_n$ is a renormalized  Ginibre matrix and that for $\tta=1.5+i$, $\tta'=3+i$   and for $\ka\in \R\backslash\{-1,1\}$, $\bPo$ is given by $$\bPo\eg
   \bQ\bpm\tta&0\\ 0&\tta'\epm\bQ^{-1}\; ,\qquad  \qquad\bQ\eg \bpm 1&\ka\\ \ka&1\epm\,.$$ In this case, $q=2$, $\al_1=\al_2=p_{1,1}=p_{2,1}=\bet_{1,1}=\bet_{2,1}=1$ and $\ell=k=1$, $\ell'=k'=2$. Thus $\bA_n+\bP_n$ has two outliers $\wt{\lambda}_n\egd\lbinom{\tta}{p_{1,1}}{\wt{\lambda}_{1,1}}$ and $\wt{\lambda'}_n\egd \lbinom{\tta'}{p_{2,1}}{\wt{\lambda}_{1,1}}$ and one can compute  the numbers  $K, \si, \si'$ of \eqre{2321413h01}, \eqre{2321413h01bis} and get \bbe\la{figure_2421415h051}\si^2\eg \f{(1+\ka^2)^2}{(1-|\tta|^{-2})(1-\ka^2)^2} \quad{\si'}^2\eg \f{(1+\ka^2)^2}{(1-|\tta'|^{-2})(1-\ka^2)^2}\quad \E[Z\ovl{Z'}]\eg \f{-4\ka^2 }{(1-(\tta\ovl{\tta'})^{-1})(1-\ka^2)^2}.\ee  
We see that for $\ka=0$,  $Z_n=\sqrt{n}(\wt{\lambda} -\tta)$ and $Z_n'=\sqrt{n}(\wt{\lambda'} -\tta')$ are asymptotically independent, but that for  $\ka\ne 0$, $Z_n$ and $Z_n'$    are not asymptotically independent anymore. This phenomenon and  the accuracy of the approximation $(Z_n,Z_n')\approx (Z,Z')$ for $n\gg 1$  are illustrated by 
Table \ref{label_tableau} and Figure \ref{figure_2421415h05}, where $10^3$ samples of $(Z_n,Z_n')$ have been simulated for $n=10^3$. \begin{table}[ht]
\begin{tabular}{|l|l|l|l|l|l|l|}
\hline 
 & \multicolumn{2}{|c|}{$\E [|Z|^2 ]$} & \multicolumn{2}{|c|}{$\E [|Z'|^2 ]$} & \multicolumn{2}{|c|}{$\E [Z \ol{Z'}  ]$} \\ 
 \cline{2-7}
 & $\kappa=0$  & $\kappa=2^{-1/2}$ & $\kappa=0$  & $\kappa=2^{-1/2}$  &  $\kappa=0$  & $\kappa=2^{-1/2}$    \\
\hline 
\hline
Theorical & $1.444$ & $13.0$ & $1.111$ & $10.0$ & $0.0$  &  $-8.755-1.358i$ \\
\hline
Empirical & $1.492$  & $12.72$ & $1.107$ & $10.04$ & $0.00616-0.00235i$   & $-8.917-1.317i$\\
\hline
\end{tabular}
\caption{Comparison between theoretical asymptotic  formulas \eqre{2321413h01ter} and  \eqre{figure_2421415h051} and a Monte-Carlo numerical computation made out of $10^3$ matrices with size $n=10^3$.}\la{label_tableau}
\end{table}
\begin{figure}[ht]
\centering
\subfigure[$\ka=0$ : uncorrelated case.]{
\includegraphics[scale=.35]{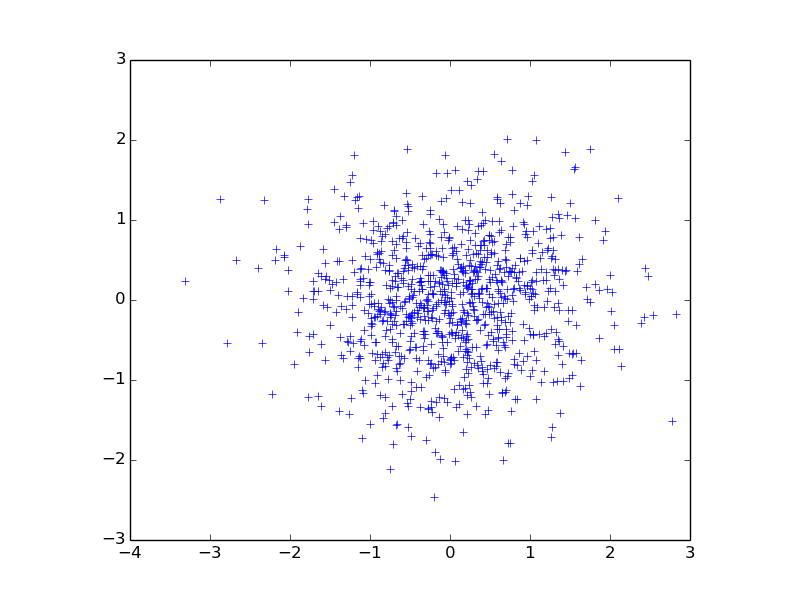} 
\label{FigCorNC1}} 
\subfigure[$\ka=2^{-1/2}$  :  correlated case.  The straight line is the theoretical optimal regression line (i.e. the line with equation $y=ax$ where $a$ minimizes the variance of $Y-aX$, computed thanks to the asymptotic  formulas \eqre{2321413h01ter} and  \eqre{figure_2421415h051}): one can notice that it fits well with the empirical datas.]
{\includegraphics[scale=.35]{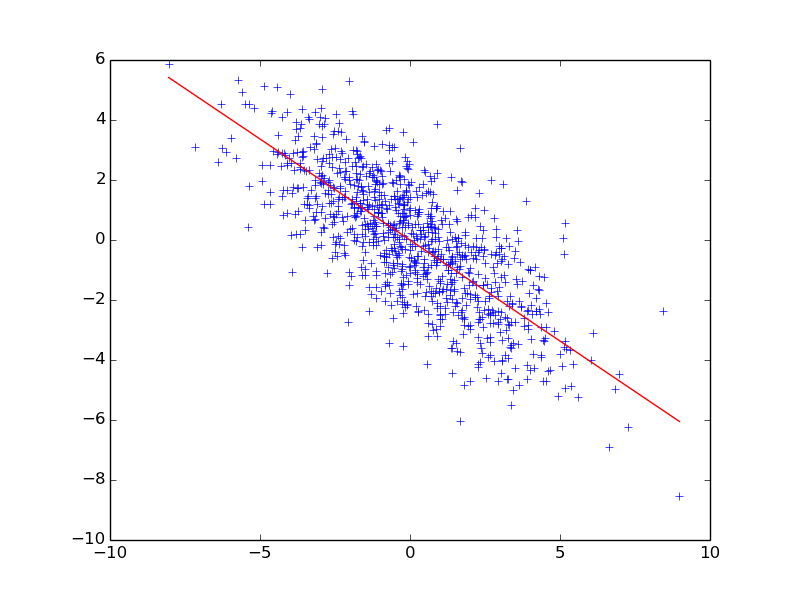} 
\label{Fig2CorC}}
\caption{{\bf Lack of correlation}/{\bf correlation} between outliers with {\bf different limits} :    abscissas (resp. ordinates) of the dots are $X\egd\Im(Z_n)$ (resp. $Y\egd\Im(Z_n')$) for  $10^3$ independent copies of $(Z_n,Z_n')$ (computed thanks to matrices with size $n=10^3$ as for Table  \ref{label_tableau}).}\label{figure_2421415h05}
\end{figure}\end{ex}

\subsection{Preliminaries to the proofs}\la{8813226214} First, for notational brevity, from now on, $n$ will be an implicit parameter ($\bA:=\bA_n$, $\bP:=\bP_n$, \ldots), except in case of ambiguity. 

Secondly,  from now on, we shall suppose that $\bT$ is deterministic. Indeed, once the results established with $\bT$ deterministic,  as $\bT$ is independent from the others   random variables and the only relevant parameter $b$ is deterministic, we can condition on $\bT$ and   apply the deterministic result.   So we suppose that $\bT$ is deterministic and that there is a constant  $M$ independent of $n$ \st for all $n$, $$\|\bT \|_{\op} \le  M.$$ 

 Thirdly,  as the set of probability measures supported by $[0, M]$ is compact, up to an extraction, one can suppose that there is a  probability measure $\Theta$ on $[0,M]$ \st the ESD of $\bT$ converges to $\Theta$ as $n\to\infty$. We will work within this subsequence. This could seem to give a partial convergence result, but in fact, what is proved is that from any subsequence, one can extract a subsequence which converges to the limit given by the theorem. This is of course enough for the proof. Note that by Hypothesis \textbf{1}, we have $\ds b^2=\int x^2  \Theta(dx)$. Having supposed that  the ESD of $\bT$ converges to $\Theta$ insures that  $\bA$ satisfies the hypotheses\footnote{There is actually another assumption in the \emph{Single Ring Theorem} \cite{GUI}, but Rudelson and   Vershynin recently showed in  \cite{RUD} that it was unnecessary.   In \cite{BD13},  Basak  Dembo also weakened the hypotheses  (roughly allowing Hypothesis \textbf{3} not to hold on a small enough set, so that $\nu$ is allowed to have some atoms). 
 As it follows from the recent preprint \cite{BEN} that the 
 convergence of the extreme eigenvalues first established in \cite{GUI2} also works in this case, we could harmlessly weaken our hypotheses down to the ones of  \cite{BD13}.} 
 of the Single Ring Theorem  of \cite{GUI} and of the paper \cite{GUI2}. We will use it once, in the proof of Lemma \ref{lem2}, where we need one of the preliminary results of \cite{GUI2}.

At last, notice that $\bA+\bP$ and $\bV (\bA+\bP) \bV^*$ have the same spectrum,  that  
\begin{equation}\label{88131}
  \ \bV (\bA + \bP) \bV^* \ = \ \bV \bU \bT + \bV \bP  \bV^*,
\end{equation} 
and that as $\bU$ and $\bV$ are independent Haar-distributed matrices,    $\bV \bU$ and $\bV$ are also Haar-distributed and independent. It follows that we shall, instead of the hypotheses made above the statement of Hypotheses \textbf{1}, \textbf{2}, and \textbf{3}, suppose that:
\begin{equation}\label{88132}
 \bA = \bU\bT \qquad\trm{ with \qquad  $\bT$   deterministic and $\bU$ Haar-distributed}\end{equation}\begin{equation*}   \textrm{ and }\qquad \textrm{$\bP$ is independent of $\bA$ and  invariant, in law, by conjugation by any unitary matrix.}  
\end{equation*} 
In the sequel $\E_\bU$  will denote   the expectation with respect to the randomness of $\bU$ and not to the one of $\bP$. 
In the same way, $\E_\bP$ will denote   the expectation with respect to the randomness of $\bP$. 

\subsection{Sketch of the proofs}
We  start  with the following trick, now quite standard in spiked models. Let $\bB \in \M_{n\ti r}(\C)$ and $\bC \in \M_{r\ti n}(\C)$ \st $\bP = \bB \bC$ (where $\M_{p\ti q}(\C)$ denotes the rectangular complex matrices of size $p\ti q$). Then 
\begin{eqnarray}
\det(z\bI - \wt{\bA}) & = & \det(z\bI - (\bA + \bP)) \nonumber\\
                    & = & \det(z\bI - \bA) \det \big( \bI - (z\bI - \bA)^{-1}\bP \big)\nonumber \\
                    & = & \det(z\bI - \bA) \det \big( \bI - (z\bI - \bA)^{-1}\bB\bC \big) \nonumber\\
        \ \ \ & = & \det(z\bI - \bA) \det \big( \bI - \bC (z\bI - \bA)^{-1}\bB \big).\label{1971318:48}
\end{eqnarray} 
 For the last step, we used the fact that for all $\bM \in \M_{r\ti n}$ and $\bN \in \M_{n\ti r}(\C)$, 
$\det \left(\bI_r + \bM \bN\right)  =  \det\left(\bI_n + \bN \bM\right)$. Therefore, the eigenvalues $z$ of $\wt{\bA}$ which are not  eigenvalues of $\bA$ are characterized by 
\bbe\la{2621411h5700}
\det \left( \bI - \bC (z\bI - \bA)^{-1}\bB \right) \ = \ 0. 
\ee
 In view of \eqref{1971318:48}, as previously done by Tao in \cite{TAO1}, we introduce the meromorphic functions (implicitly depending on $n$)
\begin{eqnarray}
f(z) & : = & \det \left( \bI - \bC (z\bI - \bA)^{-1}\bB \right) \ = \ \f{\det(z\bI - \wt{\bA})}{\det(z\bI - \bA)}, \label{def_function_f_21713}\\
g(z) & : = & \det \left( \bI - \bC (z\bI)^{-1} \bB\right) \ = \  \f{\det(z\bI - \bP)}{\det(z\bI)} \label{def_function_g_21713}
\end{eqnarray}
and aim to study the  zeros of $f$.\\

$\bullet$ The proof of Theorem \ref{t1} (eigenvalues outside the outer circle) relies on the fact that on the domain $\{|z|>b+2\ep\}$, $f(z)\approx g(z)$. This follows from the fact that for $|z|>b+2\ep$, the $n\ti n$ matrix $(z\bI-\bA)^{-1}-z^{-1}\bI$ has small entries, and even satisfies \bbe\la{2621411h}\mathbf{x}^*((z\bI-\bA)^{-1}-z^{-1}\bI)\mathbf{y}\ll 1\ee for deterministic unitary column vectors $\mathbf{x}, \mathbf{y}$.\\

$\bullet$ The proof of Theorem \ref{t2} (lack of eigenvalues inside the inner circle) relies on the fact that for $|z|<a-\delta$, $\left\| \bC (z\bI - \bA)^{-1}\bB \right\|_{\op} < 1$. We will see that it  follows from estimates as the one of \eqre{2621411h} for $\bA$ replaced by $\bA^{-1}$.\\

$\bullet$ The most difficult part of the article is the proof of Theorem \ref{t3} about the fluctuations of the outliers around their limits $\tta_i$ ($1\le i\le q$). As the outliers are the zeros of $f$, we shall   expand $f$ around  any fixed  $\tta_i$. Specifically,  for each block size $p_{i,j}$ ($1\le j\le  \al_i$), we prove at Lemma \ref{lemaux0.1} that for $\pi_{i,j} \egd \sum_{l > j} \bet_{{i,l}}p_{i,l}$ and $\bM^{\tta_i}_j$ the matrix with size\footnote{Recall the $\bet_{i,j}$ is the \emph{number} of blocks $\bR_{p_{i,j}}(\tta_i)$ in the JCF of $\bP$.} $\bet_{i,j}$  defined above, we have \bbe\la{2621411h57}f\left( \tta_i + \frac{z}{n^{1/(2p_{i,j})}} \right)\ \approx \ z^{\pi_{i,j}}  \cdot \det\left( z^{p_{i,j}}-\bM^{\tta_i}_j\right).\ee
This proves that $\bA+\bP$ has $\pi_{i,j}$ outliers tending to $\tta_i$ at rate  $\ll n^{-1/(2p_{i,j})}$,  has  $p_{i,j}\ti \bet_{i,j}$ outliers tending to $\tta_i$ at rate  $n^{-1/(2p_{i,j})}$ and that these $p_{i,j}\ti \bet_{i,j}$ outliers are distributed as the ${p_{i,j}}^{\trm{th}}$ roots of the  eigenvalues of $\bM^{\tta_i}_{j}$.
We see that the key result in this proof  is the estimate \eqre{2621411h57}. To prove it, we first specify the choice of the already introduced  matrices $\bB \in \M_{n\ti r}(\C)$ and $\bC \in \M_{r\ti n}(\C)$ \st $\bP = \bB \bC$ by imposing moreover that $\bC\bB=\bJ$ (recall that $\bJ$ is the $r\ti r$ Jordan Canonical Form of $\bP$ of \eqre{EqDecP27114}). Then,  for $$\tilde{z}\egd \tta_i + \frac{z}{n^{1/(2p_{i,j})}}\, ,\qquad \bX_n^{\tilde{z}}\egd \sqrt{n}\bC((\tilde{z}\bI-\bA)^{-1}-\tilde{z}^{-1}\bI)\bB\, ,$$ we write 
\beqy\nonumber
f\left( \tilde{z} \right) & = & \det \left(\bI - \frac{1}{\tilde{z}} \bJ - \ff{\sqrt{n}} \bX_n^{\tilde{z}} \right) \\
 \nonumber & = & \det \left( \bI - \tta_i^{-1}\bJ + \tta_i^{-1}\left( 1 - \frac{1}{1+n^{-1/(2p_{i,j})} z \tta_i^{-1}} \right)\bJ -  \ff{\sqrt{n}} \bX_n^{\tilde{z}} \right)\\
\la{2621412h20}& \approx  & \det \left(\bI - \tta_i^{-1}\bJ +   \f{z   \tta_i^{-2}}{n^{1/(2p_{i,j})}}\bJ -  \ff{\sqrt{n}} \bX_n^{\tilde{z}}  \right)
\eeqy
At this point, one has to note that (obviously) $\det \left(\bI - \tta_i^{-1}\bJ\ri)= 0$ and that (really  not obviously) the $r\ti r$ random array $\bX_n^{\tilde{z}}$ converges in distribution to a Gaussian array as $n\to \infty$ (this is proved thanks to the Weingarten calculus).   Then  the result will follow from a Taylor expansion of \eqre{2621412h20} and a careful look at the main contributions to the determinant.

\section{Eigenvalues outside the outer circle : proof of Theorem \ref{t1}} \label{pt1}
\indent We start with Equations \eqre{1971318:48} and \eqre{2621411h5700}, established in the previous Section, and the functions $f$ and $g$, introduced at \eqre{def_function_f_21713} and \eqre{def_function_g_21713}.

\begin{lemme} \label{p1}
\indent As $n$ goes to infinity, we have
\begin{eqnarray*}
  \sup_{|z| \geq b + 2\ep} \left| f(z) - g(z) \right| \ \cvp \ 0.
\end{eqnarray*}
\end{lemme}

Before proving the lemma, let us explain how it allows to conclude the proof of Theorem   \ref{t1}. The poles of $f$ and $g$ are respectively eigenvalues of the $\bA$ and of the null matrix, hence for $n$ large enough, they have no pole in the region $\left\{z \in \C\ste |z| > b+ 2\ep \right\}$, whereas their zeros in this region are precisely the eigenvalues of respectively $\wt{\bA}$ and $\bP$  that are in this region.
But $|g|$ admits the following lower bound on the circle with radius $b+\ep$  : as we assumed that any eigenvalue of $\bP$ is at least at distance at least $\ep$ from $\left\{ z \in \C \ste \ |z| = b + 2\ep \right\}$, one has
\begin{eqnarray*}
\inf_{|z| = b + 2 \ep }|g(z)| \ = \ \inf_{|z| = b + 2 \ep } \frac{\prod_{i=1}^n \left|z - \lambda_i(\bP) \right|}{|z|^n} 
 & \geq & \left(\frac{\ep}{b + 2 \ep}\right)^r    ,
\end{eqnarray*}
so that by the previous lemma, with probability tending to one, $$\forall z\in \C, |z|=b+2\ep\implies |f(z)-g(z)|<|g(z)|,$$
and so, by   Rouch\'e's Theorem \cite[p. 131]{Beardon}, we know that inside the region $\left\{z \in \C, |z| \leq b+ 2\ep \right\}$, $f$ and $g$ have the   same number of zeros (since they both have $n$ poles). Therefore, as their total number of zeros is $n$, $f$ and $g$ have the same number of zeros outside this region.\\

 Also, Lemma \ref{p1} allows   to conclude   that, after a proper labeling
$$
\forall i \in \{1,\ldots,r_b \}, \ \ \lambda_i(\wt{\bA}) - \lambda_i(\bP) \ \cvp \ 0 .
$$
 Indeed, for each fixed $i\in \{1, \ldots, r_b\}$,
\begin{eqnarray*}
 \prod_{j=1}^{r} \left| 1 - \frac{\lambda_j(\bP)}{\lambda_i(\wt{\bA})}\right|\ = \ \left|g(\lambda_i(\wt{\bA})) \right|
& = &\left| f(\lambda_i(\wt{\bA})) - g(\lambda_i(\wt{\bA}))\right| \\
& \leq & \sup_{|z| \geq b + 2\ep} \left|f(z) - g(z) \right| \ \cvp \ 0.
\end{eqnarray*}  

 Let us now explain how to prove Lemma \ref{p1}. One can notice at first that it suffices to prove that
 \begin{eqnarray}\la{1421315h}
\sup_{|z| \geq b + 2\ep} \left\|\bC (z\bI - \bA)^{-1}\bB - \bC (z\bI)^{-1}\bB  \right\|_{\op} \ \cvp \ 0,
\end{eqnarray}
simply because the function $\det : \M_r(\C) \to \C$ is Lipschitz over every bounded set of $\M_r(\C)$. Then, the proof of Lemma \ref{p1} is based on both following  lemmas (whose proofs are   postponed to Section \ref{tec}).
\begin{lemme} \label{lem1}
\indent There exists a constant $C_1>0$ such that the event $$
\EE_n:=\{\forall k \geq 1, \ \ \|\bA^k\|_{\op} \ \leq \ C_1 \cdot (b+ \ep)^k\}
$$ has  probability tending to one as $n$ tends to infinity.
\end{lemme}

\begin{lemme} \label{lem0}
\indent For all $k \geq 0$, as $n$ goes to infinity, we have 
$$
\| \bC \bA^k \bB \|_{\op} \ \cvp \ 0.
$$
\end{lemme}
On the event $\EE_n$ defined at  Lemma \ref{lem1} above,  we write, for $|z| \ge  b + 2\ep$,
\begin{eqnarray*}
\bC (z\bI - \bA)^{-1}\bB - \bC (z\bI)^{-1}\bB  & = & \bC \sum_{k=1}^{+\infty} \frac{\bA^k}{z^{k+1}} \bB.
\end{eqnarray*}
 and   it suffices to write that  for any $\delta>0$, 
\begin{eqnarray*}
\pro \left( \sup_{|z|\ge b+2\ep}  \left\|\bC (z\bI - \bA)^{-1}\bB - \bC (z\bI)^{-1}\bB \right\|_{\op} > \delta\right) 
& \leq & \pro\left(\EE_n^c\right)+ \pro\left(  \sum_{k=1}^{k_0} \frac{\left\|\bC\bA^k\bB\right\|_{\op}}{(b+2\ep)^{k+1}}   > \f{\delta}{2}\right) \\&& +\pro\left( \EE_n\trm{ and }\left\|\bC \sum_{k=k_0+1}^{+\infty} \frac{\bA^k}{(b+2\ep)^{k+1}} \bB\right\|_{\op} > \f{\delta}{2} \right).\\
\end{eqnarray*}
By to Lemma \ref{lem1} and the fact that $\bC$ and $\bB$ are uniformly bounded (see Remark \ref{165130092014}), we can find   $k_0$ so that the last event has a  vanishing probability. Then, by Lemma \ref{lem0},  the probability of the last-but-one event goes to zero as $n$ tends to infinity.  This gives \eqre{1421315h} and then Lemma \ref{p1}.\\

\section{Lack of eigenvalues inside the inner circle : proof of Theorem \ref{t2}}

\ind  Our goal here is to show that   for all $\delta \in ]0,a[$, with probability tending to one, the function  $f$ defined at \eqref{def_function_f_21713} has no zero  in the region $\left\{z \in \C, |z| < a - \delta \right\}$. Recall that
\begin{eqnarray*}
f(z) &=& \det \left( \bI - \bC (z\bI - \bA)^{-1}\bB \right),
\end{eqnarray*}
so that a simple sufficient condition would be $\left\| \bC (z\bI - \bA)^{-1}\bB \right\|_{\op} < 1$ for 
all $|z| < a - \delta$. Thus, it suffices to prove that with probability tending to one as $n$ tends to infinity, 
$$
  \sup_{|z| < a - \delta} \left\| \bC (z\bI - \bA)^{-1}\bB \right\|_{\op} < 1.
$$
 By Remark \ref{apos}, we know that $\bA$ is invertible. As   in  Section \ref{pt1}, we write,  for all $|z| < a-\delta$,
\begin{eqnarray*}
\bC \left(z\bI - \bA\right)^{-1}\bB & = & -\bC \bA^{-1} \left(\bI - z \bA^{-1} \right)^{-1}\bB \\
& = & - \bC \sum_{k=1}^{\infty} z^{k-1} \bA^{-k} \bB. \\
\end{eqnarray*}
 The idea is to see $\bA^{-1}$ as an isotropic random matrix such as $\bA$, since 
$
\bA^{-1}  =  \bV^* \diag(\frac{1}{s_1},\ldots, \frac{1}{s_n}) \bU^*
$, and satisfies the same kind of hypothesis. Indeed,   Hypotheses \textbf{1} and \textbf{2} are automatiquelly satisfied because $a>0$ (see Remark \ref{apos}), and the following lemma, proved in Section \ref{ProofLemmahypSRT},  insures us that Hypotheses \textbf{3}  is also satisfied.
\begin{lemme} \label{hypSRT}
\indent There exist a constant   $\wt{\kappa} >0$ such that
\begin{eqnarray*}
\im(z) \ > \ n^{-\wt{\kappa}} & \Rightarrow & \left|\im\left(G_{\mu_{\bT^{-1}}}(z)\right)\right| \ \leq \ \ff{\wt{\kappa}}.
\end{eqnarray*}
\end{lemme}
\indent Thus, according to \cite{GUI2}, the support of $\mu_{\bA^{-1}}$ converges in probability to the annulus \\ $\left\{z \in \C, \ b^{-1} \leq |z| \leq a^{-1} \right\}$ as $n\to\infty$, and so, according to \eqref{1421315h}, 
$$
\sup_{|\xi| > a^{-1} + \ep} \bC \sum_{k=1}^{\infty}  \frac{\bA^{-k}}{\xi^{k+1}} \bB \ \cvp \ 0.
$$
\indent Therefore
$$
\pro \left(\sup_{|z| < a - \delta} \left\| \bC (zI - \bA)^{-1}\bB \right\|_{\op} < 1 \right) \ \geq \ 1 - \pro \left( \sup_{|\xi| > a^{-1} + \ep} \left\|\bC \sum_{k=1}^{\infty}  \frac{\bA^{-k}}{\xi^{k-1}} \bB \right\|_{\op}> 1\right) \ \tto \ 1 \ \ ,
$$
with a proper choice for $\ep$. $\square$


\section{Proof of Theorem \ref{t3}}\label{t33}

\subsection{Lemma \ref{lemaux0.1} granted proof   of Theorem \ref{t3}}\la{2621411h50}
Recall that we write $\bP = \bB \bC$ and we know that 
\bbe\label{equation030214}
\sup_{|z|>b+2\ep} \left\| \bC \left(z\bI -\bA  \right)^{-1}\bB - z^{-1}\bC\bB\right\|_{\op}  \ \cvp \ 0,
\ee
(again, for notational brevity,  $n$ will be an implicit parameter, except in case of ambiguity). \\
\indent Following the ideas of \cite{BEN2}, we shall need to differentiate the function $f$ defined at \eqref{def_function_f_21713} to understand the fluctuations of $\wt{\lambda} - \theta$, and to do so, we shall need to be more accurate in the convergence in \eqref{equation030214}. 

Let us first state our key lemma, whose proof is postponed in Section \ref{1421416h59}. Recall from \eqref{EqDecP27114} that we supposed the JCF of $\bP$ to have,  for the eigenvalue $\tta_i$, $\bet_{{i,1}}$  blocks with size $p_{i,1}$, \ld\ld, $\bet_{{i,\alp_i}}$  blocks with size $p_{i,\alp_i}$.  {Recall also that
\begin{eqnarray*}
f(z) &=& \det \left( \bI - \bC (z\bI - \bA)^{-1}\bB \right).
\end{eqnarray*}
}
\begin{lemme} \label{lemaux0.1}
\indent For all $j \in \{1,\ldots,\alp_i\}$, let $F^{\tta_i}_{j}(z)$ be the rational function defined by
 \bbe\la{162142}
F^{\tta_i}_{j}(z) \ : = \ f\left( \tta_i + \frac{z}{n^{1/(2p_{i,j})}} \right).
\ee
\indent Then, there exists a collection  of  positive constants $(\gamma_{i,j})_{\ds^{1\le i\le q}_{1\le j\le \al_i}}$ and a collection of non vanishing random variables $(C_{i,j})_{\ds^{1\le i\le q}_{1\le j\le \al_i}}$ independent of $z$,   such that we have the convergence in distribution (for the topology of the uniform convergence over any compact set)
$$
\lf(  n^{\gamma_{i,j}} F^{\tta_i}_{j}(\cdot) \ri)_{\ds^{1\le i\le q}_{1\le j\le \al_i}}\ \ninf \  \left( z\in \C \ \mapsto \ z^{\pi_{i,j}}\cdot  C_{i,j} \cdot \det\left( z^{p_{i,j}}-\bM^{\tta_i}_j\right) \right)_{\ds^{1\le i\le q}_{1\le j\le \al_i}}
$$
    where $\bM^{\tta_i}_j$ is the random matrix introduced at \eqre{2121414h3} and  $\pi_{i,j} \egd \sum_{l > j} \bet_{{i,l}}p_{i,l}$. 
\end{lemme}
\indent To end the proof of Theorem \ref{t3}, we make sure that we have the right number of eigenvalues of $\wt\bA$ thanks to complex analysis considerations (Cauchy formula) :\bgt\ite Eigenvalues tending to $\tta_i$ with the highest convergence rate : 
\bgt 
\ite[-]  Lemma \ref{lemaux0.1} tells us that on any compact set, $ F^{\tta_i}_{j}$ and  $z^{\pi_{i,j}}\det(z^{p_{{i,j}}}-\bM^{\tta_i}_j)$  have the exact same number of roots (for any large enough $n$, the poles of $F^{\tta_i}_{{j}}$ leave any compact set), so, for the smallest block size    $p_{i,\alp_i}$, we know that $F^{\tta_i}_{{\alp_i}}$ has exactly $\bet_{{i,\alp_i}} \times p_{i,\alp_i}$ roots which do not eventually  leave any  compact set as $n$ goes to infinity.\\ 
\ite[-] Moreover,  we know that the only roots of $F_{\alp_j}^{\tta_i}$ are the $n^{1/(2p_{i,\alp_i})}(\wt{\lambda}-\tta_i)$'s where $\wt\lambda$ are the  eigenvalues of $\wt\bA$.\\
 \ite[-] We conclude that there are exactly $\bet_{{i,\alp_i}} \times p_{i,\alp_i}$ eigenvalues $\left(\lbinom{\tta_i}{p_{i,\alp_i}}{\wt\lambda_{s,t}}\right)_{\displaystyle^{1 \leq s \leq p_{i,\alp_i}}_{1 \leq t \leq \bet_{{i,\alp_i}}}}$ of $\wt\bA$  such that 
$$
n^{1/(2p_{i,\alp_i})}\left(\lbinom{\tta_i}{p_{i,\alp_i}}{\wt\lambda_{s,t}}-\tta_i\right) \ = \ O \left(1 \right),
$$
and thanks to Lemma \ref{lemaux0.1}, we know that the $n^{1/2p_{i,\alp_i}}\left(\lbinom{\tta_i}{p_{i,\alp_i}}{\wt\lambda_{s,t}}-\tta_i\right)$'s satisfy the equation
 $$
\det(z^{p_{{i,\alp_i}}}  - \bM_{\al_i}^{\tta_i}) + o(1) \ = \ 0
$$  
and so are tighted and converge jointly in distribution to the  ${p_{i,j}}^{\trm{th}}$ roots of the eigenvalues of $\bM_{\al_i}^{\tta_i}$.  As $\bM_{\al_i}^{\tta_i}$ is a. s. invertible (recall Remark \ref{rem230214}), none of the $n^{1/2p_{i,\alp_i}}\left(\lbinom{\tta_i}{p_{i,\alp_i}}{\wt\lambda_{s,t}}-\tta_i\right)$'s converge to $0$.\ent
\ite  Then, we take the second smallest size $p_{i,\alp_i-1}$ and work likewise: we know there are exactly $$\pi_{i, \alp_i-1}+\bet_{{i,\alp_i-1}} \times p_{i,\alp_i-1}= \bet_{{i,\alp_i}} \times p_{i,\alp_i} + \bet_{{i,\alp_i-1}} \times p_{i,\alp_i-1}$$ eigenvalues of $\wt\bA$ such that
$$
n^{1/2p_{i,\alp_i-1}} \left( \wt\lambda - \tta_i \right) \ = \ O(1).
$$
 We know that the eigenvalues  $ \lbinom{\tta_i}{p_{i,\alp_i}}{\wt\lambda_{s,t}}$ ($1 \leq s \leq p_{i,\alp_i}$, $1 \leq t \leq \bet_{{i,\alp_i}}$) are among them (because $ p_{i,\alp_i-1}>p_{i, \alp_i}$) so there are $\bet_{{i,\alp_i-1}} \times p_{i,\alp_i-1}$ other eigenvalues  $\left(\lbinom{\tta_i}{p_{i,\alp_i-1}}{\wt\lambda_{s,t}}\right)_{\displaystyle^{1 \leq s \leq p_{i,\alp_i-1}}_{1 \leq t \leq \bet_{{i,\alp_i-1}}}}$ of $\wt\bA$ such that 
$$
n^{1/2p_{i,\alp_i-1}}\left(\lbinom{\tta_i}{p_{i,\alp_i-1}}{\wt\lambda_{s,t}}-\tta_i\right) \ = \ O \left(1 \right).
$$ It follows that $\left(\lbinom{\tta_i}{p_{i,\alp_i-1}}{\wt\lambda_{s,t}}\right)_{\displaystyle^{1 \leq s \leq p_{i,\alp_i-1}}_{1 \leq t \leq \bet_{{i,\alp_i-1}}}}$ 
  converges jointly in distribution to the  $p_{i,\al_{i-1}}^{\trm{th}}$ roots of the eigenvalues of $\bM_{\al_{i-1}}^{\tta_i}$  (which are almost surely non zero). \\
\ite  At each step, $\pi_{p_{i,j}}$ corresponds to the number of eigenvalues we have already ``discovered'' and which go to $\tta_i$ faster than $n^{-1/(2p_{i,j})}$ (because   $p_{i,\alp_i}<\cdots<p_{i,1}$), and so it explains the presence of the factor   $z^{\pi_{{i,j}}}$  before  $\det(z^{p_{i,j}}-\bM^{\tta_i}_j)$ the previous lemma. So one can continue this induction and conclude. that way, we  get the exact number of eigenvalues of $\wt\bA$.
\ent 

It remains now to prove Lemma \ref{lemaux0.1}. We begin with the convergence of $z \mapsto \bX_n^z$.

\subsection{Convergence of $z \mapsto \bX_n^z$.}\label{biggestpart}
 \indent Recall that in order to simplify, we wrote, at \eqre{192140},  
$$
\bP \ = \ \bW \bpm \bPo & 0 \\ 0 & 0 \epm \bW^* \ = \ \bW \bpm \bQ \bJ \bQ^{-1} & 0 \\ 0 & 0 \epm \bW^*,
$$
where $\bJ$ is a Jordan Canonical Form and $\bW$ is supposed to be Haar-distributed from \eqref{88132}. We also wrote $\bP = \bB\bC$ without specifying any choice. For now on, we shall set down
\bbe 
\bB \egd \bW\bpm \bQ\bJ \\ 0 \epm \ \in \M_{n\ti r}(\C) \quad \text{ and }\quad \bC \egd \bpm \bQ^{-1} & 0 \epm \bW^* \ \in \M_{r\ti n}(\C).
\ee
One can easily notice that
\bbe \la{rem240214}
\bC \bB \eg \bJ \quad ; \quad \bB^* \bB \eg \bJ^* \bQ^* \bQ \bJ \quad ; \quad \bC \bC^* \eg \bQ^{-1}(\bQ^{-1})^*,
\ee
so that all these matrix products do not depend on $n$.
 
 \begin{rem}\la{165130092014}
\indent With this specific choice, the norm of the matrix $\bB$ (resp. $\bC$) is uniformly bounded by $\|\bQ\bJ\|_{\op}$ (resp. $\|\bQ^{-1}\|$) which doesn't depend on $n$. 
\end{rem}
\indent For $|z|>b+2\ep$, we define the $\M_r(\C)$-valued random variable \bbe\la{162141}\bX_n^z:=\sqrt{n}\bC \left(\left(z\bI -\bA  \right)^{-1}- z^{-1}\right)\bB.\ee


\begin{lemme}\label{lem040214}
\indent As $n$ goes to infinity, the finite dimensional marginals of $(\bX_n^z)_{|z|>b+2\ep}$ converge  to the ones of a centered complex Gaussian process $(\bX^z=[x_{i,j}^z]_{1\le i,j\le r})_{|z|>b+2\ep}$ such  that for all $\tta, \tta'$ in $\{|z|>b+2\ep\}$, \bgt
\ite $x_{i,j}^{\tta} \ \sim \ \NN \left( 0, \ \frac{b^2}{|\tta|^2} \ff{|\tta|^2 - b^2} \cdot \be_i^* \bC \bC^* \be_i \cdot \be_j^* \bB^* \bB \be_j \right)$,
\ite $ \E \left(x_{i,j}^{\tta}x_{k,l}^{\tta'} \right)  \ = \ 0$,\; $  \E \left(x_{i,j}^{\tta}\ovl{x_{k,l}^{\tta'} }\right)  \ = \ \f{b^2}{\tta\ol{\tta'}} \ff{\tta\ol{\tta'}-b^2} \cdot   \be_i^* \bC \bC^* \be_k \cdot  \be_l^* \bB^* \bB \be_j$.
\ent
\end{lemme}
 

Recall now that the event $\EE_n$ has been defined at  Lemma \ref{lem1} and has probability tending to one. 
 
\begin{lemme}\label{lem040214bis}
 There is $C$ finite \st for $n$ large enough, on  $\{|z| >b+2\ep\}$, $$   \E\lf(\one_{\EE_n}\lf\|\f{\partial}{\partial z} \bX_n^z\ri\|^4\ri)\le C,$$ where $\|\cdot\|$ denotes a norm on $\M_r(\C)$.\end{lemme}

\indent We deduce, by e.g.   \cite[Cor. 14.9]{KallenbergFoundations} (slightly modified because of    the presence of $\one_{\EE_n}$),  that as $n\to\infty$,  the random process $(\bX_n^z)_{|z|>b+2\ep}$ converges weakly, for the topology of uniform convergence on compact subsets,  to the random process $(\bX^z)_{|z|>b+2\ep}$ 

\subsubsection{Proof of lemma \ref{lem040214}} 
\indent  Let us fix an integer $p$, some complex numbers $z_1,\ldots,z_p$ from $\{|z|>b+2\ep\}$, some complex numbers $\nu_1,\ldots,\nu_p$ and some integers $i_1,j_1,\ldots,i_p,j_p$ in $\{1,\ldots,r\}$ and define  
$$
G_n \ := \ \sum_{t=1}^p \nu_t \be_{i_t}^* \bX_n^{z_t} \be_{j_t}.$$ 
At first, we notice that on the event $\EE_n$ of Lemma \ref{lem1}, we can rewrite $G_n$ this way
$$
G_n \ = \ \sqrt{n}\sum_{t=1}^p \nu_t \be_{i_t}^* \bC \sum_{k\geq 1}\frac{\bA^k}{z_t^{k+1}} \bB\be_{j_t} \ = \ \sqrt{n}\sum_{t=1}^p \nu_t \bc_{t}^* \sum_{k\geq 1}\frac{\bA^k}{z_t^{k+1}} \bb_{t}.
$$
where $\bb_t$ designates the $j_t$-th column of $\bB$ and $\bc_t$ the $i_t$-th column of $\bC^*$. As $\pro(\EE_n)\longrightarrow 1$,  $\EE_n^c$ is irrelevant to weak convergence (see details below at \eqref{ExplicationCovergence}), here is what we shall do : \\ \\
$\bullet$ \textbf{Step one :} We set 
  
\begin{equation} \label{182141}
\sigma^2 \ := \  \sum_{i,i'} \nu_i \ol{\nu}_{i'} \frac{b^2}{z_i \ol{z}_{i'}} \frac{\bb_{i'}^* \bb_i\bc_i^* \bc_{i'}}{z_i \ol{z}_{i'}-b^2}>0, 
\ee
and prove that for all fixed integer $k_0$, there is $\eta_{k_0}$ \st 
\begin{eqnarray} \label{goal1}
G_{n,k_0} :=\sqrt{n}\sum_{t=1}^p \nu_t \sum_{k=1}^{k_0} \frac{\bc^*_{t} \bA^k \bb_{t} }{z_t^{k+1}}  & \cloi & Z_{k_0} \ \defeloi \ \NN\Big(0, \sigma^2 - \eta_{k_0} \Big) ,
\end{eqnarray}
ant that 
$\eta_{k_0}  \to 0$ when $k_0 \to \infty$.
 Note that $\si^2$ doesn't depend on   $n$ thanks to \eqref{rem240214}. \\ \\
$\bullet$ \textbf{Step two :} We show that the rest shall be neglected for  large enough $k_0$. More precisely, for all $\delta>0$, we prove that there exists a large enough integer $k_0$ such that
\begin{eqnarray} \label{goal2}
\limsup_{n \to \infty} \E \left(\one_{\EE_n}\times \Big|\sqrt{n} \sum_{t=1}^p \nu_t  \sum_{k>k_0} \frac{\bc^*_{t} \bA^k \bb_{t} }{z_t^{k+1}}\Big|\right) & \leq & \delta.
\end{eqnarray}
(for $\EE_n$ the event of Lemma \ref{lem1} above).
 After that, we shall easily conclude. Indeed, to prove that $G_n$ converges in distribution to $\NN\left(0, \sigma^2 \right)$   it   suffices to prove that, for any Lipstichtz bounded test function $F$ with Lipschitz constant  $\mathcal{L}_F$, 
\begin{eqnarray*}
\E \left[F(G_n)\right] & \tto & \E \left[ F(Z) \right],
\end{eqnarray*} 
where $Z$ is a random variable such that $Z \eloi \NN\left(0, \sigma^2 \right)$. So, we write 
\begin{eqnarray}\label{ExplicationCovergence} 
\ \ \left|\E \left[F(G_n) - F(Z)  \right] \right| & \leq  & \left| \E \left[ F(G_n) - F(G_{n,k_0}) \right] \right|+\left| \E \left[ F(G_{n,k_0}) - F(Z_{k_0}) \right]\right| +\\ \nonumber &&   \left| \E \left[ F(Z_{k_0})- F(Z) \right] \right|\\
 \nonumber  &  \leq&   2\|F\|_\infty\pro\lf(\EE_n^c\ri)+\mathcal{L}_F   \E \left(\one_{\EE_n}\times  \Big|\sqrt{n} \sum_{t=1}^p \nu_t  \sum_{k>k_0} \frac{\bc^*_{t} \bA^k \bb_{t} }{z_t^{k+1}}\Big|\right)+ \\ \nonumber &&\left| \E \left[ F(G_{n,k_0})\right] - \E\left[F(Z_{k_0}) \right]\right| + \mathcal{L}_F \E \left| Z_{k_0} - Z \right|\nonumber \end{eqnarray}
which can be made as small as needed by \eqref{goal1} and \eqref{goal2} if $Z$ and $Z_{k_0}$ are coupled in the right way.\\ \\
$\bullet$ \textbf{Proof of step one : }\emph{Convergence of the finite sum}. \\
\ind Let us fix a positive integer $k_0$. Our goal here is to determine the limits of all the moments of the r.v. $G_{n,k_0}$ defined at \eqre{goal1} to conclude it is indeed asymptotically Gaussian. More precisely, we have
\begin{lemme} \label{lem1527040214}
\indent There exists $\sigma >0$ and $\eta_{k_0} $ such that $\lim_{k_0 \to \infty} \eta_{k_0}  = 0$ and  such that for all large enough $k_0$ and all non negative distinct integers $q,s$,
 $$
\E \left[ \left|G_{n,k_0} \right|^{2q} \right] \ = \  q!\cdot (\sigma^2-\eta_{k_0} )^q + o(1)\qquad \trm{ and }\qquad 
\E \left[ G_{n,k_0}^q \ol{G_{n,k_0}^s}\right] \ = \ o(1). $$
\end{lemme}
To prove Lemma \ref{lem1527040214}, we need to recall a main result about integration with respect to the Haar measure on unitary group, (see \cite[Cor. 2.4 and Cor. 2.7]{COL}),
\begin{propo} \label{wg}
\indent Let $k$ be a positive integer and $U=(u_{i,j})$ a Haar-distributed matrix. Let $(i_1,\ldots,i_k)$, $(i'_1,\ldots,i'_k)$, $(j_1,\ldots,j_k)$ and $(j'_1,\ldots,j'_k)$ be four $k$-tuple of $\left\{1,\ldots,n \right\}$. Then 
\begin{eqnarray} \label{wg1}
\E \left[ u_{i_1,j_1}  \cdots u_{i_k,j_k} \ol{ u_{ i'_1 ,j'_1 } }  \cdots \ol{u_{i'_k,j'_k}} \right] \ = \ \sum_{\sigma, \tau \in S_k} \delta_{i_1,i'_{\sigma(1)}} \ldots \delta_{i_k,i'_{\sigma(k)}} \delta_{j_1,j'_{\tau(1)}} \ldots \delta_{j_k,j'_{\tau(k)}} \wg( \tau \sigma^{-1}),
\end{eqnarray}
where $\wg$ is a function called  the \emph{Weingarten function}. Moreover, for $\sigma \in S_k$,   the asymptotical behavior of $\wg(\sigma)$ is given by   
\begin{eqnarray} \label{wg2}
n^{k+|\sigma|}\wg(\sigma) &=& \Moeb(\sigma) + O\left( \frac{1}{n^2}\right),
\end{eqnarray}
where $|\sigma|$ denotes the minimal number of factors necessary to write $\sigma$ as a product of transpositions, and $\Moeb$ denotes a function called  the \emph{M\"obius function}.
\end{propo}
\begin{rem} \label{moebius}
\indent a) The permutation $\sigma$ for which $\wg(\sigma)$ will have the largest order is the only one satisfying $|\sigma|=0$, i.e. $\sigma=id$. As a consequence, the only thing we have to know here about the M\"obius function is that $\Moeb(id) = 1$ (see \cite{COL}).\\ 
\indent b)  Notice that if for all $p \neq q$,   $i_p \neq i_q$ and $j_p \neq j_q$, then there is at most one non zero term in  the RHT of (\ref{wg1}).
\end{rem}
Lemma \ref{lem1527040214} follows from the following technical lemma  (we use the index $m$ in $\{ \cdot \}_m$ to denote 
a \emph{multiset}, i.e.  $\{x_1, \ld, x_k\}_m$ is the class of the $k$-tuple $(x_1, \ld, x_k)$ under the action of the symmetric group $S_k$).

\begin{lemme} \label{lemo(1)bis28114}
\indent Let $k_1,\ldots,k_q$ and $l_1,\ldots,l_s$ be some positive integers, let $i_1, \ld, i_q,i_1',  \ld, i'_s$ be some integers of $\{1, \ld, r\}$. Then :
\begin{enumerate}
\ite 
  If  $\left\{k_1,\ldots,k_q \right\}_m \neq \left\{l_1,\ldots,l_s \right\}_m$, we have
\beq
\E \left[\sqrt{n}\bc_{i_{1}}^* \bA^{k_1} \bb_{i_{1}}\cdots \sqrt{n}\bc_{i_{q}}^* \bA^{k_q} \bb_{i_{q}} \ol{\sqrt{n}\bc_{i'_{1}}^* \bA^{l_1} \bb_{i'_{1}} }\cdots \ol{\sqrt{n}\bc_{i'_{s}}^* \bA^{l_s} \bb_{i'_{s}} }\right] & = & o \left( 1 \right)
\eeq
\ite In the other case, $s=q$ and one can suppose that $l_1=k_1, \ld, l_q=k_q$. Under such an assumption, we have 
\beq
&&
\E \left[\sqrt{n}{\bc_{i_{1}}^*} \bA^{k_1} \bb_{j_{1}}\cdots \sqrt{n}{\bc_{i_{q}}^*} \bA^{k_q} \bb_{j_{q}} \ol{\sqrt{n}\bc_{i'_{1}}^* \bA^{l_1} \bb_{i'_{1}} }\cdots \ol{\sqrt{n}\bc_{i'_{s}}^* \bA^{l_s} \bb_{i'_{s}} }\right]\\ 
&  & \qquad\qquad\qquad \ = \  
b^{2(k_1+\cd+k_q)}\sum_{\si\in S_{k_1, \ld, k_q}}\prod_{t=1}^q\bb_{i'_{\si(t)}}^*\bb_{i_{t}}\bc_{i_{t}}^*\bc_{i'_{\si(t)}}+
o \left( 1 \right)
\eeq
where $S_{k_1, \ld, k_q}$ is the set of permutations of $\{1, \ld, q\}$ \st for each $t=1,\ld, q$, $k_t=k_{\si(t)}$. 

\ite Moreover,
\begin{eqnarray*}
&&\sum_{\displaystyle^{1\le k_1,\ldots,k_q =k_0}_{1\le k'_1,\ldots,k'_q \le k_0}}  \E \left[ \sqrt{n}\bc_{i_{1}}^* \frac{\bA^{k_1}}{z_{i_1}^{k_1+1}} \bb_{i_{1}}\ol{\sqrt{n}\bc_{i'_{1}}^* \frac{\bA^{k'_1}}{z_{i'_1}^{k'_1+1}} \bb_{i'_1}}\cdots \sqrt{n}\bc_{i_{q}}^* \frac{\bA^{{k}_q}}{z_{i_q}^{k_q+1}} \bb_{i_{q}}\ol{\sqrt{n}\bc_{i'_{q}}^* \frac{\bA^{k'_q}}{z_{i'_q}^{k'_q+1}} \bb_{i'_{q}}} \right]\\
 & = & \sum_{\si \in S_q} \prod_{t=1}^q \frac{b^{2}}{z_{i_t} \ol{z_{i'_{\si(t)}}}} \frac{1 - \left(\frac{b^{2}}{z_{i_t} \ol{z_{i'_{\si(t)}}}}\right)^{k_0}}{z_{i_t} \ol{z_{i'_{\si(t)}}}-b^2} \bb^*_{i'_{\si(t)}} \bb_{i_t} \bc^*_{i_t} \bc_{i'_{\si(t)}}   + o\left(1\right).
\end{eqnarray*}
\end{enumerate}
\end{lemme}

\indent Let us briefly explain  the main ideas of the proof of this lemma (detailed proof is given in   Section \ref{tec}). First, let us recall that $\bA = \bU \bT$, so that these expectations expand as sums of   terms   as
\begin{eqnarray*}
&&  \E \left[u_{i_{0,1},i_{1,1}}\cdots u_{i_{k_1-1,1},i_{k_1,1}} u_{i_{0,2},i_{1,2}} \cdots u_{i_{k_r-1,r},i_{k_r,r}} \ol{u_{j_{0,1},j_{1,1}}\cdots u_{j_{l_1-1,1},j_{l_1,1}} u_{j_{0,2},j_{1,2}} \cdots u_{j_{l_s-1,s},j_{l_s,s}}}\right].
\end{eqnarray*}
If the $u_{i,j}$'s were  independent and distributed as $\NN\left(0,\frac{1}{{n}}\right)$, the result would be easily proved because most of these expectations would be equal to zero. In our case, the difficulty is that, according to Proposition \ref{wg}, lots of expectations do not vanish and they are expressed with the Weingarten function (which is a very complicated function). However, we notice that when these expectations   do not vanish as in the Gaussian case, $\wg(id)$ never occurs in \eqref{wg1}, so that they  are negligible thanks to \eqref{wg2}.    \\

\indent At last, it is easy to conclude the proof of Lemma \ref{lem1527040214} thanks to Lemma \ref{lemo(1)bis28114}. Indeed, for any integers $q \neq s$, we have from (1) of Lemma \ref{lemo(1)bis28114} that $\E \left[G_{n,k_0}^q \ol{G_{n,k_0}^s} \right] = o(1)$. Moreover, we have
\beq
\E \left[\big|G_{n,k_0}\big|^{2q} \right] & = & n^{q} \E \left[ \left(\sum_{i} \nu_{i}\bc_i^* \sum_{k=1}^{k_0}  \frac{\bA^{k}}{z_i^{k+1}} \bb_i\right)^q \left(\ol{\sum_{i'} \nu_{i'}\bc_{i'}^* \sum_{k'=1}^{k_0}  \frac{\bA^{k'}}{\tta_{i'}^{k'+1}} \bb_{i'}}\right)^q\right] \\
& = &  \sum_{\displaystyle^{i_1,\ldots,i_q}_{i'_1,\ldots,i'_q}} \prod_{t=1}^q \nu_{i_t}\ol{\nu}_{i'_t}  
 \sum_{\displaystyle^{{k}_1,\ldots,{k}_q}_{k'_1,\ldots,k'_q}} \E \Big[ \sqrt{n}\bc_{i_{1}}^* \bA^{k_1} \bb_{i_{1}}\ol{\sqrt{n}\bc_{i'_{1}}^* \bA^{k'_1} \bb_{i'_{1}}}\cdots \sqrt{n}\bc_{i_{q}}^* \bA^{k_q} \bb_{i_{q}}\ol{\sqrt{n}\bc_{i'_{q}}^* \bA^{k'_q} \bb_{i'_{q}}} \Big]\\
 & = &\sum_{\displaystyle^{i_1,\ldots,i_q}_{i'_1,\ldots,i'_q}}  \prod_{t=1}^q \nu_{i_t}\ol{\nu}_{{i'_t}}   \sum_{\si \in S_q} \prod_{t=1}^q \frac{b^{2}}{z_{i_t} \ol{z}_{i'_{\si(t)}}} \frac{1 - \left(\frac{b^{2}}{z_{i_t} \ol{z}_{i'_{\si(t)}}}\right)^{k_0}}{z_{i_t} \ol{z}_{i'_{\si(t)}}-b^2} \bb^*_{i'_{\sigma(t)}} \bb_{i_t} \bc^*_{i_t} \bc_{i'_{\sigma(t)}} + o\left( 1\right) \\
  & = & q! \times \left( \sigma^2 -\eta_{k_0}  \right)^q + o\left(1\right) \\ 
\eeq
where for $\si$ is given by \eqre{182141} and $|\eta_{k_0} | < \sigma^2 \cdot (\frac{b}{b+2\ep})^{2k_0}$. \\ \\
$\bullet$
\textbf{Proof of step two : }\emph{Vanishing of the tail of the sum}. \\
\ind Our goal here is to prove that the rest can be neglected, i.e.  that for all $\delta>0$, there exists a large enough integer $k_0$ such that for any $t \in \{1,\ldots,p\}$ and for $\EE_n$ the event of Lemma \ref{lem1} above,
\begin{eqnarray}\label{aprouver13813}
\limsup_{n \to \infty}\E \left(\one_{\EE_n}\times \Big|\sqrt{n} \sum_{k>k_0} \frac{\bc^*_{t} \bA^k \bb_{t} }{z_{t}^{k+1}}\Big|\right) & \leq & \delta. 
\end{eqnarray}
\indent First, using the fact that
 $$
 \E \left[\one_{\EE_n}\times \Big| z_t^{-k-1}   \bc^*_{t} \bA^k \bb_{t}  \Big|\right] \leq \| \bB 
\|_{\op} \| \bC \|_{\op} C_1 \f{(b+\ep)^{k}}{(b+2\ep)^{k+1}},
$$
it is easy to show that for a large enough positive constant $C$ (depending only on $\ep$), we have
$$
\sqrt{n}\sum_{k > C \log n} \E \left[\one_{\EE_n}\times \Big|  z^{-k-1}_t  \bc^*_{t} \bA^k \bb_{t}  \Big|\right]  \ \eg \ o\left( 1\right).
$$

\indent Now, we only need to prove that
$$
\forall \delta>0, \ \exists k_0, \ \textrm{ for all $n$ large enough, } \ \sum_{k_0 < k < C \log n} \f{\sqrt{n}}{|z_t|^{k+1}} \sqrt{\E \left[\one_{\EE_n}\times \Big|    \bc^*_{t} \bA^k \bb_{t}  \Big|^2\right]}\le \delta.
$$
 At first, we notice that 
$$
\E \left(\one_{\EE_n}\times \Big|\sqrt{n} \sum_{k=k_0}^{C\log n} \frac{\bc^*_{t} \bA^k \bb_{t} }{z_{t}^{k+1}}\Big|\right) \ \le \   \sum_{k=k_0}^{C\log n}\f{\sqrt{n}}{|z_t|^{k+1}} \E \left(\one_{\EE_n}\times \Big|    \bc^*_{t} \bA^k \bb_{t}  \Big|\right) \ \le \   \sum_{k=k_0}^{C\log n}\f{\sqrt{n}}{|z_t|^{k+1}} \sqrt{\E \left[\one_{\EE_n}\times \Big|    \bc^*_{t} \bA^k \bb_{t}  \Big|^2\right]}
$$ 
Then we condition with respect to the $\sigma$-algebra of $\bU$, i.e. write $$\E \left[\one_{\EE_n}\times \Big|    \bc^*_{t} \bA^k \bb_{t}  \Big|^2\right]= \E\left[\one_{\EE_n}\times\E_{\bP}\lf( \Big|    \bc^*_{t} \bA^k \bb_{t}  \Big|^2\right)\ri]=\E\left[\one_{\EE_n}\times\E_{\bP}\lf(    \bc^*_{t} \bA^k \bb_{t}  \bb^*_{t}(\bA^*)^k\bc_{t}\right)\ri].$$
Let us now remember that we have supposed, at \eqref{88132}, that $\bP=\bB\bC$ is invariant, in law, by conjugation by any unitary matrix. Hence one can introduce a Haar-distributed unitary matrix $\bV$, independent of all other random variables, and write $\bP\eloi \bV\bP\bV^*$, so that \beqy\label{148131}
 \E_{\bP}\lf(    \bc^*_{t} \bA^k \bb_{t}  \bb^*_{t}(\bA^*)^k\bc_{t}\right) &= & \E_{\bP}\lf(   \Tr \bA^k \bb_{t}  \bb^*_{t}(\bA^*)^k\bc_{t}  \bc^*_{t}\right) \\ \nonumber &=& \E_{\bP}\lf(\E_{\bV}\lf( \Tr \bA^k \bV\bb_{t}  \bb^*_{t}\bV^*(\bA^*)^k\bV\bc_{t}  \bc^*_{t}\bV^*\ri)\ri),\eeqy where $\E_{\bV}$ denotes the expectation with respect to the randomness of $\bV$.

 Then, we shall use the following lemma, whose proof is postponed to Section \ref{preuvelemmecalculexpectation13813}.
 \begin{lemme}\label{lemmecalculexpectation13813}
 Let $\bV$ be an $n\times n$ Haar-distributed  unitary matrix and let  $\bA$, $\bB$, $\bC$, $\bD$ be some deterministic $n\times n$ matrices. Then \beqy\label{LHT13813} \E\Tr \bA\bV\bB\bV^*\bC\bV\bD\bV^* &=&\ff{n^2-1}\lf\{\Tr \bA\bC\Tr\bB\Tr\bD+\Tr \bA\Tr\bC\Tr\bB\bD\ri\}\\ \nonumber&&\qquad\qquad-\ff{n(n^2-1)}\lf\{\Tr \bA\bC\Tr \bB \bD+\Tr \bA\Tr\bC\Tr \bB \Tr\bD\ri\}.\eeqy
 \end{lemme}
By this lemma, one easily gets  
$$|\E_{\bV}\lf( \Tr \bA^k \bV\bb_{t}  \bb^*_{t}\bV^*(\bA^*)^k\bV\bc_{t}  \bc^*_{t}\bV^*\ri)|\le\f{2}{n-1}\left(\|\bA^k\|_{\op}^2 + \big|\tr\big( \bA^k\big)\big|^2 \right)\|\bB\|^2_{\op}\|\bC\|^2_{\op} $$ 
hence as $\bB$ and $\bC$ are supposed to be bounded, there is a constant $C$ such that  
$$ \E_{\bP}\lf(    \bc^*_{t} \bA^k \bb_{t}  \bb^*_{t}(\bA^*)^k\bc_{t}\right)\le \f{C}{n}\left(\|\bA^k\|_{\op}^2+\big|\tr \big( \bA^k\big)\big|^2\right).$$
Then, we use the following lemma, a weaker version of \cite[Theorem 1]{BEN}.
\begin{lemme} \label{lemmeFlorent}
\indent There exists a positive constant $K$ such that for all $k \leq C\log n$,   for all large enough $n$,
$$
\E \left[ \big|\tr\big( \bA^k\big)\big|^2 \right] \ \leq \ K\left( b+\ep\right)^{2k}.
$$
\end{lemme}

By \eqref{148131} and Lemma \ref{lemmeFlorent}, for all $k\le C \log n$, there exists some positive constant $C'$ such that $$\sqrt{\E \left[\one_{\EE_n}\times \Big|    \bc^*_{t} \bA^k \bb_{t}  \Big|^2\right]}\le \f{C'(b+\ep)^k}{\sqrt{n}}.$$ Hence  as $|z_t|\ge b+2\ep$ for $n$ large enough, \eqref{aprouver13813} is proved.

\subsubsection{Proof of Lemma \ref{lem040214bis}} The proof  relies on the same tricks of the proof of Lemma \ref{lem040214},  using the already noticed fact that 
 for $|z|>b+2\ep$,    $$\one_{\EE_n}\bX_n^z \ = \ \sqrt{n} \one_{\EE_n}\sum_{k\ge 1}\bC\f{\bA^k}{z^{k+1}}\bB,$$ so that $$\f{\partial}{\partial z}\one_{\EE_n}\bX_n^z \ = \ -\sqrt{n} \one_{\EE_n}\sum_{k\ge 1}(k+1)\bC\f{\bA^k}{z^{k+2}}\bB.$$ 

\subsection{Proof of Lemma \ref{lemaux0.1}}\la{1421416h59}
\indent To prove Lemma \ref{lemaux0.1}, we shall need to do a Taylor expansion of $F^{\tta_i}_{j}(z)$. From now on, we fix a compact set $K$ and consider $z \in K$. Recall that $F_j^{\tta_i}(z)$ and $\bX_n^z$ have been defined respectively at \eqre{162142} and \eqre{162141} as
 $$F_j^{\tta_i}(z)=\det\lf(\bI-\bC\lf(\lf(\tta_i+n^{-1/(2p_{i,j})}z \ri)\bI-\bA\ri)^{-1}\bB\ri)\qquad \qquad \bX_n^z=\sqrt{n}\bC \left(\left(z\bI -\bA  \right)^{-1}- z^{-1}\right)\bB,$$ 
 hence, using Lemma \ref{lem040214bis} and the convergence of $\bX_n^z$ to $\bX^z$ established at Section \ref{biggestpart},  
\beq
F^{\tta_i}_{j}(z) & = & \det \left(\bI - \frac{1}{\tta_i + n^{-1/(2p_{i,j})}z} \bJ - \frac{1}{\sqrt{n}}\bX_n^{\tta_i + n^{-1/(2p_{i,j})}z} \right) \\
& = & \det \left( \bI - \tta_i^{-1}\bJ + \tta_i^{-1}\left( 1 - \frac{1}{1+n^{-1/(2p_{i,j})} z \tta_i^{-1}} \right)\bJ - \frac{1}{\sqrt{n}}\bX^{\tta_i} + o\left( \frac{1}{\sqrt{n}}\right) \right)\\
& = & \det \left(\bI - \tta_i^{-1}\bJ + z \delta_n \tta_i^{-2}\bJ + \frac{1}{\sqrt{n}} \bG \right)
\eeq
where we define $$\delta_n: =\f{\tta_i}{z}\left( 1 - \frac{1}{1+n^{-1/(2p_{i,j})} z \tta_i^{-1}} \right)= n^{-1/(2p_{i,j})}(1+o(1))\qquad\trm{  and }\qquad \bG=-\bX^{\tta_i}+o(1).$$ \\   
\indent Let us write $\bJ$ by blocks
$$
\bJ \ = \ \left(
\begin{tabular}{c|c|c}
$\ast$ &  $(0)$  & $(0)$ \\
\hline 
$(0)$ & $ \bJ(\tta_i)$ &  $(0)$ \\
\hline 
$(0)$  & $(0)$ & $ \ast$\\
\end{tabular}
\right)
$$
where $\bJ(\tta_i)$ is the part with the blocks associated to $\tta_i$. And so, we write
$$
\bI - \tta_i^{-1}\bJ \ = \ \left(
\begin{tabular}{c|c|c}
$\bN'$ &  $(0)$  & $(0)$ \\
\hline 
$(0)$ & $ \bN$ &  $(0)$ \\
\hline 
$(0)$  & $(0)$ & $ \bN''$\\
\end{tabular}
\right)$$
where $\bN'$ and $\bN''$ are invertible matrices and $\bN$ is the diagonal by blocks matrix 
\beqy\nonumber \bN&=&\bI - \tta_i^{-1}\diag(\underbrace{\bR_{p_{i,1}}(\tta_i), \ld, \bR_{p_{i,1}}(\tta_i)}_{\trm{$\bet_{i,1}$ blocks}}, \ld\ld, \underbrace{\bR_{p_{i,\alp_i}}(\tta_i), \ld, \bR_{p_{i,\alp_i}}(\tta_i)}_{\trm{$\bet_{i,\alp_i}$ blocks}})\\ \la{2811415h}&=&-\tta_i^{-1}\diag(\underbrace{\bR_{p_{i,1}}(0), \ld, \bR_{p_{i,1}}(0)}_{\trm{$\bet_{i,1}$ blocks}}, \ld\ld, \underbrace{\bR_{p_{i,\alp_i}}(0), \ld, \bR_{p_{i,\alp_i}}(0)}_{\trm{$\bet_{i,\alp_i}$ blocks}})\eeqy
with   $\bR_p(\tta)$ as defined at \eqre{2811410h21} for $p$ an integer and $\tta\in \C$. 
 

Let us now expand  the determinant $\det \left(\bI - \tta_i^{-1}\bJ + z \delta_n \tta_i^{-2}\bJ + n^{-1/2} \bG \right)$ using the \emph{columns replacement approach} of following formula, where the $M_k$'s and the $H_k$'s are the columns of two $r\ti r$ matrices $\bM$ and $\bH$ (that one will think of as an error term, even though the formula below is exact)
\beq
 \det\left( \bM + \bH \right) & = & \det \left(\bM \right) + \sum_{k=1}^r \det\left(M_1 \big| M_2 \big|\ldots \big| H_k \big| \ldots \big| M_r \right) 
 +  \sum_{1 \leq k_1<k_2 \leq r} \det \left(M_1 \big| \ldots \big|H_{k_1} \big| \ldots \big|H_{k_2} \big|\ldots \big| M_r \right)  \\
 & & + \quad \ldots \quad  +  \sum_{1 \leq k_1<k_2<\cdots<k_s \leq r}  \det \left(M_1\big| \ldots\big|H_{k_1}\big| \ldots \big| H_{k_s}\big| \ldots \big| M_r \right) \\ 
& &+ \quad \ldots \quad + \sum_{k=1}^r \det\left(H_1\big| H_2\big| \ldots \big| M_k \big| \ldots\big| H_r \right) + \det\left( \bH\right) \\
\eeq
We shall use this formula with  $\bM = \bI - \tta_i^{-1}\bJ$ and $\bH = z \delta_n \tta_i^{-2}\bJ + n^{-1/2} \bG$, and we shall keep only higher terms. It means that the determinant is a summation of determinants of $\bM$ where some of the columns are replaced by the corresponding column of $z \delta_n \tta_i^{-2}\bJ$ or of $n^{-1/2}\bG$. Recall that $\bM$ has several columns of zeros (the ones corresponding to null columns of $\bN$), so  we know that we have to replace at least these columns to get a non-zero determinant. Moreover, we won't replace the columns of $\bN'$ or $\bN''$ because  this would  necessarily make appear negligible terms (recall that $\bN'$ and $\bN''$ are invertible), so all the non-negligible determinants will be factorizable by $\det (\bN')\det (\bN'')$. So now, let us understand what are the non-negligible terms in the summation. 

To make things clear, let us start with an example. We choose $p_{i,j} = 3$ and the matrix $\bN$ given, via \eqre{2811415h}, by 
$$
\bN \ = \ -\tta_i^{-1} \left( 
\begin{tabular}{cccccccc}
0 & 1 & 0 & \multicolumn{1}{c|}{0} &  &  &  &   \\
0 & 0 & 1 & \multicolumn{1}{c|}{0} &  & \multicolumn{3}{c}{(0)}   \\
0 & 0 & 0 & \multicolumn{1}{c|}{1} & &  &  &   \\
0 & 0 & 0 & \multicolumn{1}{c|}{0} &  &  &  &  \\
\cline{1-7} 
 &  &  & \multicolumn{1}{c|}{} & 0 & 1 & \multicolumn{1}{c|}{0} &  \\
 &  &  & \multicolumn{1}{c|}{} & 0 & 0 & \multicolumn{1}{c|}{1} &   \\
\multicolumn{3}{c}{(0)}   & \multicolumn{1}{c|}{} & 0 & 0 & \multicolumn{1}{c|}{0} &  \\
\cline{5-8} 
 &  &  &  & &  & \multicolumn{1}{c|}{} & 0  \\
\end{tabular}
\right), 
$$
we know we have to replace at least $3$ columns (the first, the fifth and the last ones) which correspond to the first column of each diagonal blocks, and we shall deal with one block at the time. Let us deal with the first one. If we replace this column by the corresponding column of $z \delta_n\tta_i^{-2}\bJ$, we get


$$
{z \delta_n}\left|
\begin{array}{cccc}
\ff{\tta_i} & \ff{\tta_i} & 0 & 0 \\
0 & 0 & \ff{\tta_i} & 0 \\
\multicolumn{2}{c}{\multirow{2}{*}{(0)}} & 0 & \ff{\tta_i} \\
\multicolumn{2}{c}{}                     & 0 & 0 \\
\end{array}
 \right|
$$

We see that in this case, some   non linearly independent columns appear. 
 It follows that once one has replaced a null column by a column from   $z \delta_n \tta_i^{-2}\bJ$, 
  the whole block needs to be replaced to get a non zero determinant :  
$$
 {z \delta_n}\left|
\begin{array}{llll}
\ff{\tta_i} & \ff{\tta_i} & 0 & 0 \\
0 & 0 & \ff{\tta_i} & 0 \\
\multicolumn{2}{c}{\multirow{2}{*}{(0)}}& 0 & \ff{\tta_i} \\
\multicolumn{2}{c}{}                    & 0 & 0 \\
\end{array}
 \right|  \tto \ 
 {(z \delta_n)^2}\left|
\begin{array}{llll}
\ff{\tta_i} & \ff{\tta_i^2} & 0 & 0 \\
0 & \ff{\tta_i} & \ff{\tta_i} & 0 \\
\multicolumn{2}{c}{\multirow{2}{*}{(0)}}& 0 & \ff{\tta_i} \\
\multicolumn{2}{c}{}                    & 0 & 0 \\
\end{array}
 \right| \tto \ {(z \delta_n)^3}\left|
\begin{array}{llll}
\ff{\tta_i} & \ff{\tta_i^2} & 0 & 0 \\
0 & \ff{\tta_i} & \ff{\tta_i^2} & 0 \\
\multicolumn{2}{c}{\multirow{2}{*}{(0)}} & \ff{\tta_i} & \ff{\tta_i} \\
\multicolumn{2}{c}{}                     & 0 & 0 \\
\end{array}
 \right| $$ $$ \tto \ {(z \delta_n)^4}\left|
\begin{array}{llll}
\ff{\tta_i} &  \ff{\tta_i^2} & 0 & 0 \\
0 & \ff{\tta_i} &  \ff{\tta_i^2}& 0 \\
\multicolumn{2}{c}{\multirow{2}{*}{(0)}}& \ff{\tta_i} &  \ff{\tta_i^2} \\
\multicolumn{2}{c}{}                    & 0 & \ff{\tta_i}\\
\end{array}
 \right|, 
$$


 Another possibility to fill a null column would be to replace it by the corresponding  one in $n^{-1/2} \bG$:
$$
\left|
\begin{array}{llll}
g_{1,1} & \ff{\tta_i} & 0 & 0 \\
g_{2,1} & 0 & \ff{\tta_i} & 0 \\
g_{3,1} & 0 & 0 & \ff{\tta_i} \\
g_{4,1} & 0 & 0 & 0 \\
\end{array}
 \right|
$$
 We obtain  an invertible block directly (i.e. without having to replace the whole block as above). However,   in this example, $\delta_n \gg n^{-1/2}$  (because $p_{i,j}=3$), this term might be negligible. If $\delta_n^4 \gg n^{-1/2}$, then first choice is relevant (the other would be negligible), or else, if $\delta_n^4 \ll n^{-1/2}$, we would make the second choice. \\
\indent  Our strategy is to choose $\bJ$ on the blocks of size $p<p_{i,j}$ (because $\delta_n^{p} \gg n^{-1/2}$) and $\bG$ on the blocks of size $p>p_{i,j}$ (because $\delta_n^{p} \ll {\frac{1}{\sqrt{n}}}$). For the blocks of size $p=p_{i,j}$, we can choose both (because $\delta_n^{p_{i,j}} \approx  {\frac{1}{\sqrt{n}}}$). So in our example, the non negligible terms are
\beq
  \det\left( \frac{G_1}{\sqrt{n}}\big| N_2 \big| N_3 \big| N_4 \big|\f{z\delta_n}{\tta_i^2} J_5 \big|\f{z\delta_n}{\tta_i^2} J_6  \big| \f{z\delta_n}{\tta_i^2} J_7 \big| \f{z\delta_n}{\tta_i^2} J_8 \right)  & = &
-z^4 \cdot {\frac{\delta_n^4}{\sqrt{n}}}   \left| 
\begin{array}{cccccccccc}
g_{1,1} & \ff{\tta_i} & 0 & 0 & \multicolumn{4}{c}{\multirow{4}{*}{(0)}}  \\
g_{2,1} & 0 & \ff{\tta_i} & 0 & \multicolumn{4}{c}{}  \\
g_{3,1} & 0 & 0 & \ff{\tta_i} & \multicolumn{4}{c}{}   \\
g_{4,1} & 0 & 0 & 0 & \multicolumn{4}{c}{}  \\
g_{5,1} & \multicolumn{3}{c}{\multirow{4}{*}{(0)}}  & \ff{\tta_i} &  \ff{\tta_i^2} & 0 & 0 \\
g_{6,1} & \multicolumn{3}{c}{}                      & 0 & \ff{\tta_i} &  \ff{\tta_i^2} & 0  \\
g_{7,1} & \multicolumn{3}{c}{}                      & 0 & 0 & \ff{\tta_i} & 0 \\
g_{8,1} & \multicolumn{3}{c}{}                      & 0 & 0 & 0 & \ff{\tta_i} \\
\end{array}
\right|
\eeq
and
\beq
 \det\left( \frac{G_1}{\sqrt{n}}\big| N_2 \big| N_3 \big| N_4 \big|\frac{G_5}{\sqrt{n}} \big|N_6  \big| N_7 \big| \f{z \delta_n}{\tta_i^2} J_8 \right) & = & -z\cdot {\frac{\delta_n}{n}}  \cdot
\left| 
\begin{array}{cccccccccc}
g_{1,1} & \ff{\tta_i} & 0 & 0 & g_{1,5} &  \multicolumn{3}{c}{\multirow{4}{*}{(0)}}  \\
g_{2,1} & 0 & \ff{\tta_i} & 0 & g_{2,5} &  \multicolumn{3}{c}{}                      \\
g_{3,1} & 0 & 0 & \ff{\tta_i} & g_{3,5} &  \multicolumn{3}{c}{}                      \\
g_{4,1} & 0 & 0 & 0 & g_{4,5} &            \multicolumn{3}{c}{}                      \\
g_{5,1} &  \multicolumn{3}{c}{\multirow{4}{*}{(0)}}  & g_{5,5} & \ff{\tta_i} & 0 & 0 \\
g_{6,1} &  \multicolumn{3}{c}{}                      & g_{6,5} & 0 & \ff{\tta_i} & 0  \\
g_{7,1} &  \multicolumn{3}{c}{}                      & g_{7,5} & 0 & 0 & 0  \\
g_{8,1} &  \multicolumn{3}{c}{}                      & g_{8,5} & 0 & 0 & \ff{\tta_i}  \\
\end{array}
\right|
\eeq
 and one can easily notice that the  sum of the non negligible terms is
\beq
&&\det\left( \frac{G_1}{\sqrt{n}}\big| N_2 \big| N_3 \big| N_4 \big|\f{z\delta_n}{\tta_i^2} J_5 \big|\f{z\delta_n}{\tta_i^2} J_6  \big| \f{z\delta_n}{\tta_i^2} J_7 \big| \f{z\delta_n}{\tta_i^2} J_8 \right) + \det\left( \frac{G_1}{\sqrt{n}}\big| N_2 \big| N_3 \big| N_4 \big|\frac{G_5}{\sqrt{n}} \big|N_6  \big| N_7 \big| \f{z \delta_n}{\tta_i^2} J_8 \right) \\
&& \eg \f{z}{\tta_i^6} \ff{n^{1+\ff{2p_{i,j}}}} \bdet g_{4,1} & g_{4,5} \\ g_{7,1} & g_{7,5} - \f{z^3}{\tta_i} \\ \edet \ + \ o\left(\ff{n^{1+\ff{2p_{i,j}}}}\right).
\eeq

\indent Now that this example is well understood, let us treat the general case : \bgt
\ite[--] We know that there are $\bet_{{i,1}}+\cdots+\bet_{{i,j-1}}$ blocks of size larger than $p_{i,j}$ so we will replace the first column of each of these blocks by the corresponding column of $n^{-1/2}\bG$. \\
\ite[--] For all the blocks of lower size, we replace all the columns by the corresponding column of $z \delta_n \tta_i^{-2}\bJ$. The number of such columns is   $\pi_{i,j}:=\beta_{i,j+1} \times p_{i,j+1}+\cdots \beta_{i,\alp_i}\times p_{i,\alp_i}$. \\
\ite[--] We also know that there are $\bet_{{i,j}}$ blocks  of size $p_{i,j}$ and for each block, we have two choices so that represents $2^{\bet_{{i,j}}}$ non negligible terms.\\
\ent

 And so, we conclude that : \bgt
\ite  The statement holds for $\displaystyle \gamma_{{i,j}} = \frac{1}{2}\sum_{l=1}^{j-1} \bet_{{i,l}}+\frac{\pi_{{i,j}}}{2p_{i,j}}+\frac{1}{2} \bet_{{i,j}}$. \\
\ite  All the non negligible terms are factorizable by $z^{\pi_{{i,j}}}$. \\
\ite  Using notations from \eqref{section2.3}, we define the matrices 
 $$\ope{M}^{\tta_i,\mathrm{I}}_{j}\egd[g^{\tta_i}_{k,\ell}]_{\ds^{k\in K(i,j)^-}_{\ell \in L(i,j)^-}} \qquad \qquad \qquad \ope{M}^{\tta_i,\dI\dI}_{j}\egd[g^{\tta_i}_{k,\ell}]_{\ds^{k\in K(i,j)^-}_{\ell\in L(i,j)}} $$ $$\ope{M}^{\tta_i,\dI\dI\dI}_{j}\egd[g^{\tta_i}_{k,\ell}]_{\ds^{k\in K(i,j)}_{\ell\in L(i,j)^-}} \qquad \qquad \qquad \ope{M}^{\tta_i,\dI\dV}_{j}\egd[g^{\tta_i}_{k,\ell}]_{\ds^{k\in K(i,j)}_{\ell \in L(i,j)}}$$
 
and with a simple calculation, one can sum up all the non-negligible terms by
$$
C \cdot \frac{z^{\pi_{i,j}}}{n^{\gamma_{i,j}}} \bdet \ \ \ope{M}^{\tta_i,\dI}_{j} & \ope{M}^{\tta_i,\dI\dI}_{j} \\ \ \ \ope{M}^{\tta_i,\dI\dI\dI}_{j} & \ope{M}^{\tta_i,\dI\dV}_{j} - \f{z^{p_{i,j}}}{\tta_i}\bI_{\bet_{i,j}} \\
\edet + o \left(\ff{n^{\gamma_{i,j}}} \right) 
$$
where $C$ is a deterministic constant equal to $\pm$ a power of $\tta_i^{-1}$.  Then, using a well-know formula (see for example Eq. (A1) of \cite{agz} p. 414), we have
$$
 \bdet \ \ \ope{M}^{\tta_i,\dI}_{j} & \ope{M}^{\tta_i,\dI\dI}_{j} \\ \ \ \ope{M}^{\tta_i,\dI\dI\dI}_{j} & \ope{M}^{\tta_i,\dI\dV}_{j} - \f{z^{p_{i,j}}}{\tta_i}\bI_{\bet_{i,j}} \\
\edet \eg \tta_i^{-\bet_{i,j}}\det \left(\ope{M}^{\tta_i,\dI}_{j} \right) \det \left( \tta_i(\ope{M}^{\tta_i,\dI\dV}_{j} - \ope{M}^{\tta_i,\dI\dI\dI}_{j}(\ope{M}^{\tta_i,\dI}_{j})^{-1}\ope{M}^{\tta_i,\dI\dI}_{j})- z^{p_{i,j}}\bI_{\bet_{i,j}} \right).
$$
\ite Thanks to Lemma \ref{lem040214}, we know that 
\beq
\E\lf( {m}^{\tta_i}_{k,\ell} \; {m}^{\tta_{i'}}_{k',\ell'}\ri)\ =\ 0, \qquad \E\lf( {m}^{\tta_i}_{k,\ell} \; \ovl{{m}^{\tta_{i'}}_{k',\ell'}}\ri)\eg \f{b^2}{\tta_i\ovl{\tta_{i'}}}\; \ff{\tta_i\ovl{\tta_{i'}}-b^2} \;\be_{ k}^*\bC\bC^*\be_{ {k'}}\;\be_{ \ell'}^*\bB^*\bB\, \be_{ {\ell}},
\eeq
and from \eqref{rem240214}, we write
$$
\be_{ k}^*\bC\bC^*\be_{ {k'}}\;\be_{ \ell'}^*\bB^*\bB \be_{\ell} \eg \be_{ k}^*\bQ^{-1}(\bQ^{-1})^*\be_{ {k'}}\;\be_{ \ell'}^*\bJ^* \bQ^* \bQ \bJ \be_{\ell}
$$
then, from the definition of the set $L(i,j)$, we know that if $\ell \in L(i,j)$ then $\bJ \be_{\ell} = \tta_i \be_{\ell}$, so, finally, 
\beq
\E\lf( {m}^{\tta_i}_{k,\ell} \; {m}^{\tta_{i'}}_{k',\ell'}\ri)\ =\ 0, \qquad \E\lf( {m}^{\tta_i}_{k,\ell} \; \ovl{{m}^{\tta_{i'}}_{k',\ell'}}\ri)\eg  \f{b^2}{\tta_i\ovl{\tta_{i'}}-b^2} \;\be_{ k}^*\bQ^{-1}(\bQ^{-1})^*\be_{ {k'}}\;\be_{ \ell'}^*\bQ^*\bQ\, \be_{ {\ell}}.
\eeq

\ent

 \section{Proofs of the technical results} \label{tec}
\subsection{Proofs of Lemmas \ref{lem1} and \ref{lem0}} \label{part1}
\begin{lemme} \label{lem2}
\indent There exists a constant $C_1$, independent of $n$, such that with  probability tending to one, 
\begin{eqnarray*}
\sup_{ |z| = b + \ep} \left\|(z\bI - \bA)^{-1}\right\|_{\op} &\leq& C_1.  
\end{eqnarray*}
\end{lemme}
\bpr
Note that for any $\eta>0$, 
$$ \|(z\bI-\bA)^{-1}\|_{\op} \ \leq \ \frac{1}{\eta} \iff \nu^z ([-\eta, \eta])=0 ,$$
where $\nu^z: = \frac{1}{2n} \sum_{i=1}^n (\delta_{-s_i^z} + \delta_{s_i^z})$ and 
the $s_i^z$'s are the singular values of $z\bI - \bA$.  

By Corollary 10 of \cite{GUI2}, for any $z$ \st $|z|>b$, there is $\bet_z$ \st 
with  probability tending to one, 
$$ \nu^z ([- \bet_z ,  \bet_z ])=0.$$ It follows from standard perturbation inequalities that  with  probability tending to one,  for any $z'$ \st $|z'-z|< \f{\bet_z}{2}$, $$ \nu^{z'} ([-\f{\bet_z}{2}, \f{\bet_z}{2}])=0.$$ 
Then with a  compacity argument, one concludes easily.\epr

\subsubsection{Proof of Lemma \ref{lem1}}  
Note first that thanks to the Cauchy formula, for all $x\in \C$,     \begin{eqnarray*}
|x| < b+\ep  &\implies& \forall k\ge 0,\quad x^k=\frac{1}{2 i \pi} \int_{|z|=b+\ep} \frac{z^k}{z-x} dz.
\end{eqnarray*} Moreover, by \cite[Th. 2]{GUI2}, the spectral radius of $\bA$ converges in probability to $b$, so that with probability tending to one, by application of the holomorphic functional calculus (which is working even for non-Hermitian matrices) to $\bA$,
\begin{eqnarray*}  \forall k\ge 0,\qquad
\bA^k &=& \frac{1}{2 i \pi} \int_{|z|=b+\ep} z^k \left(z-\bA\right)^{-1} dz.
\end{eqnarray*} Thus  with probability tending to one, \begin{eqnarray*}  \forall k\ge 0,\qquad 
\|\bA^k \|_{\op} & \leq  &  \frac{1}{2 \pi} \sup_{|z|=b+\ep} \|(z\bI-\bA)^{-1}\|_{\op} \times \int_{|z|=b+\ep} |z|^k dz.
\end{eqnarray*}
Then one concludes using the previous lemma.

\subsubsection{Proof of Lemma  \ref{lem0}} Since $\bC \bA^k \bB$ is a square  $r \times r$ matrix, it suffices to prove that each entry tends, in probability, to $0$. And since $\bB$ and $\bC$ are uniformly bounded (see Remark \ref{165130092014}) , one just has to show that for all unit vectors $\bb$ and $\bc$,  
\begin{eqnarray} \label{2}
\bc^* \bA^k \bb & \cvp & 0.
\end{eqnarray}
\indent Recall that $\bA = \bU \bT$ and 
\begin{eqnarray*}
\bc^* \bA^k \bb &=&  \sum_{i_0,i_1,\ldots,i_{k}} \overline{c_{i_0}} b_{i_k} u_{i_0,i_1} s_{i_1} u_{i_1,i_2} s_{i_2} \cdots u_{i_{k-1},i_k} s_{i_k}  ,
\end{eqnarray*}
and so we have
\begin{eqnarray*}
\E_\bU \big|\bc^* \bA^k \bb \big|^2
& = &  \sum_{ \overset{i_0,\ldots,i_{k}}{\underset{ j_0,\ldots,j_{k}}{}}}  \overline{c_{i_0}} c_{j_0} b_{i_{k}} \overline{b_{j_{k}}} s_{i_1} s_{j_1}  \ldots s_{i_{k}} s_{j_{k}} \E \left[ u_{i_0,i_1}u_{i_1,i_2}  \cdots  u_{i_{k-1},i_{k}}\overline{u}_{j_0,j_1}\overline{u}_{j_1,j_2} \ldots \overline{u}_{j_{k-1},j_{k}}  . \right] \\
\end{eqnarray*}
\indent Let $(i_1,\ldots,i_k),(i'_1,\ldots,i'_k),(j_1,\ldots,j_k)$ and $(j'_1,\ldots,j'_k)$ be $k$-tuples of intergers lower than $n$. By Proposition \ref{wg}, we know that
$$
\E \left[u_{i_1,j_1} \cdots u_{i_k,j_k} \ol{u}_{i'_1,j'_1} \cdots \ol{u}_{i'_k,j'_k} \right] \ \neq \ 0
$$ 
if and only if there are two permutations $\sigma$ and $\tau$ so that for all $p \in \left\{1,\ldots,k \right\}$, $i_{\sigma(p)} = i'_p$ and $j_{\tau(p)} = j'_p$. In our case, we know that for a $(i_0,\ldots,i_k)$ fixed, there will be no more than $(k+1)!$ tuples $(j_0,\ldots,j_k)$ leading to a non-zero expectation. By Proposition \ref{wg} again, we know that all these expectations are  $O\left(n^{-k} \right)$. So, one concludes with the following computation   
\begin{eqnarray*}
 \E_\bU |\bc^* \bA^k \bb|^2 
& \leq &  \sum_{\mu \in S_{k+1}} \sum_{ i_0,\ldots,i_{k}}  \left|\ol{c}_{i_0} c_{i_{\mu(0)}} b_{i_{k}} \ol{b}_{{i_{\mu(k)}}} \right| s_{i_1}^2  \ldots s_{i_{k}}^2 \times O \left( \frac{1}{n^k} \right) \\
& \leq & \sum_{\mu \in S_{k+1}} \sum_{ i_0,\ldots,i_{k}} \frac{1}{2} \left[|c_{i_0}|^2  |b_{i_{k}}|^2  + |c_{i_{\mu(0)}}|^2 |b_{i_{\mu(k)}}|^2 \right] s_{i_1}^2  \ldots s_{i_{k}}^2\times O \left( \frac{1}{n^k} \right) \\
& \leq & (k+1)! \sum_{ i_0,\ldots,i_{k}} |c_{i_0}|^2  |b_{i_{k}}|^2  s_{i_1}^2 \cdots s_{i_{k}}^2 O\left( n^{-k}\right) \\
& \leq &   \left(\frac{1}{n}\sum_{i=1}^n s_{i}^2\right)^{k-1} \left(  \sum_{j=1}^n |b_{j}|^2 s_{j}^2\right) \times O\left( \frac{1}{n}\right) 
\ = \ O\left( \frac{1}{n}\right). \\
\end{eqnarray*}

\subsection{Proof of Lemma \ref{hypSRT}}\label{ProofLemmahypSRT}
Lemma  \ref{hypSRT} is a direct consequence of the following lemma.

\begin{lemme}Let $\mu$ be a probability measure whose support is contained in an interval $[m, M]\subset]0,+\infty[$. Let $\mu^{-1}$ denote the push-forward of $\mu$ by the map $t\longmapsto 1/t$. 
Then for all $x\in \R$, $y>0$,  \begin{equation}\label{majStielInverse22713}|\im G_{\mu^{-1}}(x+i y)|\le \begin{cases} M&\textrm{ if $x\notin[1/(2M),2/m]$,}\\
\frac{8M^4}{m^2}\left|\im G_\mu\left(\frac{1}{x}+i\frac{m^2 y}{2}\right)\right|&\textrm{ otherwise.}
\end{cases}\end{equation}
\end{lemme}

\bpr Note that $$|\im G_{\mu^{-1}}(x+i y)|=\int\frac{y}{(x-1/t)^2+y^2}\mu(dt).$$
If $x\notin[1/(2M),2/m]$, then for all $t\in [m, M]$, $|x-1/t|\ge 1/(2M)$, and \eqref{majStielInverse22713} follows from the fact that for all $y>0$, we have $$\frac{y}{y^2+(1/(2M))^2}\le M.$$ If $x\in[1/(2M),2/m]$, then for all $t\in [m, M]$, $$\frac{1}{2M^2} \ \le \ \frac{x}{t} \ \le \ \frac{2}{m^2}$$ hence $$\left(x-\frac{1}{t}\right)^2+y^2 \ =  \ \frac{x^2}{t^2} \left(\left(\frac{1}{x} - t\right)^2+\left(\frac{yt}{x}\right)^2\right)\ \ge \ \frac{1}{4M^4} \left(\left(\frac{1}{x} - t\right)^2+\left(\frac{m^2y}{2}\right)^2\right)$$ and \eqref{majStielInverse22713}  follows directly.
\epr

\subsection{Proof of Lemma \ref{lemo(1)bis28114}} 
 First of all, as $\|\bB\|_{\op}$ and $\|\bC\|_{\op}$ are bounded (see Remark \ref{165130092014}) and for any unitary matrix $\bV$, $\bV \bB \eloi \bB$ and $\bC \bV^* \eloi \bC$, we know, by for example Theorem 2 of \cite{JIA}, that there is constant $C$ such that  with a probability tending to one, 
\begin{eqnarray} \label{14}
\forall n \geq 1, \ \ \max_{\displaystyle^{1 \leq i \leq n}_{1 \leq j \leq r}} |b_{i,j}| \ \leq \ C\sqrt{\frac{\log n}{n}} \ \ \text{ and } \ \ \max_{\displaystyle^{1 \leq i \leq n}_{1 \leq j \leq r}} |c_{j,i}| \ \leq \  C\sqrt{\frac{\log n}{n}}.
\end{eqnarray}

\subsubsection{Outline of the proof}
\indent If we expand the following expectation
$$
\Ec{\bc_{i_{1}}^* \bA^{k_1} \bb_{i_{1}}\cdots \bc_{i_{q}}^* \bA^{k_q} \bb_{i_{q}} \ol{\bc_{i'_{1}}^* \bA^{l_1} \bb_{i'_{1}} }\cdots \ol{\bc_{i'_{s}}^* \bA^{l_s} \bb_{i'_{s}} }}
$$
(where the expectation is with respect to the randomness of $\bU$), we get a summation of terms such as
\begin{equation}\label{8101421h}
\E \Big[ \underbrace{u_{t_{1,1} t_{1,2}} \cdots u_{t_{1,k_1} t_{1,k_1+1}}}_{k_1 \text{ factors}} \ccdots \underbrace{u_{t_{q,1} t_{q,2}} \cdots u_{t_{q,k_q} t_{q,k_q+1}}}_{k_q \text{ factors}} \ \underbrace{\ol{u}_{t'_{1,1} t'_{1,2}} \cdots \ol{u}_{t'_{1,l_1} t'_{1,l_1+1}}}_{l_1 \text{ factors}} \ccdots \underbrace{\ol{u}_{t'_{s,1} t'_{s,2}} \cdots \ol{u}_{t'_{s,l_s} t'_{s,l_s+1}}}_{l_s \text{ factors}} \Big].
\end{equation}

Our goal is to find out which of these terms will be negligible before the others. First, we know from   Proposition \ref{wg} that the expectation vanishes unless the set of the first indices (resp. second) of the $u_{ij}$'s is the same as the set of the first indices (resp. second) of the $\ol{u}_{ij}$'s. Secondly, each expectation  is computed thanks to the following formula   (see    Proposition \ref{wg}):
\begin{eqnarray} \label{wg1bis81014}
\E \left[ u_{i_1,j_1}  \cdots u_{i_k,j_k} \ol{ u_{ i'_1 ,j'_1 } }  \cdots \ol{u_{i'_k,j'_k}} \right] \ = \ \sum_{\sigma, \tau \in S_k} \delta_{i_1,i'_{\sigma(1)}} \ldots \delta_{i_k,i'_{\sigma(k)}} \delta_{j_1,j'_{\tau(1)}} \ldots \delta_{j_k,j'_{\tau(k)}} \wg( \tau \sigma^{-1}),
\end{eqnarray}
Then, by \eqre{wg2} of Proposition \ref{wg},  we know that the prevailing terms are the ones involving  $\wg(id)$, i.e. those allowing to match the $i$'s in the $u_{ij}$'s with the $i'$'s in the $\ol{u}_{i'j'}$'s thanks to    a permutation which also matches $j$'s in the $u_{ij}$'s with the $j'$'s in the $\ol{u}_{i'j'}$'s. To prove the second part of Lemma \ref{lemo(1)bis28114}, we shall characterize such terms among those of the type of \eqre{8101421h} and prove that the other ones are negligible. To prove the first part of the lemma, we shall prove that only negligible terms occur, i.e. that if $\left\{k_1,\ldots,k_q \right\}_m \neq \left\{l_1,\ldots,l_s \right\}_m$, then $\wg(id)$ can never occur in \eqre{8101421h}. Then, the third part of the lemma is only a straightforward summation following from  the first and second parts.

\subsubsection{Proof of $(2)$ of Lemma \ref{lemo(1)bis28114}: } Now, we reformulate the $(2)$ from Lemma \ref{lemo(1)bis28114} this way : let $k_1 > k_2 > \cdots > k_q$ be distinct positive integers and $m_1, \ldots, m_q$ positive integers, and let $\left(i_{\alp,\bet}\right)_{\displaystyle^{1 \leq \bet \leq q}_{1 \leq \alp \leq m_\bet}}$ and $\left(i'_{\alp,\bet}\right)_{\displaystyle^{1 \leq \bet \leq q}_{1 \leq \alp \leq m_\bet}}$ be some integers of $\{1,\ldots,r\}$. Our goal is to prove that
\beq
 &&\E \Big[ \sqrt{n}{\bc^*_{i_{1,1}}} \bA^{k_1} \bb_{i_{1,1}}\ol{\sqrt{n}{\bc^*_{i'_{1,1}}} \bA^{k_1} \bb_{i'_{1,1}}}\cdots \sqrt{n}\bc^*_{i_{m_1,1}}\bA^{k_1} \bb_{i_{m_1,1}} \ol{\sqrt{n}{\bc^*_{i'_{m_1,1}}} \bA^{k_1} \bb_{i'_{m_1,1}}}\ti\qquad\qquad\qquad\qquad\qquad \\
&&\qquad\qquad\qquad\qquad\qquad  \sqrt{n}{\bc^*_{i_{1,2}}} \bA^{k_2} \bb_{i_{1,2}}\ol{{\sqrt{n}\bc^*_{i'_{1,2}}}\bA^{k_2} \bb_{i'_{1,2}}} \cdots\cdots \sqrt{n}{\bc^*_{i_{m_q,q}}} \bA^{k_q} \bb_{i_{m_q,q}} \ol{\sqrt{n}{\bc^*_{i'_{m_q,q}}} \bA^{k_q} \bb_{i'_{m_q,q}}} \Big] \\
& = & b^{2\sum k_i m_i} \times \prod_{t=1}^q \left[\sum_{\mu_t \in S_{m_t}} \prod_{s=1}^{m_t}\left({\bb^*_{i'_{\mu_t(s),t}}} \bb_{i_{s,t}}  \cdot  {\bc^*_{i_{s,t}}} \bc_{i'_{\mu_t(s),t}}\right) \right] + o(1).\\ 
\eeq
\indent  We will denote the coordinate of $\bb_{i_{\alp,\beta}}$ : $\left(b^{i_{\alp,\beta}}_{t}\right)_{1 \leq t \leq n}$. We write 
\beqy \la{094317398456}
\bc_{i_{\al,\bet}}^* \bA^{k_1} \bb_{i_{\al,\bet}} \ = \ \sum_{1 \leq t_{0},\ldots,t_{k_1}\leq n} \ol{c_{t_0}^{i_{\al,\bet}}}u_{t_0,t_1} s_{t_1} u_{t_1,t_2}s_{t_2} \cdots u_{t_{k_1-1},t_{k_1}} s_{t_{k_1}} b_{t_{k_1}}^{i_{\al,\bet}} .
\eeqy

In order to simply notation, we shall use bold letters to designate tuples of consecutive indices. For example, we set $\bt_{i,j} : =( t_{0,i,j},t_{1,i,j},\ld,t_{k_j,i,j})$ and  write 
\beqy \label{notation1026}
u_{\bt_{i,j}} \ : = \ u_{t_{0,i,j},t_{1,i,j}} \cdots u_{t_{k_j-1,i,j},t_{k_j,i,j}} &\quad ; \quad & s_{\bt_{i,j}} \ := \ s_{t_{1,i,j}}\cdots s_{t_{k_j,i,j}},   
\eeqy

\noindent so if we expand the whole expectation in \eqref{094317398456} with respect to the randomness of $\bU$, we get terms as

%

\begin{eqnarray*} 
&  & \E \, u_{\bt_{1,1}} \cdots u_{\bt_{m_1,1}} \ \ol{u}_{\bt'_{1,1}}   \cdots \ol{u}_{\bt'_{m_1,1}} \ \ccdots \ u_{\bt_{1,q}}  \cdots u_{\bt_{m_q,q}} \ \ol{u}_{\bt'_{1,q}}   \cdots \ol{u}_{\bt'_{m_q,q}}   \\
& = & \E \prod_{1 \leq c \leq q} \ \prod_{1 \leq b \leq m_c} \ \prod_{0 \leq a \leq k_c-1} u_{t_{a,b,c}, t_{a+1,b,c}} \ol{u}_{t'_{a,b,c}, t'_{a+1,b,c}} , 
\end{eqnarray*}
and by Proposition \ref{wg}, for this expectation to be non-zero, the set of the first indices (resp. second) of the $u_{i,j}$'s must be the same than the set of the first indices (resp. second) of the $\ol{u}_{i,j}$'s, which can be expressed by the following  equalities of multisets\footnote{Recall that the notion of \emph{multiset} has been defined before Lemma \ref{lemo(1)bis28114}:  a multiset is roughly a collection of elements with possible repetitions (as in tuples) but where the order of appearance is insignificant (contrarily to tuples). For example,
$$
\left\{ 1,2,2,3 \right\}_m \eg \left\{ 3,2,1,2 \right\}_m \ \neq \ \left\{1,2,3 \right\}_m.
$$}
\beqy \label{ast1}
\begin{array}{c} \big\{ t_{a,b,c}\ste  \ 1 \leq c \leq q, \ 1 \leq b \leq m_c, \ 0 \leq a \leq k_c-1 \big\}_m \ =   \qquad\qquad\qquad\qquad  \\
 \qquad\qquad\qquad\qquad \big\{ t'_{a,b,c}\ste  \ 1 \leq c \leq q, \ 1 \leq b \leq m_c, \ 0 \leq a \leq k_c-1 \big\}_m  \\
 \end{array}
\eeqy
\beqy \label{ast2}
\begin{array}{c} \big\{ t_{a,b,c}\ste  \ 1 \leq c \leq q, \ 1 \leq b \leq m_c, \ 1 \leq a \leq k_c \big\}_m \ =     \qquad\qquad\qquad\qquad \\  
\qquad\qquad\qquad\qquad\quad  \big\{ t'_{a,b,c}\ste  \ 1 \leq c \leq q, \ 1 \leq b \leq m_c, \ 1 \leq a \leq k_c \big\}_m 
\end{array} 
\eeqy

 For now on, we denote by $N = \sum m_i k_i$ and by $\rho = \sum m_i$.

To make a better use of these multisets equalities, we shall need to reason on the $t_{a,b,c}$'s which are pairwise distincts and so, in a first place, we prove that the summation over the non pairwise indices is negligible. \\

To do so, we deduce first from \eqref{ast1} and \eqref{ast2} that 
for any fixed collection of indices $(t_{a,b,c})$, there is only a $O(1)$ choices of collection of indices $(t'_{a,b,c})$ leading to a non vanishing expectation. Then, noticing that 
 \begin{eqnarray*}
\card \left\{t_{a,b,c} \ste 1\le c\le q,\, 1\le b\le m_c,\, 0\le a\le  k_c\right\} & = & O \left( n^{N+\rho-1}\right),
\end{eqnarray*}
(where $N = \sum m_i k_i$ and $\rho = \sum m_i$), we know that the summation contains a $O \!\left( n^{N+\rho-1}\right)$ of terms. At last, we use the fact that the expectation over the $u_{i,j}$'s and the $\ol{u}_{i,j}$'s is at most $O\! \left({n^{-N}}\right)$, that $\sup_i s_i < M$ and that $|b_i|$ and $|c_i| = O\! \left(\sqrt{\frac{ \log n}{n}} \right)$ (recall \eqref{14}), to claim that each term of the sum is at most a $\OO{(\log n)^{2\rho} n^{-N-\rho}}$ (one should not forget that each $\bc \bA^k \bb$ is multiply by $\sqrt{n}$). 
We conclude that the summation over the non pairwise distinct indices is a $O\! \left( \frac{(\log n )^{2\rho}}{n} \right)$.\\

\indent For now on, we consider only the pairwise distinct indices so that \eqref{ast1} and \eqref{ast2} can be seen as set equalities (instead of multiset). Also, if one sees the sets as $N$-tuple, the equalities \eqref{ast1} and \eqref{ast2} means that there exists two permutations $\sigma_1$ and $\sigma_2$ in $S_N$ so that for all $1 \leq c \leq q$, $ 1 \leq b \leq m_q$, $ 0 \leq a \leq k_q-1$ (resp. $1\leq a \leq k_q$), we have $t'_{a,b,c} = t_{\sigma_1(a,b,c)}$ (resp. $t'_{a,b,c} = t_{\sigma_2(a,b,c)}$). 
\begin{rem}
\indent The notation $t_{\sigma_1(a,b,c)}$ is a little improper: the set $$\left\{(a,b,c)\ste \ 1 \leq c \leq q, \ 1 \leq b \leq m_c, \ 0 \leq a \leq k_c-1\right\}$$ is identified with $\{1,\ldots,N\}$ thanks to the colexicographical order (where $N = \sum m_i k_i$).
\end{rem}

\indent Thanks to the Proposition \ref{wg} and the Remark \ref{moebius}, we know that the expectation of the $u_{t_{a,b,c}}$'s and the $\ol{u}_{t'_{a,b,c}}$'s is equal to $\wg(\sigma_1 \circ \sigma_2^{-1})$ and so, we know that we can neglect all of these with $\sigma_1 \neq \sigma_2$. For now on, we suppose $\sigma_1 = \sigma_2$. \\
\indent One needs to understand that the sets  $\left\{  t_{a,b,c}\ste  \ 1 \leq c \leq q, \ 1 \leq b \leq m_c,\ 0\le a< k_c \right\}$ and \\
 $\left\{  t_{a,b,c}\ste  \ 1 \leq c \leq q, \ 1 \leq b \leq m_c,\ 1\le a\le  k_c \right\}$ are very similar except for the shift for the first index. Due to this likeness and the fact that they are both mapped onto the $t'_{a,b,c}$'s in the   same way (i.e. $\sigma_1 = \sigma_2$), we prove that the choice of $\sigma_1$ is very specific : \\ \\
~--~ First, using the distinctness of the indices, it easy to see that the equalities \eqref{ast1} and \eqref{ast2} lead us to these new equalities of sets 
\beqy \label{star1}
 \left\{ t_{0,b,c} \ste \ 1 \leq c \leq q, \ 1 \leq b \leq m_c \right\} \ = \  \left\{ t'_{0,b,c}\ste \ 1 \leq c \leq q, \ 1 \leq b \leq m_c \right\}, 
\eeqy
and
\beqy \label{star2}
 \left\{ t_{k_c,b,c}\ste \ 1 \leq c \leq q, \ 1 \leq b \leq m_c \right\} \ = \  \left\{ t'_{k_c,b,c}\ste \ 1 \leq c \leq q, \ 1 \leq b \leq m_c \right\} ,
\eeqy
~--~ According to the equality \eqref{star1}, we know that $\left\{ t_{0,b,c}, \ 1 \leq c \leq q, \ 1 \leq b \leq m_c \right\}$ is an invariant set of $\sigma_1$. Indeed, we know that 
$$
\left\{ t_{\sigma_1(0,b,c)}\ste \ 1 \leq c \leq q, \ 1 \leq b \leq m_c \right\} \ \subset \ \left\{ t'_{a,b,c}\ste \ 1 \leq c \leq q, \ 1 \leq b \leq m_c, \ 0 \leq a \leq k_c-1 \right\},
$$
and with the condition \eqref{star1}, to avoid non pairwise distinct indices, we must have
$$
\left\{ t_{\sigma_1(0,b,c)}\ste \ 1 \leq c \leq q, \ 1 \leq b \leq m_c \right\} \ = \ \left\{ t'_{0,b,c}\ste \ 1 \leq c \leq q, \ 1 \leq b \leq m_c \right\},
$$
and so, we deduce that
$$
\left\{ t_{0,b,c}\ste \ 1 \leq c \leq q, \ 1 \leq b \leq m_c \right\} \ = \ \left\{ t_{\sigma_1(0,b,c)}\ste \ 1 \leq c \leq q, \ 1 \leq b \leq m_c \right\}.
$$
~--~ As $\sigma_1 = \sigma_2$,  $\sigma_2$ permutes $\left\{ t_{1,b,c}\ste \ 1 \leq c \leq q, \ 1 \leq b \leq m_c \right\}$ in the same way (actually, the sets 
$
\left\{ (a,b,c)\ste \ 1 \leq c \leq q, \ 1 \leq b \leq m_c, \ 0 \leq a \leq k_c-1 \right\}$ and $\left\{ (a,b,c)\ste \ 1 \leq c \leq p, \ 1 \leq b \leq m_c, \ 1 \leq a \leq k_c \right\}
$
are indentified to the same set (with cardinality $N$) thanks to the colexicographical order, and so, the action of $\sigma_1$ and $\sigma_2$ must be seen on this common set).\\ \\
~--~ As each element of $\left\{ t_{1,b,c}\ste \ 1 \leq c \leq q, \ 1 \leq b \leq m_c\right\}$ has only one corresponding $t'_{d,e,f}$ (indeed by \eqref{ast1} and \eqref{ast2} and as the $t$'s and the  $t'$'s are pairwise   distinct, to each $t$ corresponds a  unique $t'$), we deduce that $\sigma_1$ permutes $\left\{ t_{1,b,c}\ste \ 1 \leq c \leq q, \ 1 \leq b \leq m_c \right\}$  in  the same way (indeed, it allows to claim that   
$$
\left( t_{\sigma_1(1,b,c)}\ste \ 1 \leq c \leq q, \ 1 \leq b \leq m_c \right) \ = \ \left( t_{\sigma_2(1,b,c)}\ste \ 1 \leq c \leq q, \ 1 \leq b \leq m_c \right)
$$
as $N$-tuples.\\
\\
~--~ As $\sigma_1 = \sigma_2$, we know that $\sigma_2$ permutes $\left\{ t_{2,b,c}\ste \ 1 \leq c \leq q, \ 1 \leq b \leq m_c \right\}$ in  the same way, and so on until one shows that $\sigma_2$ permutes $\left\{ t_{k_q,b,c}\ste \ 1 \leq c \leq q, \ 1 \leq b \leq m_c \right\}$  in  the same way.\\  \\
~--~ However, according to \eqref{star2}, we know that $\left\{ t_{k_c,b,c}\ste \ 1 \leq c \leq q, \ 1 \leq b \leq m_c \right\} $ is an invariant set of $\sigma_2$. \\

\indent Therefore, as $\left\{ t_{k_q,b,c}\ste \ 1 \leq c \leq q, \ 1 \leq b \leq m_c \right\}$  and $\left\{ t_{k_c,b,c}\ste \ 1 \leq c \leq q, \ 1 \leq b \leq m_c \right\} $ are invariant sets by $\sigma_2$, we know that
$$
\big\{ t_{k_q,b,c}\ste \ 1 \leq c \leq q, \ 1 \leq b \leq m_c \big\} \cap \big\{ t_{k_c,b,c}\ste \ 1 \leq c \leq q, \ 1 \leq b \leq m_c \big\} \ = \ \left\{ t_{k_q,b,q}\ste \ 1 \leq b \leq m_q \right\}
$$ 
is also an invariant set of $\sigma_2$ and we deduce that $\sigma_1$ permutes in the same way every set of the form $\left\{ t_{l,b,q}, \ 1 \leq b \leq m_q \right\}$ for $l \in \left\{0,k_q-1\right\}$. And so, we rewrite the equalities \eqref{star1} and \eqref{star2} 
\beqy \label{starstar1}
 \left\{ t_{0,b,c}\ste \ 1 \leq c \leq q-1, \ 1 \leq b \leq m_c \right\} \ = \  \left\{ t'_{0,b,c}\ste \ 1 \leq c \leq q-1, \ 1 \leq b \leq m_c \right\},  
\eeqy
and
\beqy \label{starstar2}
 \left\{ t_{k_c,b,c}\ste \ 1 \leq c \leq q-1, \ 1 \leq b \leq m_c \right\} \ = \  \left\{ t'_{k_c,b,c}\ste \ 1 \leq c \leq q-1, \ 1 \leq b \leq m_c \right\}, 
\eeqy
and   one can make an induction on $q$ to show that there exist $\mu_1 \in S_{m_1}$, $\mu_2 \in S_{m_2}, \ldots, \mu_q \in S_{m_q}$ such that for all $1 \leq c \leq q,$ $1 \leq b \leq m_c$, $1 \leq a \leq k_c$, we have
$$
t'_{a,b,c} \ = \ t_{\sigma_1(a,b,c)} \ = \ t_{a,\mu_c(b),c} 
$$ 
and so, to sum up, we deduce that the non negligible terms that we get when we expand the whole expectation are terms such as
$$
\E \prod_{1 \leq c \leq q} \ \prod_{1 \leq b \leq m_c} \ \prod_{0 \leq a \leq k_c-1} u_{t_{a,b,c}, t_{a+1,b,c}} \ol{u}_{t_{a,\mu_c(b),c}, t_{a+1,\mu_c(b),c}}
$$
where for all $c$, $\mu_c$ belongs to $S_{m_c}$ and so, one can easily deduce that
\beq
 &&\E \Big[ \sqrt{n}{\bc^*_{i_{1,1}}} \bA^{k_1} \bb_{i_{1,1}}\ol{\sqrt{n}{\bc^*_{i'_{1,1}}} \bA^{k_1} \bb_{i'_{1,1}}}\cdots \sqrt{n}\bc^*_{i_{m_1,1}}\bA^{k_1} \bb_{i_{m_1,1}} \ol{\sqrt{n}{\bc^*_{i'_{m_1,1}}} \bA^{k_1} \bb_{i'_{m_1,1}}}\ti  \\
&& \qquad\qquad \qquad  \sqrt{n}{\bc^*_{i_{1,2}}} \bA^{k_2} \bb_{i_{1,2}}\ol{{\sqrt{n}\bc^*_{i'_{1,2}}}\bA^{k_2} \bb_{i'_{1,2}}} \cdots\cdots \sqrt{n}{\bc^*_{i_{m_q,q}}} \bA^{k_q} \bb_{i_{m_q,q}} \ol{\sqrt{n}{\bc^*_{i'_{m_q,q}}} \bA^{k_q} \bb_{i'_{m_q,q}}} \Big] \\
&  & \ = \ \left( \frac{1}{n} \sum_{i=1}^n s_i^2\right)^{N-\rho } \cdot \prod^q_{u=1} \left[\sum_{\mu_u \in S_{m_u}} \prod_{s=1}^{m_u}\left( \sum_{\displaystyle^{1 \leq t_{0,s,u} \leq n}_{1 \leq t_{k_u,s,u} \leq n}} \ol{b^{i'_{\mu_u(s),u}}_{t_{k_u,s,u}}} b^{i_{s,u}}_{t_{k_u,s,u}} s_{t_{k_u,s,u}}^2  \ol{c^{i_{s,u}}_{t_{0,s,u}}} c^{i'_{\mu_u(s),u}}_{t_{0,s,u}}  \right) \right]+ o(1)
\\
&  & \ = \ b^{2N} \times \prod_{u=1}^q \left[\sum_{\mu_u \in S_{m_u}} \prod_{s=1}^{m_u}\left({\bb^*_{i'_{\mu_u(s),u}}} \bb_{i_{s,u}}  \cdot  {\bc^*_{i_{s,u}}} \bc_{i'_{\mu_u(s),u}}\right) \right] + o(1),\\ 
\eeq
and   we can conclude. \\
\begin{rem}
\indent We used the fact that 
\begin{eqnarray}
\frac{1}{n} \sum_{i=1}^n s_i^2 & = & b^2 +o(1)  \label{eq10022014}\\
\sum_{j=1}^n \ol{b^{i_{\al,\bet}}_j} b^{i_{\gamma,\delta}}_j s_j^2 & = & b^2 \bb^*_{i_{\al,\bet}} \bb_{i_{\gamma,\delta}} + o(1). \label{eq10022014bis}
\end{eqnarray}
The relation \eqref{eq10022014} is obvious and the \eqref{eq10022014bis} can be proved using the fact $\bP$ is invariant, in law, by conjugation by any unitary matrix (we explained at Section \ref{8813226214} that we can add this hypothesis). 
\end{rem}
\subsubsection{Proof of $(1)$ of Lemma \ref{lemo(1)bis28114}: }
 The proof of $(1)$ goes along the same lines as the previous proof. Our goal is to show that 
\beq
\E \left[\sqrt{n}\bc_{i_{1}}^* \bA^{k_1} \bb_{i_{1}}\cdots \sqrt{n}\bc_{i_{q}}^* \bA^{k_q} \bb_{i_{q}} \ol{\sqrt{n}\bc_{i'_{1}}^* \bA^{l_1} \bb_{i'_{1}} }\cdots \ol{\sqrt{n}\bc_{i'_{s}}^* \bA^{l_s} \bb_{i'_{s}} }\right] &=& o \left( 1 \right).
\eeq
\indent At first, one can notice that if $\sum k_i \neq \sum l_j$, the expectation is equal to zero. We assume now  that $\sum k_i = \sum l_j$, and let $N$ denote the common value. Then, we distinguish two cases. \\ \\
\noindent $\bullet$ \underline{First case : $q=s$} \\
\indent Then we can also focus on  the ``pairwise distinct indices'' summation, by  similar argument as in the previous proof. We suppose that there exists $j$ such that $k_j \neq l_j$ (otherwise, one should read the previous proof). Our goal is to show that there is no expectation equal to $\wg(id)$ (which means that we cannot have $\sigma_1 = \sigma_2$) in that case and so we shall conclude that 
\beq
\E \left[\sqrt{n}\bc_{i_{1}}^* \bA^{k_1} \bb_{i_{1}}\cdots \sqrt{n}\bc_{i_{q}}^* \bA^{k_q} \bb_{i_{q}} \ol{\sqrt{n}\bc_{i'_{1}}^* \bA^{l_1} \bb_{i'_{1}} }\cdots \ol{\sqrt{n}\bc_{i'_{q}}^* \bA^{l_q} \bb_{i'_{q}} }\right] & = & O \left( \frac{1}{n}\right).
\eeq
\indent Let us gather the $k_i$'s which are equal and in order to simply the expressions, we shall use notations in the same spirit than \eqref{notation1026}
\beqy
\begin{array}{rcl} \sqrt{n} (\bc^* \bA^{k_\al} \bb)_{\bi_\al} & : = & \sqrt{n}{\bc^*_{i_{1,\al}}} \bA^{k_\al} \bb_{i_{1,\al}}\cdots \sqrt{n}\bc^*_{i_{m_\al,\al}}\bA^{k_\al} \bb_{i_{m_\al,\al}}, \\
\vspace{2mm}\sqrt{n} (\bc^* \bA^{\ell_\bet} \bb)_{\bi'_\bet} & : = & \sqrt{n}{\bc^*_{i'_{1,\bet}}} \bA^{\ell_\bet} \bb_{i'_{1,\bet}}\cdots \sqrt{n}\bc^*_{i'_{n_\bet,\bet}}\bA^{\ell_\bet} \bb_{i'_{n_\bet,\bet}}, 
\end{array}
\eeqy
 so that we rewrite our expectation
\beq
 &&\E \Big[ \sqrt{n} (\bc^* \bA^{k_1} \bb)_{\bi_1} \cdots \sqrt{n} (\bc^* \bA^{k_r} \bb)_{\bi_r}  \ol{\sqrt{n} (\bc^* \bA^{\ell_1} \bb)_{\bi'_1}} \cdots \ol{\sqrt{n} (\bc^* \bA^{\ell_s} \bb)_{\bi'_s}} \Big]  
\eeq
with $\sum_{i=1}^r m_i k_i = \sum_{j=1}^s n_j l_j$ and $k_1 > \cdots > k_r$ and $l_1 > \cdots > l_s$. Without loss of generality, we shall assume that $(k_r,m_r) \neq (l_s,n_s)$ (indeed, otherwise, we can start the   induction from the previous proof until we find an integer $x$ such that $(k_{r-x},m_{r-x}) \neq (l_{s-x},n_{s-x})$ to show that the expectation is equal to
\beq
&&\Ec{\sqrt{n} (\bc^* \bA^{k_1} \bb)_{\bi_1} \cdots \sqrt{n} (\bc^* \bA^{k_{r-x}} \bb)_{\bi_{r-x}}  \ol{\sqrt{n} (\bc^* \bA^{\ell_1} \bb)_{\bi'_1}} \cdots \ol{\sqrt{n} (\bc^* \bA^{\ell_{s-x}} \bb)_{\bi'_{s-x}}} }\\
&\times& \prod_{t=r-x+1}^r \ \sum_{\mu_t \in S_{m_t}} \ \prod_{s=1}^{m_t} \ \Ec{\sqrt{n}\bc^*_{i_{s,t}} \bA^{k_t} \bb_{i_{s,t}}\ol{\sqrt{n}\bc^*_{i'_{\mu_t(s),t}} \bA^{k_t} \bb_{i'_{\mu_t(s),t}}}} + \oo{1},
\eeq
and the following of the proof is the same). We shall also assume that $k_r \leq l_s$. \\
\indent According to Proposition \ref{wg}, we have the following equalities
\beqy \label{astbis1}
\begin{array}{c} \big\{ t_{a,b,c}\ste \ 1 \leq c \leq r, \ 1 \leq b \leq m_c, \ 0 \leq a \leq k_c-1 \big\} \qquad\quad\qquad\qquad \\ \qquad\quad\qquad\qquad \  = \  \big\{ t'_{a,b,c}\ste \ 1 \leq c \leq s, \ 1 \leq b \leq n_c, \ 0 \leq a \leq l_c-1 \big\}, 
\end{array}
\eeqy
and
\beqy \label{astbis2}
 \begin{array}{c} \big\{ t_{a,b,c}\ste \ 1 \leq c \leq r, \ 1 \leq b \leq m_c, \ 1 \leq a \leq k_c \big\} \qquad\quad\qquad\qquad \\\quad\qquad\qquad\qquad \ = \  \big\{ t'_{a,b,c}\ste \ 1 \leq c \leq s, \ 1 \leq b \leq n_c, \ 1 \leq a \leq l_c \big\}, 
 \end{array}
\eeqy
and let $\sigma_1$ and $\sigma_2$ the two permutations describing these equalities. Let us prove by contradiction that $\sigma_1 \neq \sigma_2$ and so let us suppose that $\sigma_1 = \sigma_2$. As we consider only pairwise distinct indices, we have also
\beqy \label{starbis1}
 \left\{ t_{0,b,c}\ste \ 1 \leq c \leq r, \ 1 \leq b \leq m_c \right\} \ = \  \left\{ t'_{0,b,c}\ste \ 1 \leq c \leq s, \ 1 \leq b \leq n_c \right\}, 
\eeqy
and
\beqy \label{starbis2}
 \left\{ t_{k_c,b,c}\ste \ 1 \leq c \leq r, \ 1 \leq b \leq m_c \right\} \ = \  \left\{ t'_{l_c,b,c}\ste \ 1 \leq c \leq s, \ 1 \leq b \leq n_c \right\}, 
\eeqy
\indent According to the fact that $\sigma_1 = \sigma_2$ and \eqre{starbis1}, we can deduce that 
\beqy \label{starstar}
 \left\{ t_{k_r,b,c}\ste \ 1 \leq c \leq r, \ 1 \leq b \leq m_c \right\} \ = \  \left\{ t'_{k_r,b,c}\ste \ 1 \leq c \leq s, \ 1 \leq b \leq n_c \right\}, 
\eeqy
and here comes the contradiction. Indeed, if $k_r < l_s$, then    
$$
\left\{ t'_{l_c,b,c}\ste \ 1 \leq c \leq s, \ 1 \leq b \leq n_c \right\} \cap \left\{ t'_{k_r,b,c}\ste \ 1 \leq c \leq s, \ 1 \leq b \leq n_c \right\} = \emptyset,
$$
otherwise, $k_r=l_s$ (which means $m_r \neq n_s$), let us suppose $m_r < n_s$, so that  
$$
\card\,  \left\{ t'_{l_c,b,c}\ste \ 1 \leq c \leq s, \ 1 \leq b \leq n_c \right\} \cap \left\{ t'_{k_r,b,c}\ste \ 1 \leq c \leq s, \ 1 \leq b \leq n_c \right\}  \ = \ n_s
$$
however,
$$
\card \, \left\{ t_{k_c,b,c}\ste \ 1 \leq c \leq r, \ 1 \leq b \leq m_c \right\} \cap \left\{ t_{k_r,b,c}\ste \ 1 \leq c \leq r, \ 1 \leq b \leq m_r \right\}  \ = \ m_r,
$$
which is, according to \eqref{starbis2} and \eqref{starstar}, impossible.\\ \\

\noindent $\bullet$ \underline{Second case : $q \neq s$} \\
\indent Without loss of generality, we suppose that   $q>s$. We cannot consider here the pairwise distinct indices simply because the cardinal of the $t_{i,j}$'s is different than the one of the $t'_{i,j}$'s. \\
  Expanding  the product
$$
\E \left[\sqrt{n}\bc_{i_{1}}^* \bA^{k_1} \bb_{i_{1}}\cdots \sqrt{n}\bc_{i_{q}}^* \bA^{k_q} \bb_{i_{q}} \ol{\sqrt{n}\bc_{i'_{1}}^* \bA^{l_1} \bb_{i'_{1}} }\cdots \ol{\sqrt{n}\bc_{i'_{s}}^* \bA^{l_s} \bb_{i'_{s}} }\right],
$$  
we get terms such as
\begin{eqnarray*}
&&  \E_\bU \left[u_{t_{0,1},t_{1,1}}\cdots u_{t_{k_1-1,1},t_{k_1,1}} u_{t_{0,2},t_{1,2}} \cdots u_{t_{k_q-1,q},t_{k_q,q}} \ol{u_{t'_{0,1},t'_{1,1}}\cdots u_{t'_{l_1-1,1},t'_{l_1,1}} u_{t'_{0,2},t'_{1,2}} \cdots u_{t'_{l_s-1,s},t'_{l_s,s}}}\right]
\end{eqnarray*}

\indent According to Proposition \ref{p1}, for the expectation to be non zero, one needs to have the equality of sets
\begin{eqnarray*}
\{t_{0,1},t_{1,1},\ldots,t_{k_1-1,1},t_{0,2},\ldots,t_{k_2-1,2},\ldots,t_{k_q-1,q}\} & = & \{t'_{0,1},t'_{1,1},\ldots,t'_{l_1-1,1},t'_{0,2},\ldots,t'_{l_2-1,2},\ldots,t'_{l_s-1,s}\}, \\
\{t_{1,1},t_{2,1},\ldots,t_{k_1,1},t_{1,2},\ldots,t_{k_2,2},t_{1,3},\ldots,t_{k_q,q}\} & = & \{t'_{1,1},t'_{2,1},\ldots,t'_{l_1,1},t'_{1,2},\ldots,t'_{l_2,2},t'_{1,3},\ldots,t'_{l_s,s}\}. 
\end{eqnarray*}

Set $\mathcal{A}:=\left\{ t_{a,b} , \ 1 \leq b \leq r, \ 0 \leq a \leq k_b\right\}$ and $\mathcal{B} = \left\{ t'_{a,b} , \ 1 \leq b \leq s, \ 0 \leq a \leq l_b \right\}$. 
According to the previous inequalities, to \eqref{14}  and to the fact that all the expectations are $O \left( n^{-N}\right)$, we write
\begin{eqnarray*}
&&\left|\E \left[\sqrt{n}\bc_{i_{1}}^* \bA^{k_1} \bb_{i_{1}}\cdots \sqrt{n}\bc_{i_{q}}^* \bA^{k_q} \bb_{i_{q}} \ol{\sqrt{n}\bc_{i'_{1}}^* \bA^{l_1} \bb_{i'_{1}} }\cdots \ol{\sqrt{n}\bc_{i'_{s}}^* \bA^{l_s} \bb_{i'_{s}} }\right]\right| \\
& \leq &C^{4q} \frac{\log^{2q} (n)}{ n^{2q}} \cdot n^{\frac{q+s}{2}}  \cdot \sum_{\substack{  t_{0,1},\ldots,t_{k_1,1}\\  \vdots\\  t_{0, q},\ldots,t_{k_q,q}}} \sum_{ \mu \in \mathcal{B}^{ \mathcal{A}}
} \prod_{a=1}^q  s_{t_{1,a}} s_{\mu(t_{1,a})} \cdots s_{\mu(t_{k_a,a})} \times O \left( n^{-N} \right) \\
& \leq & O\left(\log^{2q}(n) \cdot n^{\frac{q+s}{2}-N-2q}\right) \sum_{\substack{\mu \in \mathcal{B}^{\mathcal{A}}  \\ t_{0,1},\ldots,t_{k_1,1}\\  \vdots\\  t_{0, q},\ldots,t_{k_q,q}}} \frac{1}{2}\left[\prod_{a=1}^q s_{t_{1,a}}^2 \cdots s^2_{t_{k_a,a}} + \prod_{a=1}^q    s^2_{\mu(t_{1,a})} \cdots s^2_{\mu(t_{k_a,a})} \right] 
\end{eqnarray*}
\indent On the one hand, 
\begin{eqnarray*}
\sum_{\substack{\mu \in \mathcal{B}^{\mathcal{A}}  \\ t_{0,1},\ldots,t_{k_1,1}\\  \vdots\\  t_{0, q},\ldots,t_{k_q,q}}} \prod_{a=1}^q s_{t_{1,a}}^2 \cdots s^2_{t_{k_a,a}} 
& \leq & \card(\mathcal{B}^{\mathcal A})\times \card\left(\{1,\ldots,n\}^{N+q}\right) \times M^{2N}  
\ = \ O \left( n^{N+q}\right).
\end{eqnarray*}
\indent On the other, 
\begin{eqnarray*}
\sum_{\substack{\mu \in \mathcal{B}^{\mathcal{A}}  \\ t_{0,1},\ldots,t_{k_1,1}\\  \vdots\\  t_{0, q},\ldots,t_{k_q,q}}} \prod_{a=1}^q    s^2_{\mu(t_{1,a})} \cdots s^2_{\mu(t_{k_a,a})} 
& = & O \left( 1 \right)\times \sum_{\substack{ t'_{0,1},\ldots,t'_{l_1,1} \\ \vdots \\ t'_{0, s},\ldots,t'_{l_s,s}}} \prod_{a=1}^s s^2_{t'_{l_1,a}}\cdots s^2_{t'_{l_a,a}} 
\ = \ O \left(n^{N+s} \right).
\end{eqnarray*}
\indent Indeed, for any fixed $J=\{t'_{0,1},t'_{1,1},\ldots,t'_{l_1,1},t'_{0,2},\ldots,t'_{l_2,2},\ldots,t'_{l_s,s}\}$, there are $  O (1)$ of $\mu$'s in $\mathcal{S}_{\mathcal{I} \to \mathcal{J}}$ and $I = \{t_{0,1},t_{1,1},\ldots,t_{k_1,1},t_{0,2},\ldots,t_{k_2,2},\ldots,t_{k_q,q}\}$ such that $\mu(I) = J$.\\
\indent Therefore, 
\begin{eqnarray*}
\left|\E_{\bU} \left[\sqrt{n}\bc_{i_{1}}^* \bA^{k_1} \bb_{i_{1}}\cdots \sqrt{n}\bc_{i_{q}}^* \bA^{k_q} \bb_{i_{q}} \ol{\sqrt{n}\bc_{i'_{1}}^* \bA^{l_1} \bb_{i'_{1}} }\cdots \ol{\sqrt{n}\bc_{i'_{s}}^* \bA^{l_s} \bb_{i'_{s}} }\right]\right| & = & O \left( \log^{2q}(n) \cdot n^{-\frac{q-s}{2}}\right),
\end{eqnarray*}  
and, since $q>s$, it is at least a $O \left( \frac{\log^{2q}(n)}{\sqrt{n}} \right)$.

\subsubsection{Proof of $(3)$ of Lemma \ref{lemo(1)bis28114}: } \indent If $\left\{k_1,\ldots,k_q \right\}_m \neq \left\{k'_1,\ldots,k'_q \right\}_m$, we know that it contributes to the $o(1)$. So we rewrite
\beq
&&\sum_{\displaystyle^{k_1,\ldots,k_q =1}_{k'_1,\ldots,k'_q =1}}^{k_0}  \E \left[ \prod_{\al=1}^q \sqrt{n}\bc_{i_{\al}}^* \frac{\bA^{k_\al}}{z_{i_1}^{k_\al+1}} \bb_{i_{\al}} \ol{\sqrt{n}\bc_{i'_{\al}}^* \frac{\bA^{k'_\al}}{z_{i'_\al}^{k'_\al+1}} \bb_{i'_{\al}}} \right]\\
&  = & \sum_{k_1,\ldots,k_q =1}^{k_0} \sum_{\underset{{\left\{k_1,\ldots,k_q \right\}_m = \left\{k'_1,\ldots,k'_q \right\}_m}}{k'_1,\ldots,k'_q =1}}^{k_0} \E \left[ \prod_{\al=1}^q \sqrt{n}\bc_{i_{\al}}^* \frac{\bA^{k_\al}}{z_{i_1}^{k_\al+1}} \bb_{i_{\al}} \ol{\sqrt{n}\bc_{i'_{\al}}^* \frac{\bA^{k'_\al}}{z_{i'_\al}^{k'_\al+1}} \bb_{i'_{\al}}} \right] + o \left( 1 \right).\\
\eeq
Then, we fixed $(k_1,\ldots,k_q)$, and let us calculate
\begin{eqnarray} \label{expect10022014}
&&\sum_{\underset{{\left\{k_1,\ldots,k_q \right\}_m = \left\{k'_1,\ldots,k'_q \right\}_m}}{k'_1,\ldots,k'_q =1}}^{k_0} \E \left[ \prod_{\al=1}^q \sqrt{n}\bc_{i_{\al}}^* \frac{\bA^{k_\al}}{z_{i_1}^{k_\al+1}} \bb_{i_{\al}} \ol{\sqrt{n}\bc_{i'_{\al}}^* \frac{\bA^{k'_\al}}{z_{i'_\al}^{k'_\al+1}} \bb_{i'_{\al}}} \right],
\end{eqnarray}
to do so, we will use the previous notations and write 
$$\left(\underbrace{k_1,\ldots,k_1}_{m_1}, \underbrace{k_2,\ldots,k_2}_{m_2}, \ldots,\underbrace{k_s,\ldots,k_s}_{m_s} \right)$$
and we shall show that \eqref{expect10022014} doesn't depend on the $m_i$'s but depends only on $q = \sum m_i k_i$. We rewrite the summation
\begin{align}\la{1131000125}
\sum_{\{k'_1,\ldots,k'_q\}_m = \{k_1,\ldots,k_s\}_m}\E \left[ \sqrt{n}\bc_{i_{1,1}}^* \frac{\bA^{{k}_1}}{z_{i_{1,1}}^{{k_1+1}}} \bb_{i_{1,1}}\ol{\sqrt{n}\bc_{i'_{1,1}}^* \frac{\bA^{k'_1}}{z_{i'_{1,1}}^{k'_1+1}} \bb_{i'_{1,1}}} \cdots \sqrt{n}\bc_{i_{m_1,1}}^* \frac{\bA^{{k}_1}}{z_{i_{m_1,1}}^{{k_1+1}}} \bb_{i_{m_1,1}}\ol{ \sqrt{n}\bc_{i'_{m_1,1}}^* \frac{\bA^{k'_{m_1}}}{z_{i'_{m_1,1}}^{k'_{m_1}+1}} \bb_{i'_{m_1,1}}} \right.& \nonumber\\
\sqrt{n}\bc_{i_{1,2}}^* \frac{\bA^{k_2}}{z_{i_{1,2}}^{k_2+1}} \bb_{i_{1,2}}\ol{\sqrt{n}\bc_{i'_{1,2}}^* \frac{\bA^{k'_{m_1+1}}}{z_{i'_{1,2}}^{k'_{m_1+1}+1}} \bb_{i'_{1,2}}} \cdots  \sqrt{n}\bc_{i_{m_2,2}}^* \frac{\bA^{k_2}}{z_{i_{m_2,2}}^{k_2+1}} \bb_{i_{m_2,2}}\ol{\sqrt{n}\bc_{i'_{m_2,2}}^* \frac{\bA^{k'_{m_1+m_2}}}{z_{i'_{m_2,2}}^{k'_{m_1+m_2}+1}} \bb_{i'_{m_2,2}}}&\\
\left. \cdots \sqrt{n}\bc_{i_{m_s,s}}^* \frac{\bA^{k_s}}{z_{i_{m_s,s}}^{k_s+1}} \bb_{i_{m_s,s}}\ol{\sqrt{n}\bc_{i'_{m_s,s}}^* \frac{\bA^{k'_q}}{z_{i'_{m_s,s}}^{k'_q+1}} \bb_{i'_{m_s,s}}} \right] . \nonumber &
\end{align}
\indent We gather the $k_i$'s which are equal, so we rewrite the summation thanks to permutations of the set \\
$\mathcal{I} = \left\{ (\alp,\beta), \ 1 \leq \beta \leq s, \  1 \leq \alp \leq m_s \right\}$ :
$$ \sum_{\mu \in S_{\mathcal{I}}}\E \left[ \sqrt{n}\bc_{i_{1,1}}^* \frac{\bA^{{k}_1}}{z_{i_{1,1}}^{{k_1+1}}} \bb_{i_{1,1}}\ol{\sqrt{n}\bc_{i'_{\mu(1,1)}}^* \frac{\bA^{k_1}}{z_{i'_{\mu(1,1)}}^{k_1+1}} \bb_{i'_{\mu(1,1)}}} \cdots \sqrt{n}\bc_{i_{m_1,1}}^* \frac{\bA^{{k}_1}}{z_{i_{m_1,1}}^{{k_1+1}}} \bb_{i_{m_1,1}}\ol{ \sqrt{n}\bc_{i'_{\mu(m_1,1)}}^* \frac{\bA^{k_{1}}}{z_{i'_{\mu(m_1,1)}}^{k_{1}+1}} \bb_{i'_{\mu(m_1,1)}}} \right.$$
$$\qquad\qquad\qquad \left.  \sqrt{n}\bc_{i_{1,2}}^* \frac{\bA^{k_2}}{z_{i_{1,2}}^{k_2+1}} \bb_{i_{1,2}}\ol{\sqrt{n}\bc_{i'_{\mu(1,2)}}^* \frac{\bA^{k_{2}}}{z_{i'_{\mu(1,2)}}^{k_{2}+1}} \bb_{i'_{\mu(1,2)}}} \cdots \sqrt{n}\bc_{i_{m_s,s}}^* \frac{\bA^{k_s}}{z_{i_{m_s,s}}^{k_s+1}} \bb_{i_{m_s,s}}\ol{\sqrt{n}\bc_{i'_{\mu(m_s,s)}}^* \frac{\bA^{k_s}}{z_{i'_{\mu(m_s,s)}}^{k_s+1}} \bb_{i'_{\mu(m_s,s)}}} \right] $$
except that we count several times each terms. Indeed, for example, if one wants to rearrange 
\beqy \label{113651}
\Ec{\bc^*_{i_{1,1}} \f{\bA^{k_1}}{z_{i_{1,1}}^{k_1}} \bb_{i_{1,1}}\bc^*_{i_{2,1}} \f{\bA^{k_1}}{z_{i_{2,1}}^{k_1}} \bb_{i_{2,1}}\bc^*_{i_{1,2}} \f{\bA^{k_2}}{z_{i_{1,2}}^{k_2}} \bb_{i_{1,2}} \ol{\bc^*_{i'_{1,1}} \f{\bA^{k_1}}{z_{i'_{1,1}}^{k_1}} \bb_{i'_{1,1}}\bc^*_{i'_{2,1}} \f{\bA^{k_2}}{z_{i'_{2,1}}^{k_2}} \bb_{i'_{2,1}}\bc^*_{i'_{1,2}} \f{\bA^{k_1}}{z_{i'_{1,2}}^{k_1}} \bb_{i'_{1,2}}}},
\eeqy
there are two ways to do it :
\beq
\Ec{\bc^*_{i_{1,1}} \f{\bA^{k_1}}{z_{i_{1,1}}^{k_1}} \bb_{i_{1,1}} \ol{\bc^*_{i'_{1,1}} \f{\bA^{k_1}}{z_{i'_{1,1}}^{k_1}} \bb_{i'_{1,1}}} \bc^*_{i_{2,1}} \f{\bA^{k_1}}{z_{i_{2,1}}^{k_1}} \bb_{i_{2,1}}\ol{\bc^*_{i'_{1,2}} \f{\bA^{k_1}}{z_{i'_{1,2}}^{k_1}} \bb_{i'_{1,2}}} \bc^*_{i_{1,2}} \f{\bA^{k_2}}{z_{i_{1,2}}^{k_2}} \bb_{i_{1,2}} \ol{\bc^*_{i'_{2,1}} \f{\bA^{k_2}}{z_{i'_{2,1}}^{k_2}} \bb_{i'_{2,1}}}},
\eeq
or
\beq
\Ec{\bc^*_{i_{1,1}} \f{\bA^{k_1}}{z_{i_{1,1}}^{k_1}} \bb_{i_{1,1}} \ol{\bc^*_{i'_{1,2}} \f{\bA^{k_1}}{z_{i'_{1,2}}^{k_1}} \bb_{i'_{1,2}}} \bc^*_{i_{2,1}} \f{\bA^{k_1}}{z_{i_{2,1}}^{k_1}} \bb_{i_{2,1}} \ol{\bc^*_{i'_{1,1}} \f{\bA^{k_1}}{z_{i'_{1,1}}^{k_1}} \bb_{i'_{1,1}}} \bc^*_{i_{1,2}} \f{\bA^{k_2}}{z_{i_{1,2}}^{k_2}} \bb_{i_{1,2}} \ol{\bc^*_{i'_{2,1}} \f{\bA^{k_2}}{z_{i'_{2,1}}^{k_2}} \bb_{i'_{2,1}}}},
\eeq
and so \eqref{113651} would be counted twice.
%
Actually, it is easy to see that $\mu_1$ and $\mu_2$ 
give us the same terms if and only if $\sigma = \mu_1 \circ \mu_2^{-1}$ is a permutation such that for all $(i,j) \in 
\mathcal{I}$, $\sigma(i,j) = (i',j)$ (it means that $\sigma$ doesn't change the second index). Let us denote by 
$S_{k_1,\ldots,k_q}$ the  set of such $\si$'s in $S_{\mathcal{I}}$. Then the expression of  \eqref{1131000125} rewrites 
\beq
&& \frac{1}{\card S_{k_1,\ldots,k_q}}\sum_{\mu \in S_{\mathcal{I}}}  \E \left[ \prod_{\al=1}^s\prod_{\beta=1}^{m_\al} \sqrt{n}\bc_{i_{\bet,\al}}^* \frac{\bA^{{k}_\al}}{z_{i_{\bet,\al}}^{{k_\al+1}}} \bb_{i_{\bet,\al}}\ol{\sqrt{n}\bc_{i'_{\mu(\bet,\al)}}^* \frac{\bA^{k_\al}}{z_{i'_{\mu(\bet,\al)}}^{k_\al+1}} \bb_{i'_{\mu(\bet,\al)}}}\right]
\\
&= &\ff{\card S_{k_1,\ldots,k_q}} \sum_{\mu \in S_{\mathcal{I}}}\sum_{\si\in S_{k_1, \ld, k_q}}\prod_{t=1}^s \prod_{u=1}^{m_t}  \ff{z_{i_{u,t}}\ol{z}_{i'_{\mu(u,t)}}}\left(\frac{b^2}{z_{i_{u,t}}\ol{z}_{i'_{\mu(u,t)}}}\right)^{k_t} \bb_{i'_{\si\circ \mu(u,t)}}^*\bb_{i_{u,t}}\bc_{i_{u,t}}^*\bc_{i'_{\si\circ\mu(u,t)}}+
o \left( 1 \right)\\
&=& \ff{\card S_{k_1,\ldots,k_q}}  \sum_{\si\in S_{k_1, \ld, k_q}} \sum_{\mu \in S_{\mathcal{I}}}\prod_{t=1}^s \prod_{u=1}^{m_t} \ff{z_{i_{u,t}}\ol{z}_{i'_{\mu(u,t)}}}\left(\frac{b^2}{z_{i_{u,t}}\ol{z}_{i'_{\mu(u,t)}}}\right)^{k_t} \bb_{i'_{ \mu(u,t)}}^*\bb_{i_{u,t}}\bc_{i_{u,t}}^*\bc_{i'_{\mu(u,t)}}+
o \left( 1 \right)\\
&
=& \sum_{\mu \in S_{\mathcal{I}}}\prod_{t=1}^s \prod_{u=1}^{m_t} \ff{z_{i_{u,t}}\ol{z}_{i'_{\mu(u,t)}}}\left(\frac{b^2}{z_{i_{u,t}}\ol{z}_{i'_{\mu(u,t)}}}\right)^{k_t}\bb_{i'_{ \mu(u,t)}}^*\bb_{i_{u,t}}\bc_{i_{u,t}}^*\bc_{i'_{\mu(u,t)}}+
o \left( 1 \right) 
\eeq

and if we go back to the notation $\{1,\ldots,k_q \}$, we have
$$
\sum_{\mu \in S_{\mathcal{I}}}\prod_{t=1}^s \prod_{u=1}^{m_t} \ff{z_{i_{u,t}}\ol{z}_{i'_{\mu(u,t)}}}\left(\frac{b^2}{z_{i_{u,t}}\ol{z}_{i'_{\mu(u,t)}}}\right)^{k_t} \bb_{i'_{ \mu(u,t)}}^*\bb_{i_{u,t}}\bc_{i_{u,t}}^*\bc_{i'_{\mu(u,t)}}  
\ = \qquad\qquad\qquad\qquad$$
$$\qquad\qquad\qquad\qquad \sum_{\sigma \in S_q}\prod_{t=1}^q  \ff{z_{i_{t}}\ol{z}_{i'_{\si(t)}}}\left(\frac{b^2}{z_{i_{t}}\ol{z}_{i'_{\si(t)}}}\right)^{k_t} \bb_{i'_{\sigma(t)}}^* \bb_{i_t}\cdot \bc_{i_t}^* \bc_{i'_{\sigma(t)}}
$$
and so
\begin{align*}
&\sum_{\displaystyle^{k_1,\ldots,k_q =1}_{k'_1,\ldots,k'_q =1}}^{k_0}  \E \left[ \sqrt{n}\bc_{i_{1}}^* \frac{\bA^{k_1}}{z_{i_1}^{k_1+1}} \bb_{i_{1}} \ol{\sqrt{n}\bc_{i'_{1}}^* \frac{\bA^{k'_1}}{z_{i'_1}^{k'_1+1}} \bb_{i'_{1}}}\cdots \sqrt{n}\bc_{i_{q}}^* \frac{\bA^{{k}_q}}{z_{i_q}^{k_q+1}} \bb_{i_{1,q}}\ol{\sqrt{n}\bc_{i'_{q}}^* \frac{\bA^{k'_q}}{z_{i'_q}^{k'_q+1}} \bb_{i'_{q}}} \right]\\
\\
&=  \!\!\!\!\sum_{k_1,\ldots,k_q =1}^{k_0} \!\!\!\!\!\!\!\!\sum_{\underset{{\left\{k_1,\ldots,k_q \right\}_m = \left\{k'_1,\ldots,k'_q \right\}_m}}{k'_1,\ldots,k'_q =1}}^{k_0}\!\!\!\!\!\!\!\! \E \left[ \sqrt{n}\bc_{i_{1}}^* \frac{\bA^{k_1}}{z_{i_1}^{k_1+1}} \bb_{i_{1}}\ol{\sqrt{n}\bc_{i'_{1}}^* \frac{\bA^{k'_1}}{z_{i'_1}^{k'_1+1}} \bb_{i'_{1}}} \cdots \sqrt{n}\bc_{i_{q}}^* \frac{\bA^{{k}_q}}{z_{i_q}^{k_q+1}} \bb_{i_{1,q}}\ol{\sqrt{n}\bc_{i'_{q}}^* \frac{\bA^{k'_q}}{z_{i'_q}^{k'_q+1}} \bb_{i'_{q}}} \right] + o \left( 1 \right)\\
& =  \sum_{k_1,\ldots,k_q =1}^{k_0} \sum_{\sigma \in S_q}\prod_{t=1}^q  \ff{z_{i_{t}}\ol{z_{i'_{\si(t)}}}}\left(\frac{b^2}{z_{i_{t}}\ol{z_{i'_{\si(t)}}}}\right)^{k_t} \bb_{i'_{\sigma(t)}}^* \bb_{i_t}\cdot \bc_{i_t}^* \bc_{i'_{\sigma(t)}} +o(1) \\
& =    \sum_{\si \in S_q} \prod_{t=1}^q \frac{b^{2}}{z_{i_t} \ol{z_{i'_{\si(t)}}}} \frac{1 - \left(\frac{b^{2}}{z_{i_t} \ol{z_{i'_{\si(t)}}}}\right)^{k_0}}{z_{i_t} \ol{z_{i'_{\si(t)}}}-b^2} \bb^*_{i'_{\si(t)}} \bb_{i_t} \bc^*_{i_t} \bc_{i'_{\si(t)}}   + o\left(1\right).
\end{align*}
This allows to conclude directly.

\subsection{Proof of Lemma \ref{lemmecalculexpectation13813}}\label{preuvelemmecalculexpectation13813}
We want to compute $$E:=\E\Tr \bA\bV\bB\bV^*\bC\bV\bD\bV^*.$$ 
Let us denote the entries of $\bV$ by $v_{ij}$, the entries of $\bA$ by $a_{ij}$, the entries of $\bB$ by $b_{ij}$\ld Then, expanding the trace,  we have 
%
%
  $$E=\sum_{1\le \al,\bet,i,j,\ga,\tau,k,l\le n}\underbrace{\E[ a_{\al\bet}v_{\bet i}b_{ij}\ov_{\ga j}c_{\ga \tau}v_{\tau k}d_{kl}\ov_{\al l}]}_{:=\E_{\al,\bet,i,j,\ga,\tau,k,l}}$$
By the left and right invariance of the Haar measure on the unitary group (see  Proposition \ref{wg}), for  the expectation of a product of entries of $\bV$ and $\ovl{\bV}$ to be non zero, we need each row to appear as much times  in $\bV$ as in $\ovl{\bV}$ and 
each column to appear as much times in $\bV$ as in $\ovl{\bV}$. It follows that 
for $\E_{\al,\bet,i,j,\ga,\tau,k,l}$ to be non zero, we need to have the equalities of multisets:$$\{\bet, \tau\}_m=\{\al,\ga\}_m,\qquad \{i,k\}_m=\{j,l\}_m$$ 
 The first condition is equivalent to one of the three   conditions $$\al=\bet=\ga=\tau\qquad\trm{ or }\qquad \al=\bet\ne\ga=\tau\qquad\trm{ or }\qquad \al=\tau\ne\bet=\ga$$ and the second condition is equivalent to one of the three   conditions $$i=j=k=l\qquad\trm{ or }\qquad i=j\ne k=l\qquad\trm{ or }\qquad i=l\ne j=k.$$
 Hence we have 9 cases to consider below. In each one, the involved moments of the $v_{ij}$'s are computed thanks to e.g. Proposition 4.2.3 of \cite{hiai}: for any $a,b,c,d$, we have \bgt\ite[$\bullet\bullet$\qquad] 
 $\ds \E[|v_{ab}|^4]=\f{2}{n(n+1)}$,
 \ite[$\bullet\bullet$\qquad] $\ds b\ne d\implies  \E [|v_{ab}|^2  |v_{cd}|^2]=\begin{cases}\ff{n(n+1)}&\trm{ if $a=c$}\\ \\
 \ff{n^2-1}&\trm{ if $a\ne c$}
\end{cases}$ 
\ite[$\bullet\bullet$\qquad]  
 $\ds a\ne b \trm{ and } c\ne d\implies  \E [v_{ac}\ovl{v_{ad}} v_{bd} \ovl{v_{bc}}]= 
 -\ff{n(n^2-1)}   $\ent

So let us  treat the 9 cases:  
 
 $\bullet$ Under   condition $\al=\bet=\ga=\tau$ and $i=j=k=l$, we have $$\sum_{\al,\bet,i,j,\ga,\tau,k,l}\E_{\al,\bet,i,j,\ga,\tau,k,l}= \sum_{\al,i}\E[ a_{\al\al}v_{\al i}b_{ii}\ov_{\al i}c_{\al \al}v_{\al i}d_{ii}\ov_{\al i}]=\f{2}{n(n+1)}\sum_{\al}a_{\al\al}c_{\al \al}\sum_ib_{ii}d_{ii}$$
 
 $\bullet$ Under   condition $\al=\bet=\ga=\tau$ and $ i=j\ne k=l$, we get $\ \ds\f{1}{n(n+1)}\sum_{\al}a_{\al\al}c_{\al \al}\sum_{i\ne k}b_{ii}d_{kk}$
 
  $\bullet$ Under   condition $\al=\bet=\ga=\tau$ and $i=l\ne j=k$, we get $\ \ds \ff{n(n+1)}\sum_{\al}a_{\al\al}c_{\al \al}\sum_{i\ne j}b_{ij}d_{ji}$
  
   $\bullet$ Under   condition $ \al=\bet\ne\ga=\tau$ and $i=j=k=l$, we get $\ \ds \ff{n(n+1)}\sum_{\al\ne \ga}a_{\al\al}c_{\ga \ga}\sum_ib_{ii}d_{ii}$
 
 $\bullet$ Under   condition $ \al=\bet\ne\ga=\tau$ and $ i=j\ne k=l$, we get $\ \ds \f{-1}{n(n^2-1)}\sum_{\al\ne \ga}a_{\al\al}c_{\ga \ga}\sum_{i\ne k}b_{ii}d_{kk}$

  $\bullet$ Under   condition $ \al=\bet\ne\ga=\tau$ and $i=l\ne j=k$, we get $\ \ds \f{1}{n^2-1}\sum_{\al\ne\ga}a_{\al\al}c_{\ga \ga}\sum_{i\ne j}b_{ij}d_{ji}$
   
    $\bullet$ Under   condition $\al=\tau\ne\bet=\ga$ and $i=j=k=l$, we get $\ \ds \f{1}{n(n+1)}\sum_{\al\ne\bet}a_{\al\bet}c_{\bet \al}\sum_ib_{ii}d_{ii}$
 
 $\bullet$ Under   condition $\al=\tau\ne\bet=\ga$ and $ i=j\ne k=l$, we get $\ \ds \f{1}{n^2-1}\sum_{\al\ne\bet}a_{\al\bet}c_{\bet \al}\sum_{i\ne k}b_{ii}d_{kk}$
 
  $\bullet$ Under   condition $\al=\tau\ne\bet=\ga$ and $i=l\ne j=k$, we get $\ \ds \f{-1}{n(n^2-1)}\sum_{\al\ne\bet}a_{\al\bet}c_{\bet \al}\sum_{i\ne j}b_{ij}d_{ji}$
   
   Summing up the nine previous sums, we easily get the desired result:
 \beq n(n+1)E&=&2\sum_{\al}a_{\al\al}c_{\al \al}\sum_ib_{ii}d_{ii}+ \sum_{\al}a_{\al\al}c_{\al \al}\sum_{i\ne k}b_{ii}d_{kk}+ \sum_{\al}a_{\al\al}c_{\al \al}\sum_{i\ne j}b_{ij}d_{ji}\\&&+\sum_{\al\ne \ga}a_{\al\al}c_{\ga \ga}\sum_ib_{ii}d_{ii}-\ff{n-1}\sum_{\al\ne \ga}a_{\al\al}c_{\ga \ga}\sum_{i\ne k}b_{ii}d_{kk}+\f{n}{n-1}\sum_{\al\ne\ga}a_{\al\al}c_{\ga \ga}\sum_{i\ne j}b_{ij}d_{ji}\\
   &&+\sum_{\al\ne\bet}a_{\al\bet}c_{\bet \al}\sum_ib_{ii}d_{ii}+\f{n}{n-1}\sum_{\al\ne\bet}a_{\al\bet}c_{\bet \al}\sum_{i\ne k}b_{ii}d_{kk}   -\ff{n-1} \sum_{\al\ne\bet}a_{\al\bet}c_{\bet \al}\sum_{i\ne j}b_{ij}d_{ji}\\
   &=&  \sum_{\al}a_{\al\al}c_{\al \al}\sum_{i, k}b_{ii}d_{kk}+ \sum_{\al}a_{\al\al}c_{\al \al}\sum_{i, j}b_{ij}d_{ji}\\
   &&+\f{n}{n-1}\sum_{\al\ne \ga}a_{\al\al}c_{\ga \ga}\sum_{i,j}b_{ij}d_{ji}-\ff{n-1}\sum_{\al\ne \ga}a_{\al\al}c_{\ga \ga}\sum_{i, k}b_{ii}d_{kk}\\
   && +\f{n}{n-1}\sum_{\al\ne\bet}a_{\al\bet}c_{\bet \al}\sum_{i, k}b_{ii}d_{kk}   -\ff{n-1} \sum_{\al\ne\bet}a_{\al\bet}c_{\bet \al}\sum_{i, j}b_{ij}d_{ji}\\
         &=&  \Tr \bB\Tr \bD\lf\{\sum_{\al}a_{\al\al}c_{\al \al} -\ff{n-1}\sum_{\al\ne \ga}a_{\al\al}c_{\ga \ga}+\f{n}{n-1}\sum_{\al\ne\bet}a_{\al\bet}c_{\bet \al}\ri\}\\ &&+\Tr \bB\bD\lf\{ \sum_{\al}a_{\al\al}c_{\al \al} +\f{n}{n-1}\sum_{\al\ne \ga}a_{\al\al}c_{\ga \ga}   -\ff{n-1} \sum_{\al\ne\bet}a_{\al\bet}c_{\bet \al}\ri\} \\
                &=& \f{n}{n-1}\lf\{\Tr \bA\bC\Tr\bB\Tr\bD+\Tr \bA\Tr\bC\Tr\bB\bD\ri\}\\ &&-\ff{n-1}\lf\{\Tr \bA\bC\Tr \bB \bD+\Tr \bA\Tr\bC\Tr \bB \Tr\bD\ri\}.
             \eeq

  \saut
  
{\bf Acknowledgments:}   We  would like to  thank J. Novak  for   discussions on Weingarten calculus.

 \saut

\end{document}